\documentclass[10pt]{amsart} %
\usepackage{fullpage} %
\usepackage{latexsym}
\usepackage{amsfonts}
\usepackage{amsmath,amsthm,amssymb}
\usepackage{mathabx}
\usepackage{tikz}
\usetikzlibrary{matrix}
\usepackage[colorlinks=true, linkcolor=blue, citecolor=magenta]{hyperref}
\usepackage{mathrsfs}
\usepackage{soul}
\usepackage[all]{xy}

\newcommand{\cal}{\mathcal}

\def\epsilon{\varepsilon}
\def\phi{\varphi}
\def\Pr{\mathbb P}
\def\hat{\widehat}
\def\check{\widecheck}
\def\subset{\subseteq}

\newcommand{\supp}{\mbox{Supp}}
\newcommand{\Curr}{\mbox{Curr}}
\newcommand{\Out}{\mbox{Out}}

\newcommand{\Edges}{\mbox{Edges}}
\newcommand{\bvec}{\overleftarrow}

\newcommand{\Long}{\mbox{Long}}
\newcommand{\Short}{\mbox{Short}}
\newcommand{\Blow}{\mbox{Blow-up}}
\newcommand{\Contr}{\mbox{Contr}}
\newcommand{\disjoint}{\,\overset{.}{\cup}\,}
\newcommand{\mile}{\rm minlength}

\newcommand{\FN}{F_N}   

\newcommand{\CVN}{\mbox{CV}_N}

\newcommand{\PCurr}{\Pr\Curr(\FN)}

\newcommand{\R}{\mathbb R}
\newcommand{\Z}{\mathbb Z}

\newcommand{\N}{\mathbb N}

\newcommand{\A}{{\cal A}}

\newcommand{\Lsig}{L^{\Sigma}}
\newcommand{\Lused}{L^{\tiny \bvec \Gamma}}

\def\strutdepth{\dp\strutbox}
\def \ss{\strut\vadjust{\kern-\strutdepth \sss}}
\def \sss{\vtop to \strutdepth{
\baselineskip\strutdepth\vss\llap{$\diamondsuit\;\;$}\null}}

\def\strutdepth{\dp\strutbox}
\def \sst{\strut\vadjust{\kern-\strutdepth \ssss}}
\def \ssss{\vtop to \strutdepth{
\baselineskip\strutdepth\vss\llap{$\spadesuit\;\;$}\null}}

\def\strutdepth{\dp\strutbox}
\def \ssh{\strut\vadjust{\kern-\strutdepth \sssh}}
\def \sssh{\vtop to \strutdepth{
\baselineskip\strutdepth\vss\llap{$\heartsuit\;\;$}\null}}

\def\qed{\hfill\rlap{$\sqcup$}$\sqcap$\par}
\def\bar{\overline}

\def\strutdepth{\dp\strutbox}
\def \ss{\strut\vadjust{\kern-\strutdepth \sss}}
\def \sss{\vtop to \strutdepth{
\baselineskip\strutdepth\vss\llap{$\diamondsuit\;\;$}\null}}

\def\strutdepth{\dp\strutbox}
\def \sst{\strut\vadjust{\kern-\strutdepth \ssss}}
\def \ssss{\vtop to \strutdepth{
\baselineskip\strutdepth\vss\llap{$\spadesuit\;\;$}\null}}

\def\qed{\hfill\rlap{$\sqcup$}$\sqcap$\par}

\newtheorem{thm}{Theorem}[section]
\newtheorem{cor}[thm]{Corollary}
\newtheorem{lem}[thm]{Lemma}
\newtheorem{prop}[thm]{Proposition}

\newtheorem{fact}[thm]{Fact}
\theoremstyle{definition}
\newtheorem{defn}[thm]{Definition}

\newtheorem{rem}[thm]{Remark}
\newtheorem{defn-rem}[thm]{Definition-Remark}
\newtheorem{convention}[thm]{Convention}

\newtheorem{question}[thm]{Question}
\theoremstyle{remark}

\numberwithin{equation}{section}

\begin{document}

\author[N.~Bedaride]{Nicolas B\'edaride}
\author[A.~Hilion]{Arnaud Hilion}
\author[M.~Lustig]{Martin Lustig}

\address{\tt 
Aix Marseille Univ, CNRS, Centrale Marseille, I2M, Marseille, France.}
\email{\tt Nicolas.Bedaride@univ-amu.fr}
  \email{\tt Arnaud.Hilion@univ-amu.fr}
  \email{\tt Martin.Lustig@univ-amu.fr}

\title[Invariant measures for graph towers]
{Graph towers, laminations and their invariant measures}
\begin{abstract} 
In this paper we present a combinatorial machinery, consisting of a {\em graph tower} $\bvec \Gamma$ and {\em vector towers} $\bvec v$ on $\bvec \Gamma$, which allows us to efficiently describe all invariant measures $\mu = \mu^{\tiny \bvec v}$ on any given shift space over a finite alphabet.

The new technology admits a number of direct applications, in particular concerning invariant measures on non-primitive substitution subshifts, minimal subshifts with many ergodic measures, or an efficient calculation of the measure of a given cylinder. It also applies to currents on a free group $\FN$, and in particular the set of projectively fixed currents under the action of a (possibly reducible) endomorphism $\phi: \FN \to \FN$ is determined, 
when $\phi$ is represented by a train track map.
\end{abstract} 

\subjclass[2010]{Primary 20F65, 37B10, Secondary 37A25}
 
\keywords{Invariant measures on subshifts, currents on free groups, laminations, substitutions}

\maketitle

\section{Introduction}
\label{introduction} 

This paper is situated in the area where geometric group theory and symbolic dynamics overlap:  
in fact, it borrows its language from train track technology (see \cite{CB, DV, Pe, Lu-habil, Th}) and Stallings folds (see \cite{Sta}), 
but the approach chosen here is 
in spirit 
closer to standard objects from symbolic dynamics.

The central object of this paper, a {\em graph tower}, is an infinite sequence 
\begin{equation}
\label{central-object}
\bvec \Gamma 
: \quad
\ldots \overset{f_{3,4}}{\longrightarrow} \Gamma_3  \overset{f_{2,3}}{\longrightarrow}  \Gamma_2 \overset{f_{1,2}}{\longrightarrow} \Gamma_1 \overset{f_{0,1}}{\longrightarrow}  \Gamma_0
\end{equation}
of finite graphs $\Gamma_n$, and of graph maps $f_{n, n+1}$ which map vertices to vertices and for which any composition $f_{m, n} := f_{m, m+1} \circ f_{m+1, m+2} \circ \ldots \circ f_{n-1, n}$ maps edges to non-trivial reduced edge paths. 
Any graph tower $\bvec \Gamma$ defines a {\em used lamination} 
$\Lused$:  this 
is the set of all biinfinite reduced paths $\gamma$ in $\Gamma_0$ which have the property that every finite subpath $\gamma_0$ of $\gamma$ is also the subpath of some $f_{0, n}(e_i)$, where $e_i$ is an edge of the level graph $\Gamma_n$ (we say: $\gamma_0$ is ``used'' by $e_i$).

\smallskip

As a rather special case this set-up includes the possibility that every $\Gamma_n$ is a graph with a single vertex, and the edges $e_i$ are provided with a preferred orientation, such that any of the edge paths $f_{m,n}(e_i)$ crosses only over positively oriented edges. This gives a direct translation (of this special case) of (\ref{central-object}) into a {\em directive sequence}
\begin{equation}
\label{central-object2}
\sigma 
: \quad
\ldots \overset{\sigma_{3,4}}{\longrightarrow} \cal A^*_3  \overset{\sigma_{2,3}}{\longrightarrow}  \cal A^*_2 \overset{\sigma_{1,2}}{\longrightarrow} \cal A^*_1 \overset{\sigma_{0,1}}{\longrightarrow}  \cal A^*_0
\end{equation}
of homomorphisms (``substitutions'')
of free monoids $\cal A^*_n$. A directive sequence $\sigma$ determines (in a similar way as above the used lamination) a subshift $X_\sigma$ on the alphabet $\cal A_0$. Such subshifts have been termed {\em $S$-adic} by Ferenczi 
\cite{Fer96}, and have received since then a lot of attention from combinatorially minded symbolic dynamists (see the survey paper \cite{BD}), in particular if the set $S$ of all the $\sigma_{n, n+1}$ is finite. The used lamination $\Lused$ in the general geometric setting (\ref{central-object}) fits into the more general concept of {\em symbolic laminations}, which are defined for any graph $\Gamma$.

Symbolic laminations constitute a rather universal tool which englobes simultaneously general subshifts on a discrete letter set, as well as algebraic laminations for free groups $\FN$ (see \cite{CHL1-I}), but also geodesic laminations on surfaces  $S_g$ of genus $g \geq 2$ (see \cite{CB, Th}).

\smallskip

An important feature of all graph towers (\ref{central-object}) considered in this paper is that they are {\em expanding}, which means that for any of the $\Gamma_n$ there is a lower bound, for any edge $e_i$ of $\Gamma_n$, to the length $|f_{0, n}(e_i)|$, and this bound tends to $\infty$ for $n \to \infty$. The analogous notion for (\ref{central-object2}) states that the directive sequence $\sigma$ is {\em everywhere growing}. This condition, however, is not restrictive: 

\begin{prop} [see Proposition \ref{any-lamination}]
\label{all-lams}
Every symbolic lamination $\Lsig$ is the {used lamination} of some expanding graph tower $\bvec \Gamma$: 
$$\Lsig = \Lused$$
In particular, every subshift $X$ possesses an everywhere growing $S$-adic expansion $X = X_\sigma$ (with possibly infinite 
set of substitutions $S$).
\end{prop}

The main goal of this paper is to present a strong and easy to handle tool which allows one to present, study and operate invariant measures on any symbolic lamination $\Lsig$ (or on any subshift $X$), through a presentation $\Lsig = \Lused$ of the latter as used lamination of some expanding graph tower $\bvec \Gamma$ (or, respectively, through an $S$-adic expansion $X = X_\sigma$ by means of an everywhere growing directive sequence $\sigma$). 

For this purpose we define a {\em vector tower} $\bvec v = (\vec v_n)_{n \in \N \cup \{0\}}$ on $\bvec \Gamma$ by the condition that the $\vec v_n$ are non-negative vectors, with coordinates indexed by the edges of $\Gamma_n$, which satisfy the {\em compatibility property} $\vec v_m = M(f_{m, n}) \cdot \vec v_n\,$. Here for any graph map $f: \Gamma' \to \Gamma$ we denote by $M(f)$ the {\em transition matrix} which has as coefficients the number of times that the image path $f(e')$ of any edge $e'$ of $\Gamma'$ crosses over the edge $e$ of $\Gamma$ or its inverse $\bar e$ (both counted positively).

We now state 
the main result of this paper (see Theorem \ref{thm1}), in a 
simplified form. 
The reader may also want to consult Figure \ref{fig:tower} for a quick direct insight into the terms used here:

\begin{thm}
\label{main-short}
Let $\bvec \Gamma$ be an expanding graph tower, and let $\Lsig = \Lused$ denote its used lamination.

Then every 
vector tower $\bvec v$ on $\bvec \Gamma$ determines an invariant measure $\mu^{\tiny \bvec v}$ on $\Lsig$, 
and any invariant measure $\mu$ on $\Lsig$ is given via $\mu = \mu^{\tiny \bvec v}$ by some vector tower $\bvec v$ on $\bvec \Gamma$.

Furthermore, the measure of any cylinder can be read off with arbitrary preciseness from a sufficiently large finite part of $\bvec \Gamma$ and $\bvec v$.
\end{thm}

An important special case occurs if a graph tower $\bvec \Gamma$ 
as in (\ref{central-object}) 
is {\em stationary}, i.e. $\bvec \Gamma$ is given by some graph $\Gamma$ and some 
graph self-map $f: \Gamma \to \Gamma$ through the conditions $\Gamma_n = \Gamma$ and $f_{n, n+1} = f$ for all $n \geq 0$. In this  case 
we write $\bvec \Gamma_{\! \! f}$ for $\bvec \Gamma$, and in section \ref{stationary-pseudo} we show that, 
if $f$ is expanding, then the cone $\cal V(\bvec \Gamma_{\! \! f})$ of vector towers on $\bvec \Gamma_{\! \! f}$ is finite dimensional, and its dimension is equal to the number of ergodic probability measures on $\Lsig_f := L^{\tiny \bvec \Gamma_{\! \! f}}$.
The following is a simplified version of Theorem \ref{mieux-que-Moulinette} (see also see Theorem \ref{mieux-que-Moulinette-2} and the preceding two paragraphs):

\begin{thm}
\label{stationary-simple}
Let $\bvec \Gamma_{\! \! f}$ be a stationary graph tower given by some graph map $f: \Gamma \to \Gamma$.
Then there 
is a 1-1 relationship between the ergodic measures on 
the used lamination 
$\Lsig_f$ on one 
hand, 
and the 
extremal non-negative eigenvectors $\vec v$ of $M(f)$ on the other. 
\end{thm}

This bijection extends furthermore canonically (see Corollary \ref{Moulinette-eat-your-heart-out}) to {\em pseudo-stationary} graph towers, where the condition $f_{n, n+1} = f$ is replaced by the weaker assumption $M(f_{n, n+1}) = M(f)$ for all $n \geq 0$.

\medskip

A precise description of our results and the terms used there is given in section \ref{results} below. 
We will thus proceed here rather by explaining some immediate applications of this new ``tower technology'':

\subsection*{Automorphisms and endomorphisms of free groups}
For any stationary graph tower $\bvec \Gamma_{\! \! f}$ it follows that the defining graph map $f: \Gamma \to \Gamma$ has the {train track property} (see Definition-Remark \ref{train-track-need}) which is well-known in the context of free groups $\FN$ and {fully irreducible (= iwip)} automorphisms $\phi \in \Out(\FN)$, see \cite{BH92}. The statement of Theorem \ref{mieux-que-Moulinette} is used 
to describe the fixed point set of {hyperbolic automorphisms} $\phi \in \Out(\FN)$ on projectivized current space $\PCurr$ (see Theorem \ref{thm-BHL-old-new}). This complements nicely earlier results of Uyanik-Lustig \cite{LU} about the generalized North-South dynamics of such $\phi$ on the Teichm\"uller space analogue $\PCurr$.

It should be noted here that the train track map $f: \Gamma \to \Gamma$ used in the above set-up is not assumed to be a homotopy equivalence: Hence the results of Theorem \ref{mieux-que-Moulinette} 
are also valid for 
endomorphisms of free groups,
while in the context of currents many of the standard features and tools known for automorphisms of $\FN$ don't quite apply. 

We believe that the tower technology presented here has the potential for a systematic treatment of the induced action of non-surjective (or even non-injective) endomorphisms of $\FN$ on the corresponding current spaces.

\subsection*{Substitution subshifts}
The restriction of (\ref{central-object}) to (\ref{central-object2}) in the stationary case yields subshifts $X_\sigma$ associated to a single {substitution} $\sigma: \cal A^* \to \cal A^*$. Such {\em substitution subshifts} $X_\sigma$ constitute one of the most cherished classes of subshifts in symbolic dynamics, with many examples that have received a lot of detailed attention 
\cite{D-L, Buf-Sol, F-M-N}. In most cases, however, it is assumed that $\sigma$ is primitive, while in this paper we only need to assume that $\sigma$ is everywhere growing. For such substitutions a similar result as 
Theorem \ref{stationary-simple} above 
has been obtained previously by Bezuglyi-Kwiatkowski-Medynets-Solomyak (see Theorem 2.9 of \cite{BKMS} and Theorem 5.7 of \cite{BD} and its succeeding paragraph). 
However, Theorem \ref{stationary-simple} is more direct and also more general, 
in that it doesn't require as extra assumption the absence of periodic leaves in $X_\sigma$ (see Remarks \ref{BV-translation} and \ref{Arnaud's-chouchou} for details).

\subsection*{Computation of cylinder measures}

It is part of the tower technology presented in this paper (see Remark \ref{sharp-error-estimate})  that for any measure $\mu$, given via $\mu = \mu^{\tiny \bvec v}$ through a vector tower $\bvec v$ on a given graph tower $\bvec \Gamma$, and for any cylinder $C_\gamma$, given through a finite reduced path $\gamma$ in the bottom level graph $\Gamma_0$, the value of $\mu(C_\gamma)$ can be calculated up to an arbitrary small computable error term. When $\bvec \Gamma$ is stationary, $\bvec \Gamma = \bvec \Gamma_{\! \! f}$, then, as described above, $\bvec v$ can be expressed in terms of the non-negative eigenvectors of $M(f)$. This allows a further sharpening of the calculation of $\mu(C_\gamma)$ through the computation of a scalar product formula.
This formula and the fairly direct and efficient calculation of its value for several non-trivial examples are carefully described in our ``cousin paper'' \cite{BHL2}.
As a special (particularly easy) case this includes the situation where $f$ is a train track map that represents an iwip $\phi \in \Out(\FN)$ as above.

\subsection*{Many ergodic measures}
Stationary graph towers are a special case of {\em thin} graph towers, 
i.e. towers $\bvec \Gamma$ where all level graphs have the same number of edges, called the {\em tower dimension} $\dim \bvec \Gamma$. Translated into $S$-adic terminology, this means that all ``level alphabets'' $\cal A_n$ in (\ref{central-object2}) have the same number of letters (so that indeed they can be taken to be equal to a fixed alphabet $\cal A$). For thin 
$\bvec \Gamma$ we prove (see Corollary \ref{finite-ergodic}) that the number of distinct ergodic probability measures on $\Lused$ is bounded above by $\dim \bvec \Gamma$. The analogous statement for the special case of directive sequences over a fixed alphabet $\cal A$ seems to be well-known (see Remark 5 of \cite{BD}). 

In the above sketched context of automorphisms of $\FN$ or geodesic laminations on surfaces, thin graph towers occur naturally (as unfolding paths converging towards a point in the boundary of Outer space or Teichm\"uller space respectively, see section \ref{Outer-space} and \cite{CH, NPR}), but in all cases known the above bound $\dim \bvec \Gamma$ is never attained: in fact, the natural bound for the number of distinct ergodic probability measures seems in each of those cases to be (about) half of the tower dimension. 

Wondering whether for general thin $\bvec \Gamma$ there might be a better
bound than the tower dimension of $\bvec \Gamma$ we realized quickly (see Remark \ref{ergodic-bound}) that for non-minimal laminations $\Lused$ the best bound the number of distinct ergodic probability measures on $\Lused$ is indeed $\dim \bvec \Gamma$. A bit of persistence eventually led us then to the following realization result:

\begin{prop}
\label{existence-of-many}
For any integer $d \geq 1$ there exist directive sequences over an alphabet $\cal A$ with $d$ letters such that the associated subshift $X_\sigma$ is minimal and admits $d$ distinct ergodic probability measures.
\end{prop}

Somehow surprisingly, it turns out that the set $S$ of substitutions in the directive sequence of the last proposition can chosen to be finite (indeed $4$ substitutions suffice, for any value of $d$). The construction used in the proof of Proposition \ref{existence-of-many} has been ``outsourced'' again to our paper \cite{BHL2}; it is exclusively based on the tower technology developed here, plus some calculations that fit onto two pages.

\bigskip

{\it As already indicated above, the methods presented in this paper are by no means disjoint from previously existing ones. We'd like to mention in particular the use of frequency vectors in $S$-adic theory (see \cite{BD} and section \ref{comparison} below), and the work of Bezuglyi, Kwiatkowski, Medynets and Solomyak \cite{BKM, BKMS, BKMS2} on Bratteli-Vershik systems (see 
Appendix 
\ref{Bratteli-Vershik-presentation} for a brief exposition of their work). On the geometry-of-groups side we should mention 
\cite{CH1, CH, Lu-habil, NPR}, but several other authors have also worked on related issues.}

\medskip
\noindent
{\em Organization of the paper:}
In section \ref{set-up-I} we set up the notation and carefully explain the graph language used here, and in section \ref{tools-results} we give the precise statements of our main results. In section \ref{sec:S-adic} we give a brief 
introduction to standard symbolic dynamics terminology (\S \ref{S2-substitutions}) and to the more specialized $S$-adic approach (\S \ref{S-adic-intro}). We then describe the canonical translations from graph tower language to $S$-adic directive sequences 
and back (\S\ref{dictionary}), and for the special $S$-adic case we compare in \S \ref{comparison} our main results with what was known previously. Furthermore, in the Appendix \ref{Bratteli-Vershik-presentation} we give some indications about the relation of our work to an alternative approach via Bratteli diagrams and Vershik maps.

The proof of our main result (Theorem \ref{thm1}) stretches through sections \ref{graph-dialects} - \ref{sec:weight-vectors}:  In \S\ref{graph-dialects} we define three dialects for graphs and graph maps that are used in the sequel, in \S\ref{towers} we define expanding graph towers $\bvec \Gamma$ and their used symbolic lamination $\Lused$. In \S\ref{weights-currents} weight towers $\bvec \omega$ on any $\bvec \Gamma$ are introduced, and in \S\ref{measures-through-weights} it is shown that any such  $\bvec \omega$ describes an invariant measures $\mu^{\tiny \bvec \omega}$ on $\Lused$. In \S \ref{measure-to-weight} it is shown that any invariant measure $\mu$ on $\Lused$ is defined via $\mu = \mu^{\tiny \bvec \omega}$ by some weight tower on $\bvec \Gamma$. Finally in \S \ref{sec:weight-vectors} we introduce vector towers $\bvec v$ on $\bvec \Gamma$ and show that they are equivalent to (but more handy than) weight towers.

(However, the local edges in the ``blow-up dialect'' presentation of graph towers, 
although not any more visible in the vector tower approach, contain through their weights very subtle information, which in many ways encode the key innovative feature of our tower technology.) 

We then proceed in \S \ref{applications} to deduce the first (fairly general) consequences that derive immediately from our new set-up, in particular concerning the special case of graph towers with finite tower dimension (\S \ref{thin-thick-towers}) and the even more specialized but most important case of stationary graph towers (\S \ref{stationary-pseudo}). 
The last subsection, \S \ref{Outer-space}, is devoted entirely to the geometry-of-groups side of our business: although we make an effort in our presentation to give at least some 
ideas 
to the non-experts, in truth this section assumes some familiarity with Outer space and its action by $\Out(\FN)$.

For convenience of our readers we have also added two further 
appendices, 
where the basic facts for convex cones in vector spaces (\S \ref{cones}) and for non-negative matrices (\S \ref{non-negative-matrices}) are recalled, and the terminology used here is specified. 

\medskip
\noindent
{\em Acknowledgements:} During the course of this work we have received many useful comments from our friends and colleagues in and around the marseillian symbolic dynamics community. We'd like to mention in particular Julien Cassaigne, 
Thierry Coulbois, Vincent Delecroix, Fabien Durand, S\'ebastien Ferenczi, Pascal Hubert and Glenn Merlet.

 
\section{Statement of results}\label{results}

In this section we will state precisely the main results of this paper. Before doing so in subsection \ref{tools-results}, we give in section \ref{set-up-I} a brief review of the terms and notations used. None of the content of subsections \ref{S2-graphs} - \ref{ergodic-measures} is new; however, since we use graph language rather than combinatorics (as would be more standard in symbolic dynamics), we redefine below carefully some of these well-known terms.

\subsection{The set-up}
\label{set-up-I}

\subsubsection{Graphs and edge paths}
\label{S2-graphs}

${}^{}$

In this paper a graph $\Gamma$ is a topological space which consists of vertices $v$ or $v_i$ and non-oriented edges $E$ or $E_i$. 
Here an {\em edge} $E$  is always homemorphic to the interval $[0, 1]$, and its two endpoints (and no other point of $E$) belong to the set of {\em vertices} of $\Gamma$. Distinct edges can only intersect in their endpoints.

Since for practical purposes one almost always needs to work with oriented edges, we associate to every non-oriented edge $E$ of $\Gamma$ abstractly two oppositely oriented edges $e$ and $e'$, so that the set $\Edges^\pm(\Gamma)$ of {\em oriented edges of $\Gamma$} contains twice as many elements than (non-oriented) edges present in the topological space $\Gamma$.

For every oriented edge $e \in \Edges^\pm(\Gamma)$ we denote the edge in $\Edges^\pm(\Gamma)$ with reversed orientation by $\bar e$, which gives $\bar{\bar e} = e$. The map $e \mapsto \bar e$ is hence a fixpoint-free involution on the set $\Edges^\pm(\Gamma)$ of oriented edges of $\Gamma$. Whenever need be, we denote by
$$\Edges^+(\Gamma) \subseteq\Edges^\pm(\Gamma)$$
and any section of the quotient map 
$$\Edges^\pm(\Gamma) \to \Edges^\pm(\Gamma) / \langle e = \bar e \rangle\, .$$
One can think of $\Edges^+(\Gamma)$ as the choice of a ``positive orientation'' that is (arbitrarily) defined for any edge of $\Gamma$.
We denote the {\em terminal endpoint} of an oriented edge $e$ by $\tau(e)$. The initial endpoint is given by $\tau(\bar e)$.

\smallskip

Unless otherwise stated, we always assume that a given graph $\Gamma$ is finite 
(i.e. {with} finitely many edges and vertices).
It is also fair to assume that $\Gamma$ is connected, as this is the most interesting case; formally, however, the connectedness of $\Gamma$ is nowhere used in this paper.
The special cases where $\Gamma$ (or one of its connected components) has Euler characteristic $\chi(\Gamma) \geq 0$, i.e. $\Gamma$ is contractible or homotopy equivalent to the circle, are not formally excluded, but they are outside of the interest of this paper and hence not given any attention.

\smallskip

An edge path $\gamma = \ldots e_{i-1} e_i e_{i+1} \ldots$ is a finite, one-sided infinite or biinfinite sequence of edges $e_i \in \Edges^\pm(\Gamma)$ such that $\tau(e_i) = \tau(\bar e_{i+1})$ for {any} {pair of indices $i$, $i+1$}
occurring in $\gamma$. The set of indices $\mathscr I$ of an edge path is always the intersection of $\Z$ with some interval of $\R$.
For finite edge paths (i.e. $\mathscr I$ finite), we work throughout with the convention that the indexing is immaterial: for example, the paths $e_1 e_2 e_3$ and $e_4 e_5 e_6$ are equal if $e_1 = e_4, e_2 = e_5$ and $e_3 = e_6$.
The case of biinfinite edge paths (i.e. $\mathscr I=\Z$) will be specified below in Convention \ref{biinfinite-images}. One-sided infinite edge paths will not play a role in this paper. 

For any finite edge path {$\gamma = e_1 e_2 \ldots e_r$} we denote by $\bar \gamma$ {the path with the {reversed} orientation: $\bar \gamma = \bar e_{r} \bar e_{r-1} \ldots \bar e_{1}$.} 
The {\em terminal endpoint} of {$\gamma$} is defined by $\tau(\gamma) := \tau(e_r)$; {its {\em initial endpoint} is given by $\tau(\bar \gamma) = \tau(\bar e_1)$.}
Given an edge path $\gamma = \ldots e_{i-1} e_i e_{i+1} \ldots$ and some indices {$p\leq q$ occuring on} $\gamma$, the edge path $e_p \ldots e_{q}$ is a {\em subpath} of $\gamma$. An edge path $\gamma = \ldots e_{i-1} e_i e_{i+1} \ldots$ is {\em reduced} if for {any index $i$ one has} $e_{i+1} \neq \bar e_i$. Note that, in general, edge paths are not assumed to be reduced.

\begin{defn}
\label{length-defn}

(1)
For any finite edge path {$\gamma = e_1 e_2 \ldots e_r$} in a graph $\Gamma$ the {\em combinatorial length} (or simply {\em length})
$$|\gamma |= r$$ of $\gamma$ counts the number edges traversed by $\gamma$.

\smallskip
\noindent
(2)
For any second finite edge path $\gamma'$ in $\Gamma$ we denote by
$$|\gamma|_{\gamma'}$$
the number of occurrences of $\gamma'$ as subpath  in the path $\gamma$. In particular, for any edge $e$ of $\Gamma$ we denote by $|\gamma|_e$ the number of times that $\gamma$ crosses over $e$. As a consequence, we obtain:
$$|\gamma| = \sum_{e \,\in\, \text{Edges}^{\pm1}(\Gamma)} |\gamma|_e$$

\smallskip
\noindent
(3)
An edge path $\gamma$ is {\em trivial} if $|\gamma| = 0$. Otherwise $\gamma$ is called {\em non-trivial}. A trivial edge path coincides with its terminal and its initial endpoint.
\end{defn}

\smallskip

\subsubsection{Graph maps}
\label{S2-graph-maps}

${}^{}$

A {\em graph map} $f: \Gamma \to \Gamma'$ is a map between graphs that sends vertices to vertices and edges to possibly non-reduced edge paths, such that $\tau(f(e)) = f(\tau(e))$ and $f(\bar e) = \overline{f(e)}$ for any $e \in \Edges^\pm(\Gamma)$. As usual, the image $f(\gamma)$ of any finite edge path $\gamma = e_1 \ldots e_r$ is the concatenation of paths $f(e_1) \cdot f(e_2) \cdot \ldots \cdot f(e_r)$. It should be noted here that even if $\gamma$ and all $f(e_i)$ are reduced, the edge path $f(\gamma)$ may well be not reduced.

\smallskip

For any graph map $f: \Gamma \to \Gamma'$ there is a well defined {\em transition matrix} (also referred to as {\em incidence matrix}) 
$$M(f) = (m_{e', e})_{e' \,\in\, {\rm Edges}^+(\Gamma'), \,e \,\in\, {\rm Edges}^+(\Gamma)}\, ,$$
where $m_{e', e} := |f(e)|_{e'} + |f(e)|_{\bar e'}$ is equal to the number of times that $f(e)$ crosses over $e'$ or over $\bar e'$. 

Since both of these occurrences are counted positively, the matrix $M(f)$ is always non-negative, 
and it is indeed independent of the choice of the positive orientations given through $\Edges^+(\Gamma) \subset \Edges^\pm(\Gamma)$ and $\Edges^+(\Gamma') \subset \Edges^\pm(\Gamma')$. One easily verifies 
\begin{equation}
\label{matrix-product}
M(g \circ f) = M(g) \cdot M(f)
\end{equation}
for any graph maps $f: \Gamma \to \Gamma'$ and $g: \Gamma' \to \Gamma''$.

\smallskip

We'd like to emphasize again that, even if for every edge $e_i$ of $\Gamma$ and every edge $e'_j$ of $\Gamma'$ the edge paths $f(e_i)$ and $g(e'_j)$ are reduced, it is still quite possible that some of the edge paths $g(f(e_i))$ are not reduced.

\smallskip

\subsubsection{Symbolic dynamics}
\label{S2-symbolic-dynamics}

${}^{}$

For any graph $\Gamma$ we denote by $\Sigma(\Gamma)$ the set of $\Z$-{indexed} biinfinite reduced edge paths (``biinfinite words'') $\gamma = \ldots e_{n-1} e_n e_{n+1} \ldots$ in $\Gamma$. The set $\Sigma(\Gamma)$ is naturally identified with a subset of $({\Edges^\pm(\Gamma)})^\Z$ and 
thus provided 
{(i)} with a product topology, 
{(ii)} with a shift map $S$ {given through}
$$\gamma \mapsto S(\gamma) = \ldots e'_{n-1} e'_n e'_{n+1} \ldots \quad\text{ where } \quad e'_n:= e_{n+1},$$ 
and 
{(iii)}  with an inversion 
$$\gamma \mapsto \bar \gamma = \ldots e'_{n-1} e'_n e'_{n+1} \ldots  \quad\text{ where } \quad e'_n := \bar e_{-n+1}.$$ 

Throughout this paper we use the following conventions:

\begin{convention}
\label{biinfinite-images}
(1)
Two biinfinite paths $\gamma = \ldots e_{-1} e_0 e_1 \ldots$ and $\gamma' = \ldots e'_{-1} e'_0 e'_1 \ldots$ in some graph $\Gamma$ are {\em equal} if $e_k = e'_k$ for every $k \in \Z$. 

\smallskip
\noindent
(2)
Any graph map $f: \Gamma \to \Gamma'$ maps a biinfinite path $\gamma = \ldots e_{-1} e_0 e_1 \ldots$ in $\Gamma$ to a biinfinite path $f(\gamma)  = \ldots e'_{-1} e'_0 e'_1 \ldots$ in $\Gamma'$, where the indexing in $\gamma'$ is determined by the convention that the edge $e'_1$ is the initial edge of the finite path $f(e_1)$.
\end{convention}

A {\em symbolic lamination} on $\Gamma$ is a non-empty subset $\Lsig \subseteq\Sigma(\Gamma)$ which is closed, $S$-invariant
(i.e. $S(\Lsig) = \Lsig$), and invariant under inversion. In symbolic dynamics, symbolic laminations are known under the name of {\em subshift} on the ``alphabet'' $\Edges^\pm(\Gamma)$, if we treat each pair  $e_i$ and $\bar e_i$ as distinct unrelated symbols. In this view, $\Sigma(\Gamma)$ 
can be seen as {\em subshift of finite type} of the {\em full shift} $({\Edges^\pm(\Gamma)})^\Z$. 
For more symbolic dynamics terminology see subsection \ref{S2-substitutions}; in subsection \ref{towers-to-S-adic} details are given about the above ``subshift of finite type'' viewpoint.

A symbolic lamination $\Lsig \subseteq\Sigma(\Gamma)$ is called {\em minimal} if it doesn't contain a proper subset which itself is a symbolic lamination. This is equivalent to requiring that for any path $\gamma \in \Lsig$ the closure of the orbit of $\gamma$ under the action of $S$ and of the inversion is all of $\Lsig$.

\smallskip

\subsubsection{Invariant measures on symbolic laminations}
\label{sec:invariant-measures}

${}^{}$

An {\em invariant measure} $\mu_\Sigma$ for $\Gamma$ is a finite Borel measure on $\Sigma(\Gamma)$ which is invariant under shift and inversion. A concrete combinatorial understanding of such $\mu_\Sigma$ is possible through the following particular (closed and open) measurable sets:

Any reduced path $\gamma = e_1 \ldots e_r$ in $\Gamma$ defines a {\em cylinder} $C_\gamma \subseteq \Sigma(\Gamma)$ which is the set of all biinfinite reduced paths $\ldots e'_{n-1} e'_n e'_{n+1} \ldots$ which satisfy $e' _1 = e_1, \ldots, e'_r = e_r$. We denote its measure by
$$\mu_\Gamma(\gamma) := \mu_\Sigma(C_\gamma) \, .$$
It is a well-known consequence of Carath\'eodory's extension theorem
that $\mu_\Sigma$ is uniquely determined by the values of $\mu_\Gamma(\gamma)$ for all reduced paths $\gamma$ in $\Gamma$.

An invariant measure $\mu_\Sigma$ for some finite graph $\Gamma$ is said to be a {\em probability measure} if $\mu_\Sigma(\Sigma(\Gamma)) = 1$. This is equivalent to requiring that
$$\sum_{e \,\in \text{Edges}^\pm(\Gamma)} \mu_\Gamma(e) = 1$$
(where the edge $e$ is interpreted as edge path of length 1). The support of any invariant measure $\mu_\Sigma$ is a symbolic lamination, 
denoted by $L^\Sigma(\mu_\Sigma)$. It is given by:
\begin{equation}
\label{support-measure}
L^\Sigma(\mu_\Sigma) = \Sigma(\Gamma) \smallsetminus \cup \{C_\gamma \mid \mu_\Gamma(\gamma) = 0\}
\end{equation}

Here we have to assume that $\mu_\Sigma$ is not the zero-measure, as for $\mu_\Sigma = 0$ it follows 
$L^\Sigma(\mu_\Sigma) = \emptyset$.

\smallskip

\subsubsection{Ergodic measures}
\label{ergodic-measures}

${}^{}$

The set $\cal M(\Gamma)$ of invariant measures for $\Gamma$ is naturally equipped with an addition and an external multiplication with scalars $\lambda \in \R_{\geq 0}$. A measure $\mu_\Sigma \in \cal M(\Gamma)$ is 
called {\em ergodic} if any measurable set $X \subset \Sigma(\Gamma)$, which is invariant under shift and inversion, satisfies $\mu_\Sigma(X) = 0$ or $\mu_\Sigma(X) = \mu_\Sigma(\Sigma(\Gamma))$.
(Note that in the context of this paper it is natural to demand that an ``invariant set'' be invariant with respect to both, shift and inversion -- i.e. one passes from $\Z$-actions to $\Z \rtimes \Z/2\Z$-actions.)
This amounts precisely to the condition that any expression of $\mu_\Sigma$ as finite convex combination
$$\mu_\Sigma = \sum \lambda_i \, \mu_\Sigma^i \qquad \text{with all} \qquad \lambda_i > 0, \,\,\mu_\Sigma^i \in \cal M(\Gamma)$$
implies that each $\mu_\Sigma^i$ is a scalar multiple of $\mu_\Sigma$.

The following is well-known (see \cite{Walters}, Theorem 6.10 and the succeeding remark, for a general statement, and \cite{FM} for a more specific symbolic dynamics proof):

\begin{prop}
\label{Fer-Mont}
For any graph $\Gamma$ any family of ergodic measures $\mu_\Sigma^i \in \cal M(\Gamma)$, which are pairwise not scalar multiples of each other, is linearly independent.
\qed
\end{prop}

For more background and general facts about $\cal M(\Gamma)$ see Remark \ref{cone-of-measures}.

\bigskip
\subsection{New tools and results}
\label{tools-results}

${}^{}$

\smallskip

The following is the central tool for this paper:

\begin{defn}
\label{graph-towers-2}
A {\em graph tower} $\bvec \Gamma$ is given by an infinite family $(\Gamma_n)_{n \in \N \cup \{0\}}$ of finite connected {\em level graphs} $\Gamma_n$, and an infinite family $\bvec f = (f_{m, n})_{0 \leq m \leq n}$ of graph maps $f_{m, n}: \Gamma_n \to \Gamma_m$ 
with the following properties:

\begin{enumerate}
\item[(a)]
$f_{m,n}$ maps vertices to vertices.
\item[(b)]
$f_{m,n}$ maps edges to reduced non-trivial edge paths.
\item[(c)]
The family $\bvec f$ is {\em compatible}:  one has $f_{k,m} \circ f_{m,n} = f_{k, n}$ for all integers $n \geq m \geq k \geq 0$. 
\end{enumerate}
For simplicity we will use the abbreviations  $f_n := f_{0, n}$ for all $n \geq 0$.
\end{defn}

To any graph tower there are canonically associated several symbolic laminations. The most useful in our context is 
defined as follows:

\begin{defn}
\label{used-lamination+}
(a)
Let $\bvec \Gamma = ((\Gamma_n)_{n \in \N \cup \{0\}} , (f_{m,n})_{0 \leq m\leq n})$ be a graph tower. A finite edge path $\gamma_0$ in 
the bottom level graph $\Gamma_0$ is called {\em used} if 
there exist arbitrary high levels $n \geq 0$ 
such that $\gamma_0$ is ``used'' by some edge $e$ of $\Gamma_n$, i.e. the path $f_{n}(e)$ 
contains $\gamma_0$ 
or $\bar\gamma_0$ as subpath.

\smallskip
\noindent
(b)
We denote by $\Lused \subset \Sigma(\Gamma_0)$ the symbolic lamination which consists of all biinfinite paths $\gamma$ 
with the property that every finite subpath 
$\gamma_0$ of $\gamma$ is used. The lamination $\Lused$ will be called the {\em used symbolic lamination} of $\bvec \Gamma$.
\end{defn}

Any vertex $v$ of a level graph $\Gamma_m$ 
with $m \geq 1$ of some graph tower $\bvec \Gamma$ is called {\em inessential} if $v$ has valence 2 
and if $f_m$ is locally injective at $v$, and if 
furthermore the same is true for any vertex $v'$ of any level graph $\Gamma_n$ with $n \geq m$ which satisfies $f_{m, n}(v') = v$. It follows directly from the above definitions that erasing all inessential vertices from all level graphs $\Gamma_n$ of level $n \geq 1$ of $\bvec \Gamma$ (while leaving the maps $f_{m, n}$ unchanged) 
defines a new graph tower $\bvec \Gamma^*$ 
for which the used lamination contains that of
$\bvec \Gamma$.

\begin{defn}
\label{strongly-expanding}
(a)
A graph tower $\bvec \Gamma$, given by a family $\bvec f$ of graph maps as in Definition \ref{graph-towers-2}, is said to be {\em strongly expanding} if
$$\lim_{n \to \infty} \left(\min_{e_i \in {\rm Edges}^\pm(\Gamma_n)} |f_{n}(e_i)|\right)\to \infty \, .$$

\smallskip
\noindent
(b)
We say that $\bvec \Gamma$ is {\em expanding} if the associated graph tower $\bvec \Gamma^*$, defined through erasing all inessential vertices as described above, is strongly expanding.
\end{defn}

The expansiveness of a graph tower is crucial to everything we do in this paper; 
in fact, throughout the paper we will only work with expanding graph towers. 
Fortunately this is not really a restriction, as is shown by the following proposition
(proved in 
in section \ref{towers} after Definition-Remark \ref{name-used-lam}), 
which applies indeed to {\em any} symbolic lamination:

\begin{prop}
\label{any-lamination}
Let $\Gamma$ be a finite graph, and let $L^\Sigma \subset \Sigma(\Gamma)$ be an arbitrary symbolic lamination on $\Gamma$. Then there exists an expanding graph tower $\bvec \Gamma = ((\Gamma_n)_{n \in \N \cup \{0\}}, (f_{m,n})_{0\leq m\leq n})$ and a graph map $f: \Gamma \to \Gamma_0$ which induces a 
bijection $$f_{L^\Sigma}:L^\Sigma \to \Lused \, .$$
\end{prop}

\medskip

We now turn to invariant measures on symbolic laminations. The 
central tool of this paper, introduced to study such measures in the context of graph towers, is given by the following:

\begin{defn}
\label{vector-towers}
A {\em vector tower} $\bvec v$ on a given graph tower $\bvec \Gamma$ is a family $\bvec v = (\vec v_n)_{n \in \N \cup\{0\}}$ of non-negative functions 
$$\vec v_n : {\rm Edges}^+(\Gamma_n) \to \R_{\geq 0}$$
on the set of positively oriented edges of the level graphs $\Gamma_n$ of $\bvec \Gamma$.  
Here the functions $\vec v_n$ are thought of as column vectors $\vec v_n = (\vec v_n(e_i))_{e_i \in {\rm Edges}^+(\Gamma_n)}$, 
and they must satisfy the {\em compatibility equalities}
\begin{equation}
\label{defining-compatibility}
\vec v_{m} = M(f_{m,n}) \vec v_{n}
\end{equation}
for all $n \geq m \geq 0$ (where $M(f_{m,n})$ denotes the transition matrix of $f_{m,n}$ as defined in subsection \ref{S2-graph-maps}).
\end{defn}

We notice that there is a natural addition for vector towers on a given graph tower, and similarly an external multiplication with non-negative scalars $\lambda \in \R_{\geq 0}$. We are now able to state the main result of this paper:

\begin{thm}
\label{thm1}
Let $\bvec \Gamma$
be an expanding graph tower with used symbolic lamination $\Lused$, and let $\cal M := \cal M(\Lused)$ denote the set of invariant measures on $\Lused $. Let $\cal V := \cal V(\bvec\Gamma)$ denote the set of vector towers on $\bvec\Gamma$.
\begin{enumerate}
\item
Every vector tower $\bvec v$ on $\bvec \Gamma$ determines an invariant measure $\mu_\Sigma^{\tiny \bvec v}$ on $\Lused $.
\item
Conversely, every invariant measure $\mu_\Sigma$ on $\Lused $ is given via $\mu_\Sigma = \mu_\Sigma^{\tiny\bvec v}$ by some vector tower $\bvec v$ on $\bvec \Gamma$.
\item
The issuing map $\frak m_{\tiny \bvec \Gamma}: \cal V \to \cal M, \bvec v \mapsto \mu_\Sigma^{\tiny\bvec v}$ is linear
(with respect to linear combinations with non-negative scalars).
\item
For any vector tower $\bvec v = (\vec v_n)_{n \in \N \cup\{0\}}$
on $\bvec \Gamma$ the sequence of sums (using the notation of Definition \ref{length-defn} (2))
$$\sum_{e \,\in \text{Edges}^\pm(\Gamma_n)} \vec v_n(e) \cdot |f_n(e)|_\gamma$$
converges, and we obtain, for the cylinder $C_\gamma$:
\begin{equation}
\label{approximation-sum}
\mu_\Sigma^{\tiny\bvec v}(C_\gamma) = \lim_{n \to \infty} \sum_{e \,\in \text{Edges}^\pm(\Gamma_n)} \vec v_n(e) \cdot |f_n(e)|_\gamma
\end{equation}
\end{enumerate}
\end{thm}

We see from property (2) of this theorem that the map $\frak m_{\tiny \bvec \Gamma}: \cal V(\bvec \Gamma) \to \cal M(\Lused)$ is surjective. In \cite{BHL3} the natural class of ``non-repeating'' graph towers will be defined and studied, and it will be shown that for any non-repeating $\bvec \Gamma$ the map $\frak m_{\tiny \bvec \Gamma}$ is also injective.

\medskip

As is shown in property (4) of Theorem \ref{thm1} the relationship between invariant measures and vector towers given by the map $\frak m_{\tiny \bvec \Gamma}$  is not just theoretical, but can also be used very concretely. Indeed, we obtain (see 
Remarks \ref{increasing} and \ref{sharp-error-estimate}):

\begin{rem}
\label{estimate1}
Let $\bvec \Gamma$ and $\bvec v$ be as in Theorem \ref{thm1}. Then we have:
\begin{enumerate}
\item
The sequence of sums on the right hand side of equality (\ref{approximation-sum}) is increasing.
\item
For any $n \geq 1$ the difference between this sum and the limit value $\mu_\Sigma^{\tiny \bvec v}(C_\gamma)$ is bounded above by a constant that can be controlled 
from only knowing the first $n$ levels of the graph tower $\bvec \Gamma$ and the corresponding vectors of the vector tower $\bvec v$. This gives the possibility to calculate the measure of any given cylinder up to an arbitrary close precision.
\item
Even better, using ``weight functions'' as introduced in section \ref{weights-currents} below instead of vector towers, one can actually calculate the precise value of any $\mu_\Sigma^{\tiny \bvec v}(C_\gamma)$. Indeed, once the weight functions are known,
the calculation of $\mu_\Sigma^{\tiny\bvec v}(C_\gamma)$ becomes  rather practical and efficient. This is further developed in
\cite{BHL2}, where several concrete examples are studied.
\end{enumerate}
\end{rem}

In view of Proposition \ref{Fer-Mont} 
we note that statements (2) and (3) of Theorem \ref{thm1} give directly the following corollary
(compare also \cite{Dur}, where a similar result is stated in the context of $S$-adic subshifts,
as well as \cite{BKMS2}, Proposition 2.13, for Bratteli-Vershik systems).

\begin{cor}
\label{finite-ergodic}
Let $\Gamma$ be a finite graph, and let $L^\Sigma \subset \Sigma(\Gamma)$ be an arbitrary symbolic lamination on $\Gamma$. Let $\bvec \Gamma = ((\Gamma_n)_{n \in \N \cup \{0\}}, (f_{m,n})_{0\leq m\leq n})$ be an expanding graph tower 
with $0$-th level graph 
$\Gamma_0 = \Gamma$ and with 
$\Lused = L^\Sigma$.

If there exists an integer $K \geq 1$ such that for infinitely many distinct level graphs $\Gamma_n$ the number of edges in $\text{Edges}^+(\Gamma_n)$ is bounded above by $K$, then the number of distinct ergodic probability measures on $L^\Sigma$ is also bounded above by $K$.
\qed
\end{cor}

Explicit examples 
for arbitrary $K \geq 2$, where $\Lsig = \Lused$ is minimal and admits $K$ distinct ergodic probability measures, and where every level graph $\Gamma_n$ of $\bvec \Gamma$ has precisely $K$ edges, 
have been exhibited 
in the recent paper \cite{BHL2} by the authors.

\medskip

A graph tower $\bvec \Gamma$ is called {\em stationary} if there exists a graph map $f: \Gamma \to \Gamma$, such that for every level graph $\Gamma_n$ of $\bvec \Gamma$ an identification 
$\Gamma_n = \Gamma$ 
exists, 
modulo which we have $f_{n-1, n} = f$ for all $n \geq 1$. More generally, we say that $\bvec \Gamma$ is {\em pseudo-stationary} if for any $n \geq 1$ the transition matrix $M(f_{n-1, n})$ is identical to some fixed {\em level transition matrix} $M$, up to suitable permutations of the rows and columns. Thus every stationary $\bvec \Gamma$ is also pseudo-stationary. It follows directly from equality (\ref{defining-compatibility}) that every non-negative eigenvector $\vec v$ of $M$ with eigenvalue $\lambda > 1$ defines canonically an associated vector tower $\bvec v(\vec v) := (\frac{1}{\lambda^n}\vec v)_{n \in \N \cup \{0\}}$ on the given pseudo-stationary $\bvec \Gamma$.

The non-negative eigenvectors of a (possibly reducible) non-negative matrix $M$ have been studied and used in many different contexts; a summary of the definitions and results needed here is given in Appendix \ref{non-negative-matrices}, where 
in particular an explicit description of the finitely many non-negative {\em principal eigenvectors} $\vec v_i$ of any such $M$ 
is given in Proposition \ref{principal-ev} and Definition-Remark \ref{barycentric-principal}. We prove (see Theorem \ref{mieux-que-Moulinette} and Corollary \ref{Moulinette-eat-your-heart-out}):

\begin{thm}
\label{mieux-que-Moulinette-2}
Let $\bvec \Gamma$ be an expanding pseudo-stationary graph tower with level transition matrix $M$. Then every non-negative eigenvector $\vec v$ of $M$ determines an invariant measure $\mu_\Sigma^{\tiny \bvec v(\vec v)}$ on $\Lused$. The measure $\mu_\Sigma^{\tiny \bvec v(\vec v)}$ is ergodic if and only if $\vec v$ is a multiple of one of the principal eigenvectors $\vec v_i$ of $M$.

Conversely, for $\mu_\Sigma^i := \mu_\Sigma^{\tiny \bvec v(\vec v_i)}$,
any invariant measure $\mu_\Sigma$ on $\Lused$ is a non-negative linear combination
$$\mu_\Sigma = \sum c_i
\mu_\Sigma^i\, ,$$
where the coefficients $c_i 
\geq 0$ are uniquely determined by $\mu_\Sigma$. 
In particular, the $\mu_\Sigma^i$ are up to scalar multiples the only ergodic measures on $\Lused$.
\end{thm}

This generalizes and slightly improves a result of Bezuglyi, Kwiatkowski, Medynets and Solomyak (see \cite{BKMS}); for a precise account of the relation between their result and ours see Remark \ref{Arnaud's-chouchou}.

\medskip
In section \ref{Outer-space} we apply the above Theorem \ref{mieux-que-Moulinette-2} to automorphisms of a free group $\FN$ and to their action on the space $\PCurr$ of projectivized currents on $\FN$. In particular, in Theorem \ref{thm-BHL-old-new} a precise description of the ``attracting'' fixed points of $\phi$ in $\PCurr$ in terms of the principal eigenvectors of the transition matrix $M(f)$ is given, for any $\phi \in \Out(\FN)$ which is given by a (not necessarily irreducible) train track map $f$.


\section{$S$-adic symbolic dynamics}
\label{sec:S-adic}

The goal of this section is to give a short introduction to $S$-adic theory and present in a concise manner 
the known results that come close to the results of this paper (see subsection \ref{comparison}). 
We will give a precise prescription how to translate our results into $S$-adic language
(see subsection \ref{dictionary}), 
and it will become apparent that all of the above mentioned $S$-adic results can be derived 
directly as easy consequences from the main result of this paper (Theorem \ref{thm1}). The converse procedure, even for the special case of $S$-adic subshifts, seems rather far fetched.

We also give in Appendix \ref{Bratteli-Vershik-presentation} a brief review of Bratteli-Vershik technology in the context of the $S$-adic approach, and again compare there its merits and the issuing results to the ones presented in this paper.

\smallskip

Throughout this section we use \cite{BD} as standard reference; indeed, we try to use as much as possible their notations.

\subsection{Classical symbolic dynamics and substitutions}
\label{S2-substitutions}

${}^{}$

\smallskip

Let $\cal A = \{a_1, \ldots, a_N\}$ be a finite set, called {\em alphabet}. We denote by $\cal A^*$ the free monoid over $\cal A$. Its neutral element, the empty word, is denoted in this paper  by $1_\cal A$ (following the tradition in group theory rather than that in combinatorics). 

Furthermore, let
$$\Sigma_\cal A = \{\ldots x_{-1} x_0 x_1 x_2 \ldots \mid x_i \in \cal A \}$$
be the set of biinfinite words in $\cal A$, called the {\em full shift} over $\cal A$. 

For any two ``words'' $v = v_1 \ldots v_r$ and $w = w_1 \ldots w_s$ in $\cal A^*$ we define the {\em cylinder} 
$$[v, w] \subseteq\Sigma_\cal A$$
as the set of all biinfinite words $\ldots x_{-1} x_0 x_1 x_2 \ldots$ in $\cal A$ which satisfy $x_{-r+1} = v_1,  x_{-r+2} = v_2, \ldots,  x_{0} = v_r$ and $x_1 = w_1, \ldots, x_s = w_s$. The full shift $\Sigma_\cal A$, being in bijection with the set $\cal A^\Z$, is naturally equipped with the product topology, where $\cal A$ is given the discrete topology. The set of cylinders $[v,w]$, for $v,w \in \cal A^*$, form a basis of this topology. The full shift $\Sigma_\cal A$ is compact, and indeed it is a Cantor set.

The shift map $S:\Sigma_\cal A \to \Sigma_\cal A$ is defined for $x = \ldots x_{-1} x_0 x_1 x_2 \ldots$ by $S(x) = \ldots y_{-1} y_0 y_1 y_2 \ldots$, with $y_n = x_{n+1}$ for all $n \in \Z$. It is bijective and continuous with respect to the above product topology, and hence a homeomorphism. 

A {\em subshift} is a non-empty closed subset $X$ of $\Sigma_\cal A$ which is invariant under the shift map $S$. Such a subshift $X$ is called {\em minimal} if it is the closure of the shift-orbit of any $x \in X$.

Let $\mu$ be a finite Borel measure supported on a subshift $X\subseteq\Sigma_\cal A$. The measure is called {\em invariant} if for every measurable set $A\subseteq X$ one has $\mu(S^{-1}(A))=\mu(A)$.  Such a measure $\mu$ is said to be {\em ergodic} if $\mu$ can not be written in any non-trivial way as sum $\mu_1 + \mu_2$ of two invariant measures $\mu_1$ and $\mu_2$ (i.e. $\mu_1 \neq 0 \neq \mu_2$ and $\mu_1 \neq \lambda \mu_2$ for any $\lambda \in \R_{> 0}$). An invariant measure is called a {\em probability measure} if $\mu(X) = 1$, which is equivalent to $\underset{a_i \in \cal A}{\sum} \mu([1_\cal A, a_i]) = 1$.

Furthermore, from Carath\'eodory's extension theorem it is well-known (see also \cite{FM}) that any invariant measure $\mu$ is determined by its values $\mu([1_\cal A, w])$ on the cylinders $[1_\cal A, w]$ for all $w \in \cal A^*$.

\begin{defn-rem}
\label{substitution}
A {\em substitution} $\sigma$ is given by a map
$$\cal A \to \cal A^*, \,\, a_i \mapsto \sigma(a_i) \, .$$
A substitution defines both, an endomorphism of $\cal A^*$, and a continuous map from $\Sigma_\cal A$ to itself which maps $[v, w]$ to $[\sigma(v), \sigma(w)]$.  Both of these maps are also denoted by $\sigma$, and both are summarized under the name of ``substitution''.
\end{defn-rem}

A substitution $\sigma: \cal A^* \to \cal A^*$ is called {\em non-erasing} if $\sigma(a_i) \neq 1_\cal A$ for any $a_i \in \cal A$. One says that $\sigma$ is {\em everywhere growing} if each $a_i \in \cal A$ satisfies $|\sigma^n(a_i)| \to \infty$ for $n \to \infty$. 

Here 
$|w |$ denotes the word length of $w \in \cal A^*$, and for any 
further word $v \in \cal A^*$ we denote by $|w|_v$ the number of occurrences of $v$ as factor in $w$ (where ``factor'' is here meant synonymously to what is called ``subword'' in combinatorial group theory).

For any substitution $\sigma$ we define the associated language $\mathcal L_\sigma \subseteq \cal A^*$ to be the set of factors in $\cal A^*$ of the words $\sigma^n(a_i)$, with $n\geq 1$ and $a_i\in\mathcal A$.

One defines the subshift $X_\sigma \subseteq \Sigma_\cal A$ {\em associated to the substitution $\sigma$} as the set of all $x = \ldots x_{k-1} x_k x_{k+1} \dots\in \Sigma_\cal A$ 
with the property that for any integers $m \geq n$ the word $x_{n} \ldots x_m$ is an element of $\mathcal L_\sigma$.

For any substitution $\sigma: \cal A^* \to \cal A^*$ let $m_{i, j} := |\sigma(a_j)|_{a_i}$ be the number of occurrences of the letter $a_i$ in the word $\sigma(a_j)$. The non-negative matrix 
$$M_\sigma := (m_{i, j})_{a_i, a_j \in \cal A}$$
is called the {\em incidence matrix} for the substitution $\sigma$. The substitution $\sigma$ is called {\em primitive} if $M_\sigma$ is primitive, i.e. there exists an integer $k$ such that every coefficient of $M_\sigma^k$ is positive.
See Appendix \ref{non-negative-matrices} for more terminology and basic facts about non-negative matrices.

\subsection{$S$-adic sequences}
\label{S-adic-intro}

${}^{}$

In $S$-adic theory (see for instance \cite{BD, DLR}) one considers {\em directive sequences} of free monoids 
$\cal A^*_n$ and of monoid morphisms $\sigma_n: \cal A^*_{n+1} \to \cal A^*_{n}$ (for $n \geq 0$). The morphisms $\sigma_n$ 
(also called ``substitutions'') belong to a given set $S$, which in many circumstances is assumed to be finite. 

To any such a directive sequence $\sigma: = \sigma_0 \circ \sigma_1 \circ \ldots$ one associates the language $\cal L_\sigma \subset \cal A_0^*$, defined as the set of factors in $\cal A_0^*$ of the words $\sigma_0 \circ \sigma_1 \circ \ldots \circ \sigma_n(a_i)$, for any $n\geq 0$ and any $a_i\in\mathcal A_{n+1}$.

The subshift $X_\sigma \subseteq \Sigma_{\cal A_0}$ {\em associated to the directive sequence $\sigma$} is the set of all $x = \ldots x_{k-1} x_k x_{k+1} \dots\in \Sigma_{\cal A_0}$ such that for any two integers $m \geq n$ the word $x_{n} \ldots x_m$ is an element of $\mathcal L_\sigma$.

The directive sequence $\sigma$ is called {\em primitive} (or {\em weakly primitive} by some authors) if for any $m \geq 1$ there is an integer $n \geq m+1$ such the incidence matrix $M_{\sigma_{[m, n)}}$ is positive, for $\sigma_{[m, n)} := \sigma_m \circ \sigma_{m+1} \circ \ldots \circ \sigma_{n-1}$.

One says that $\sigma$ is {\em everywhere growing} if one has
$$\min_{a_i \in \cal A_n}|\sigma_{[0, n)}(a_i)| \to \infty \qquad \text{for} \qquad n \to \infty \, .$$

The above terminology coincides with that for substitution subshifts as introduced in subsection \ref{S2-substitutions}: indeed, one recovers 
the latter as special case of a stationary $S$-adic sequence, i.e. all terms $\sigma_n$ in the directive sequence $\sigma$ are equal.

\begin{prop}[Theorem 5.2 of \cite{BD}]
For any primitive everywhere growing directive sequence $\sigma$ the subshift $X_\sigma$ is minimal.
\end{prop}

If the directive sequence $\sigma$ in the last proposition is stationary (or ``strongly minimal'', see Definition 5.1 of \cite{BD}), then one can deduce furthermore that $X_\sigma$ is uniquely ergodic. In general, however, there are famous counterexamples to such a ``tempting'' conclusion (see section \ref{applications}).

\smallskip

\subsection{Dictionary}
\label{dictionary}

\subsubsection{From $S$-adic to graph towers}
\label{S-adic-to-towers}

${}^{}$

There is a canonical way how to translate the $S$-adic setting into that of graph towers as defined in section~\ref{tools-results}: 

\smallskip

We first consider any finite alphabet $\cal A$ and associate to it a {\em rose with $N$ leaves}, i.e. a graph $R(\cal A)$ with one vertex and $N =\# \cal A$ edges. We fix a bijection $a_i \mapsto e(a_i)$ between $\cal A$ and the set $\Edges^+(R(\cal A))$ of positively oriented edges. 

Any subshift $X \subset \Sigma_\cal A$ gives thus canonically rise to a symbolic lamination $\Lsig(X) \subset \Sigma(R(\cal A))$ which 
consists precisely of those biinfinite paths $\ldots e_i e_{i+1} \ldots$ in $R(\cal A)$ which have the property that for some $\ldots x_i x_{i+1}\ldots \in X$ one has $e_i = e(x_i)$ for all $i \in \Z$, or else $e_i = \bar e(x_{-i+1})$ for all $i \in \Z$. This defines a map $\theta: \Lsig(X) \to X$ which is 2 - 1: for any $x \in X$ the preimage set $\theta^{-1}(x)$ consists precisely of two biinfinite paths $\gamma$ and $\bar \gamma$.

Furthermore, an invariant measure $\mu$ on $X$ defines canonically an invariant measure $\mu_\Sigma$ on $\Lsig(X)$, which is determined by setting for any finite path $\gamma = e_1 \ldots e_q$ in $R(\cal A)$, with $e_i = e(x_i)$
and $x_i \in \cal A$, the following cylinder measures:
\begin{equation}
\label{cylinder-measures}
\mu_\Sigma(C_\gamma) = \mu_\Sigma(C_{\bar \gamma}) := \frac{1}{2} \mu([1_\cal A, x_1 \ldots x_q])
\end{equation}
The factor $\frac{1}{2}$ is introduced only to ensure that any probability measure $\mu$ gives rise to a probability measure $\mu_\Sigma$. Indeed, we have as direct consequence of the above definition, for any measurable subset $X' \subset X$, the equality
$$\mu_\Sigma(\theta^{-1}(X')) = \mu(X') \, .$$

\smallskip

Furthermore, for any two free monoids $\cal A^*_1$ and $\cal A^*_2$ each morphism $\sigma: \cal A^*_{1} \to \cal A^*_2$ is realized by a graph map $f_\sigma: R(\cal A_1) \to R(\cal A_2)$ which maps the edge $e(a_i)$ to the edge path $e(b_1) e(b_2) \ldots e(b_r)$, where $a_i \in \cal A_1$ and $\sigma(a_i) = b_1 b_2 \ldots b_r \in \cal A^*_2$.

\begin{rem}
\label{directive-seq-to-tower}
We now consider a directive sequence
$\sigma = \sigma_0 \circ \sigma_1 \circ \ldots$ as specified in subsection \ref{S-adic-intro}. 

\smallskip
\noindent
(1)
For any alphabet $\cal A_n$ of this sequence let $R_n := R(\cal A_n)$ be the rose as defined above, and for any morphism $\sigma_n: \cal A^*_{n+1} \to \cal A^*_n$ let $f_{n, n+1}:= f_{\sigma_n}: R_{n+1} \to R_n$ be the corresponding graph map.

\smallskip
\noindent
(2)
With this set-up we observe directly that the families $(R_n)_{\N \cup \{0\}}$ and $(f_m \circ \ldots \circ f_n)_{0\leq m\leq n}$ define a graph tower $\bvec R(\sigma)$, and that $\bvec R(\sigma)$ is expanding if and only if the directive sequence $\sigma$ is everywhere growing.

\smallskip
\noindent
(3)
The used lamination $L^{\tiny \bvec R(\sigma)}$
of $\bvec R(\sigma)$ is now directly related to the subshift $X_\sigma$ defined by the directive sequence $\sigma$. Indeed, one verifies easily from the above set-up:
$$L^{\tiny \bvec R(\sigma)}= \Lsig(X_\sigma)$$

\smallskip
\noindent
(4)
In particular, every invariant measure $\mu$ on $X_\sigma$ defines canonically an invariant measure $\mu_\Sigma$ on 
$L^{\tiny \bvec R(\sigma)}$ which satisfies (\ref{cylinder-measures}), for $\cal A_0$ in place of $\cal A$.
\end{rem}


\subsubsection{From graph towers to  $S$-adic}
\label{towers-to-S-adic}

${}^{}$

There is also a canonical way how to translate the setting of graph towers and symbolic laminations presented in section \ref{tools-results} into the traditional $S$-adic language from symbolic dynamics:

To any graph $\Gamma$ we associate an alphabet $\cal A_\Gamma$ which is in bijection with the set $\Edges^\pm(\Gamma)$, i.e. any (non-oriented) topological edge $E$ of $\Gamma$, which gives rise to the pair of oriented edges $e, \bar e \in \Edges^\pm(\Gamma)$, is represented by two symbols $a_e, a_{\bar e} \in \cal A_\Gamma$. Thus $\cal A_\Gamma$ is equipped with a natural fixed-point-free involution $\iota: \cal A_\Gamma \to \cal A_\Gamma, \,\, a_e \mapsto a_{\bar e}$. This map extends to an involution 
$\cal A^*_\Gamma \to \cal A^*_\Gamma$, still denoted by $\iota$, which is given by
$$\iota : x_1 x_2 \ldots x_q \mapsto \iota(x_q) \iota(x_{q-1}) \ldots \iota(x_1)$$
where $x_i \in \cal A_\Gamma$ for all $i = 1, \ldots, q$.

One then considers the subshift of finite type 
$\Sigma^\iota_{\cal A_\Gamma}\subset \Sigma_{\cal A_\Gamma}$ defined through the rule that no factor of type $a_e a_{\bar e}$ is admitted in the biinfinite words of $\Sigma^\iota_{\cal A_\Gamma}$. We now work with subshifts $X \subset \Sigma^\iota_{\cal A_\Gamma}$ which are invariant under $\iota$, and similarly with $\iota_*$-invariant measures $\mu$ on $X$, i.e. one has 
\begin{equation}
\label{iota-invariant}
\text{(a)} \quad \iota(X) = X \qquad \text{and} \qquad \text{(b)} \quad \mu(\iota(X')) = \mu(X')
\end{equation}
for any measurable subset $X'$ of $X$. The following is fairly straightforward from these definitions; details have been given in \cite{ABHS}, Annex A, 
and in \cite{LU}, section 5.2 :

\begin{prop}
\label{back-and-forth}
Let $\Gamma$ be any finite graph, and let $\A_\Gamma$ and $\iota$ be as defined above.

\smallskip
\noindent
(1)
Any symbolic lamination $\Lsig \subset \Sigma(\Gamma)$ gives canonically rise to a subshift $X \subset \Sigma^\iota_{\cal A_\Gamma}$ which satisfies 
(\ref{iota-invariant})(a), where any finite reduced path $\gamma = e_1 \ldots e_q$ in $\Gamma$ is a subpath of some biinfinite path of $\Lsig$ if and only if $a_{e_1} \ldots a_{e_q}$ (and hence also $a_{\bar e_q} \ldots a_{\bar e_1}$) 
is a factor of some biinfinite word in $X$.

\smallskip
\noindent
(2)
Similarly, to any invariant measure $\mu_\Sigma$ on $\Sigma(\Gamma)$ there is canonically associated a measure $\mu$ on $\Sigma^\iota_{\cal A_\Gamma}$ which satisfies (\ref{iota-invariant})(b), and furthermore 
$$
\mu_\Sigma(C_\gamma) = 
\mu([1_{\cal A_\Gamma}, a_{e_1} \ldots a_{e_q}])$$
for any finite reduced path $\gamma = e_1 \ldots e_q$ in $\Gamma$.
\qed
\end{prop}

It follows directly from this set-up that any graph map $f: \Gamma \to \Gamma'$ gives rise to a monoid morphism (i.e. a ``substitution'' in the slightly generalized sense of subsection \ref{S-adic-intro}) $f^{\cal A}: \cal A_\Gamma \to \cal A_{\Gamma'}$.  Hence any family of graphs $\Gamma_n$ and of graph maps $f_{m, n}: \Gamma_n \to \Gamma_m$, with $n \geq m \geq 0$, which defines a graph tower $\bvec \Gamma$, gives directly rise to a family of substitutions $\sigma_{m,n} := f^{\cal A}_{m, n}: \cal A_{\Gamma_n} \to \cal A_{\Gamma_m}$ and hence to a directive sequence 
$\sigma := \sigma_{0,1} \circ \sigma_{1,2}\circ \ldots\,\,$.
The latter is everywhere growing if and only if the graph tower $\bvec \Gamma$ is expanding.
For the subshift $X_\sigma \subset \Sigma_{\cal A_{\Gamma_0}}$ which is associated to this directive sequence (see subsection \ref{S-adic-intro}) one easily verifies:

\begin{lem}
\label{same-subshift}
The subshift $X_\sigma$ is precisely the subshift associated in Proposition \ref{back-and-forth} (1) to the used lamination $\Lused$.
\qed
\end{lem}

Formally, this gives the possibility to avoid graph maps and graph towers altogether, and do everything presented in this paper in the traditional settings of symbolic dynamics, with the extra feature that one always has to keep track 
of the involution $\iota$ that needs to be carried along. 
However, in this paper we prefer the approach through graph towers and vector towers, the reason being that this setting is much more suitable in the context of automorphisms of free groups $\FN$ and closer to the existing tools for the study of their action on current space $\PCurr$.

\begin{rem}
\label{minimality-criterion}
The reader who prefers traditional symbolic dynamics, and hence decides anyway to engage into the alternative to ``carry along $\iota$'', should be immediately given a warning that this leads to certain unpleasant phenomena, for instance:

If $X$ is the subshift associated in Propostiton \ref{back-and-forth} (1) to a given symbolic lamination $\Lsig$, then assuming that $\Lsig$ is minimal doesn't imply that $X$ is minimal. Indeed, it is quite possible that $X$ is the disjoint union of two minimal sub-subshifts $X_1$ an $X_2$ which are vice versa $\iota$-images of each other. This happens in particular if $\Lsig$ itself is obtained via $\Lsig = L^{\tiny \bvec R(\sigma)}$ from some directive sequence $\sigma$ as described in Remark \ref{directive-seq-to-tower}.
\end{rem}

\subsection{Comparison of results}
\label{comparison}

${}^{}$

The following general result for arbitrary subshifts $X \subset \Sigma_\cal A$ (with 
finite $\cal A$) is well-known:

If $\mu$ is an ergodic measure on $X$, then for any word $w \in \cal A^*$ and any $\mu$-generic biinfinite word $x = \ldots x_i x_{i+1} \ldots \in X$ the {\em frequency} 
$$\cal F_w(x) := \lim_{n \to \infty} \frac{1}{2n+1}|x_{-n} \ldots x_n|_w$$
is well defined, where as before 
for any $v \in \cal A^*$ we denote by $|v|_w$ the number of occurrences 
of $w$ as factor (= subword) of $v$. One obtains as a standard consequence of Birkhoff's ergodicity theorem:

\begin{fact}
\label{frequencies-measure}
Let $\mu$ be an ergodic measure on a subshift $X \subset \Sigma_\cal A$.
For every $w \in \cal A^*$ and any $\mu$-generic $x \in X$ one has:
$$\cal F_w(x) = \mu([1_\cal A, w])$$
\end{fact}

Turning now to $S$-adic subshifts, say with directive sequence $\sigma = \sigma_0 \circ \sigma_1 \circ \ldots$ as above, we abbreviate the composition of transition matrices by 
$$M_{[m,n)} := M(\sigma_{[m, n)}) = M(\sigma_m \circ \cdots \circ \sigma_{n-1}) = M(\sigma_m) \cdot \ldots \cdot M(\sigma_{n-1})$$ 
for any non-negative integers $m \leq n-1$. One considers the non-negative  cone $\R_{\geq 0}^{d_n}$, for $d_n = \# \cal A_n$, and the sequence of image cones $\cal C_n^m := M_{[m,n)}(\R_{\geq 0}^{d_n})$ in $\R_{\geq 0}^{d_m}$. Then one has:

\begin{prop}[\cite{BD}, Theorem 5.7]
\label{cone-intersection}
The subshift $X_\sigma$ is uniquely ergodic if 
$\sigma$ is everywhere growing and if for every $m \geq 0$ the nested intersection $\cal C^m := \bigcap\{\cal C_n^m \mid n \geq m+1 \}$ is $1$-dimensional.
\end{prop}
A coarser but more general criterion is given by the following:

\begin{prop}[\cite{DLR}, \S4]
\label{ergodic-bound}
Let $\sigma$ be an everywhere growing directive sequence of substitutions $\sigma_n: \cal A^*_{n+1} \to \cal A^*_n$. Then the number of ergodic probability measures on the $S$-adic subshift $X_\sigma$ is bounded above by 
$$\max_{n \in \N} \#\cal A_n\, .$$
\end{prop}

More detailed information is given by the following: 

\begin{prop}[\cite{BD}, Theorem 5.7]
\label{image-vectors}
Let $\sigma$ be an everywhere growing directive sequence of substitutions $\sigma_n: \cal A^*_{n+1} \to \cal A^*_n$. Then the set of vectors 
$\vec v  \in \cal C^0 = \bigcap\{\cal C_n^0 \mid n \geq 1 \}$ coincides with the set of points $(\mu([1_\cal A, a_i])_{a_i \in \cal A_0}$,  
where $\mu$ is any invariant measure on $X_\sigma$.
\end{prop}

Fact \ref{frequencies-measure} and all 
of the last three propositions are direct consequences of our Theorem \ref{thm1} via the dictionary given in subsection \ref{dictionary}. 
In the converse direction, Proposition \ref{image-vectors}, when applied not just to the cone $\cal C^0$ but also to the cones $\cal C^m$ with arbitrary $m \geq 0$, seems to give the possibility to associate (in a non-unique way) to any directive sequence $\sigma$ of substitutions and any invariant measure $\mu$ on the subshift $X_\sigma$ a vector tower $\bvec v(\mu)$ on the graph tower $\bvec R(\sigma)$, which determines the values of $\mu$ on all cylinders that are determined by the letters of any of the $\cal A_n$ of $\sigma$.
This would constitute a weak analogue of part (2) of Theorem \ref{thm1} in the special case of $S$-adic subshifts.
The approach pursued in Theorem 5.7 of \cite{BD}, however, works entirely through frequencies and is hence in practice limited to ergodic measures: For arbitrary measures $\mu$  no analogue of Fact \ref{frequencies-measure} can possibly be true, so that one has to decompose such $\mu$ as convex combination of ergodic measures (using Proposition \ref{Fer-Mont}). 

\medskip

Another alternative approach for substitutional subshifts (and in part also for $S$-adic subshifts), 
based on the machinery of Bratteli diagrams and Vershik maps, has been 
pursued 
by Bezuglyi, Kwiatkowski, Medynets and Solomyak in \cite{BKM, BKMS, BKMS2}. A brief introduction to Bratteli-Vershik technology and a discussion of their work together with a comparison to our results is given below in Appendix \ref{Bratteli-Vershik-presentation}.


\section{Graph maps in several different dialects}
\label{graph-dialects}
In this section we set up the basic terminology for the technology of graph towers, weight towers and vector towers, which will be defined and studied in the subsequent sections. 
\begin{convention}
\label{graph-convention}
\rm
Throughout this section we use the terminology defined in section \ref{set-up-I}. For the convenience of the reader we recall:
\begin{enumerate}
\item
All graphs $\Gamma$ are finite and connected, and possibly with vertices of valence 1 or valence 2.
\item
A {graph map} $f: \Gamma' \to \Gamma$ is a map between graphs $\Gamma'$ and $\Gamma$ which maps vertices to vertices and edges to edge paths.
\end{enumerate}
\end{convention}

Without further assumptions, for any graph map $f: \Gamma' \to \Gamma$ the image edge path $f(e')$ in $\Gamma$ of an edge $e'$ of $\Gamma'$ may well not be reduced. It could also be a trivial edge path, i.e. $e$ is contracted by $f$ to a single vertex. 

However, in this paper we will almost always be concerned with graph maps that maps edges to non-trivial reduced paths. To achieve this condition, starting from a general graph map $f: \Gamma' \to \Gamma$, we pass to a quotient graph of $\Gamma'$ by contracting all edges that are contracted by $f$, and we introduce additional vertices of valence 2 by subdividing edges of $\Gamma'$ at those points where $f$ is not locally injective. The resulting quotient graph of $\Gamma'$ and the graph map induced on it by $f$ has now the desired property.

\begin{defn-rem}
\label{folding-vertex}
The reader should be aware of the following ``pathology'', which may well occur for a graph map $f: \Gamma' \to \Gamma$, even if all edges have image paths that are non-trivial reduced:

At some vertex $v'$ of $\Gamma'$ the edges $e'_j$ with $v'$ as initial vertex are mapped to edge paths $f(e'_j)$ which all have the same edge $e_i$ as initial edge. In this case we say that $v'$ is an {\em $f$-folding vertex}. This happens in particular if $v'$ is mapped to a vertex of valence 1.
\end{defn-rem}

We will now define three different ``dialects'', in which graphs and graph maps can occur, as well as the formal transition between them.
This will be done below with all technical details, since it is the base for what comes in the subsequent sections.

\begin{figure}[!t]

\centering
\begin{tikzpicture}

\draw (-4,0.5) circle (0.5) ;
\draw (-4,-0.5) circle (0.5) ;
\node (A)at(-4,0) {$\bullet$};
\draw (A) node[below] {$v^*$};
\draw (-4,1) node {$>$};
\draw (-4,1.4) node {$\alpha^*$};
\draw (-4,-1) node {$>$};
\draw (-4,-1.4) node {$\beta^*$};

\draw [->,thick] (-3.2,0) -- (-2.3,0) node[midway,above] {$f^*$} ;

\draw (-1.5,0.5) circle (0.5) ;
\draw (-1.5,-0.5) circle (0.5) ;
\node (B)at(-1.5,0) {$\bullet$};
\draw (B) node[below] {$w$};
\draw (-1.5,1) node {$>$};
\draw (-1.5,1.4) node {$a$};
\draw (-1.5,-1) node {$>$};
\draw (-1.5,-1.4) node {$b$};

\draw (-0.2,2.5) -- (-0.2,-4) ;

\draw (2,1) circle (1) ;
\draw (2,-1) circle (1) ;
\node (C)at(2,0) {$\bullet$};
\draw (C) node[below] {$\check{v}$};
\node (V1)at(1.14,1.5)  {$\bullet$} ;
\draw (V1) node[left] {$\check{v_1}$};
\node (V2)at(2.86,1.5)  {$\bullet$} ;
\draw (V2) node[right] {$\check{v_2}$};
\node (V3)at(2,-2)  {$\bullet$} ;
\draw (V3) node[below] {$\check{v_3}$};

\draw (2,2) node {$>$};
\draw (2,2.4) node {$\check{\alpha}_2$};
\draw (1.14,0.5) node[rotate=120] {$>$};
\draw (0.75,0.5) node {$\check{\beta}_1$};
\draw (2.86,0.5) node[rotate=60] {$<$};
\draw (3.25,0.5) node {$\check{\beta}_3$};

\draw (1,-1) node[rotate=90] {$<$};
\draw (0.6,-1) node {$\check{\beta}_1$};
\draw (3,-1) node[rotate=90] {$>$};
\draw (3.4,-1) node {$\check{\beta}_2$};

\draw [->,thick] (3.9,0) -- (4.8,0) node[midway,above] {$\check{f}$} ;

\draw (5.7,0.5) circle (0.5) ;
\draw (5.7,-0.5) circle (0.5) ;
\node (D)at(5.7,0) {$\bullet$};
\draw (D) node[below] {$w\;$};
\draw (5.7,1) node {$>$};
\draw (5.7,1.4) node {$a$};
\draw (5.7,-1) node {$>$};
\draw (5.7,-1.4) node {$b$};

\draw (-2.5,-2.5) node[below]{$
\begin{array}{rcl}
f^*(v^*) & = & w \\
f^*(\alpha^*) & = & aba \\
f^*(\beta^*) & = & ba 
\end{array}
$};

\draw (3,-2.5) node[below]{$
\begin{array}{rcl}
\check{f}(\check{v}) = \check{f}(\check{v}_i) & = & w \\
\check{f}(\check{\alpha}_1) =  \check{f}(\check{\alpha}_3) =  \check{f}(\check{\beta}_2) & = & a \\
\check{f}(\check{\alpha}_2) = \check{f}(\check{\beta}_1) & = & b
\end{array}
$};

\end{tikzpicture}

\caption{We consider the graph $\Gamma$ with one vertex $v$ and two edges $a$, $b$ (the ``rose with two petals''), as well as the graph map $f:\Gamma\rightarrow\Gamma$ defined by $f(a)=aba$ and $f(b)=ba$.
In long-edge dialect, this gives rise to a graph map $f^*:\Gamma^*\rightarrow\Gamma$ (represented on the left hand side).
In short-edge dialect, this gives rise to a graph map $\check{f}:\check{\Gamma}\rightarrow\Gamma$ (represented on the right hand side).
The correspondence between $f^*$ and $\check{f}$ is given by:
$\alpha^* \leftrightarrow \check{\alpha}_1\check{\alpha}_2\check{\alpha}_3$ and
$\beta^* \leftrightarrow \check{\beta}_1\check{\beta}_2$.
}
\label{fig:short edge}

\end{figure}
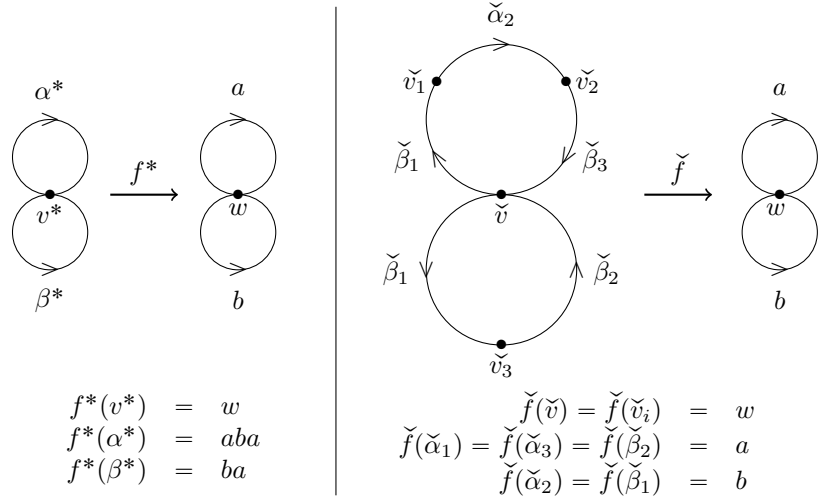

\begin{defn}[see Figure \ref{fig:short edge}]
\label{dialects}
(1)
A graph $\Gamma^*$ is given in {\em long-edge dialect} if $\Gamma^*$ has only {\em intrinsic} vertices, i.e. vertices of valence $\geq 3$
or of valence $1$. 
More generally, $\Gamma^*$ is in {\em long-edge dialect relative to some subset $V_0$ of the vertices of $\Gamma^*$} if all vertices of $\Gamma^*$ outside $V_0$ do not have valence 2.  

Edges of such a graph are called {\em long} edges, and we usually denote them by $e^*$ or $e_i^*$.

A graph map 
$f^*:\Gamma^* \to \Gamma$ is in {\em long-edge dialect} if 
the following conditions are satisfied:
\begin{enumerate}
\item[(a)]
The graph $\Gamma^*$ is in {long-edge dialect} relative to the set of all $f$-folding vertices.
\item[(b)]
For any edge $e^*$ of $\Gamma^*$ the image edge path $f^*(e^*)$ is non-trivial and reduced.
\end{enumerate}

\smallskip
\noindent
(2)
A graph map $\check f:\check\Gamma \to \Gamma$ is in {\em short-edge dialect} if for every edge $e$ of $\check \Gamma$ the image path $\check f(e)$ has length 1, or in other words: $\check f$ maps every edge to a single edge. The edges of such $\check \Gamma$ are called {\em short} edges.
\end{defn}

\begin{rem}
\label{long-short}
(1)
Let $f: \Gamma' \to \Gamma$ be a graph map that maps edges to non-trivial reduced edge paths.
The ``translation'' of $f$ and $\Gamma'$ into long-edge dialect is simply done 
by erasing all valence $2$ vertices from $\Gamma'$ which are not $f$-folding vertices.  We formalize this transition by calling the resulting graph $\Long(\Gamma')$ and the resulting map $\Long(f)$.

\smallskip
\noindent
(2)
In part (1) above we have purposefully ignored the very special case where $\Gamma'$ possesses only vertices of valence 2, i.e. $\Gamma'$ is homeomorphic to the circle. Of course, in this case one can only erase all but one vertices; however, since this case is not relevant to the goals of this paper and furthermore easy to figure out, we decided to leave out the notational extra efforts needed to formally treat this case.

\smallskip
\noindent
(3)
Similarly, the translation into short-edge dialect is given by introducing new valence 2 vertices in $\Gamma'$ for every $f$-preimage point of a vertex of $\Gamma$ (unless, of course, the preimage point is already a vertex of $\Gamma'$). Again, we formalize this transition by calling the resulting graph $\Short(\Gamma')$ and the resulting map $\Short(f)$.

The reader verifies directly the following equalities:
$$\Long(\Gamma') = \Long(\Long(\Gamma')) = \Long(\Short(\Gamma'))$$ 
$$\Long(f) = \Long(\Long(f)) = \Long(\Short(f))$$ 
$$\Short(\Gamma') = \Short(\Short(\Gamma')) = \Short(\Long(\Gamma'))$$ 
$$\Short(f) = \Short(\Short(f)) = \Short(\Long(f))$$ 
\end{rem}

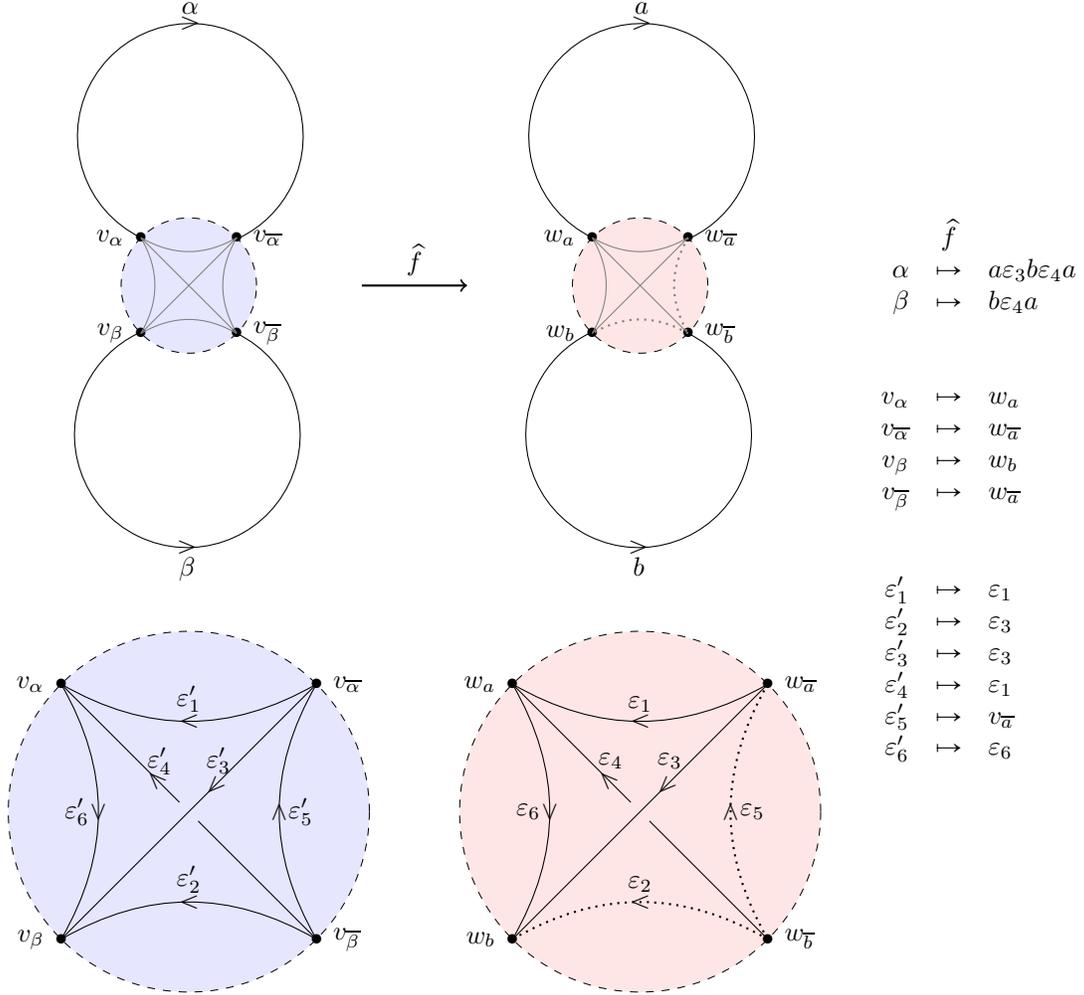
\begin{figure}[!t]

\centering

\begin{tikzpicture}[scale=1]

\node[draw,circle,fill=blue!10,ultra thin,dashed,minimum width=1.8cm] (A)at(-3,2) {};
\draw (A.north west) node{$\bullet$} node[left]{$v_\alpha\;$};
\draw (A.north east) node{$\bullet$} node[right]{$\;v_{\bar\alpha}$};
\draw (A.south west) node{$\bullet$} node[left]{$v_\beta\;$};
\draw (A.south east) node{$\bullet$} node[right]{$\;v_{\bar\beta}$};

\node[draw,circle,fill=red!10,ultra thin,dashed,minimum width=1.8cm] (B)at(3,2) {};
\draw (B.north west) node{$\bullet$} node[left]{$w_a\;$};
\draw (B.north east) node{$\bullet$} node[right]{$\;w_{\bar a}$};
\draw (B.south west) node{$\bullet$} node[left]{$w_b\;$};
\draw (B.south east) node{$\bullet$} node[right]{$\;w_{\bar b}$};

\draw [gray, thin] (A.north west) to[bend right] (A.north east);
\draw [gray, thin] (A.north east) to[bend right] (A.south east);
\draw [gray, thin] (A.south east) to[bend right] (A.south west);
\draw [gray, thin] (A.south west) to[bend right] (A.north west);
\draw [gray, thin] (A.north west) to (A.south east);
\draw [gray, thin] (A.north east) to (A.south west);

\draw [gray, thin] (B.north west) to[bend right] (B.north east);
\draw [gray, thick, dotted] (B.north east) to[bend right] (B.south east);
\draw [gray, thick, dotted] (B.south east) to[bend right] (B.south west);
\draw [gray, thin] (B.south west) to[bend right] (B.north west);
\draw [gray, thin] (B.north west) to (B.south east);
\draw [gray, thin] (B.north east) to (B.south west);

\draw [->,thick] (-0.7,2) -- (0.7,2) node[midway,above] {$\hat{f}$} ;

\draw (A.north west) arc (-116:-424:1.5) node[midway,above] {$\alpha$} node[midway] {$>$};
\draw (A.south east) arc (64:-244:1.5) node[midway,below] {$\beta$} node[midway] {$>$};

\draw (B.north west) arc (-116:-424:1.5) node[midway,above] {$a$} node[midway] {$>$};
\draw (B.south east) arc (64:-244:1.5) node[midway,below] {$b$} node[midway] {$>$};

\node[draw,circle,fill=blue!10,ultra thin,dashed,minimum width=4.8cm] (AA)at(-3,-5) {};
\node (AAA)at(-3,-5) {};
\draw (AA.north west) node{$\bullet$} node[left]{$v_\alpha\;$};
\draw (AA.north east) node{$\bullet$} node[right]{$\;v_{\bar\alpha}$};
\draw (AA.south west) node{$\bullet$} node[left]{$v_\beta\;$};
\draw (AA.south east) node{$\bullet$} node[right]{$\;v_{\bar\beta}$};

\node[draw,circle,fill=red!10,ultra thin,dashed,minimum width=4.8cm] (BB)at(3,-5) {};
\node (BBB)at(3,-5) {};
\draw (BB.north west) node{$\bullet$} node[left]{$w_a\;$};
\draw (BB.north east) node{$\bullet$} node[right]{$\;w_{\bar a}$};
\draw (BB.south west) node{$\bullet$} node[left]{$w_b\;$};
\draw (BB.south east) node{$\bullet$} node[right]{$\;w_{\bar b}$};

\draw (AA.north west) to[bend right] node[midway,above] {$\varepsilon'_1$} node[midway] {$<$} (AA.north east) ;
\draw (AA.north east) to[bend right] node[midway,right] {$\varepsilon'_5$} node[midway,rotate=90] {$>$} (AA.south east);
\draw (AA.south east) to[bend right] node[midway,above] {$\varepsilon'_2$} node[midway] {$<$} (AA.south west);
\draw (AA.south west) to[bend right] node[midway,left] {$\varepsilon'_6$} node[midway,rotate=90] {$<$} (AA.north west);
\draw (AA.north west) -- node[pos=0.65,right] {$\varepsilon'_4$} node[pos=0.8,rotate=-45] {$<$} (AAA) -- (AA.south east);
\draw (AA.north east) to node[pos=0.3,left] {$\varepsilon'_3$} node[pos=0.4,rotate=45] {$<$} (AA.south west);

\draw (BB.north west) to[bend right] node[midway,above] {$\epsilon_1$} node[midway] {$<$} (BB.north east) ;
\draw [thick,dotted] (BB.north east) to[bend right] node[midway,right] {$\epsilon_5$} node[midway,rotate=90] {$>$} (BB.south east);
\draw [thick,dotted] (BB.south east) to[bend right] node[midway,above] {$\epsilon_2$} node[midway] {$<$} (BB.south west);
\draw (BB.south west) to[bend right] node[midway,left] {$\epsilon_6$} node[midway,rotate=90] {$<$} (BB.north west);
\draw (BB.north west) -- node[pos=0.65,right] {$\epsilon_4$} node[pos=0.8,rotate=-45] {$<$} (BBB) -- (BB.south east);
\draw (BB.north east) to node[pos=0.3,left] {$\epsilon_3$} node[pos=0.4,rotate=45] {$<$} (BB.south west);

\draw (7.5,3) node[below]{$
\begin{array}{rcl}
  & \hat f & \\     
\alpha & \mapsto & a \epsilon_3 b \epsilon_4 a \\
\beta   & \mapsto & b \epsilon_4 a \\
 & & \\
 & & \\
v_\alpha           & \mapsto & w_a \\
v_{\bar \alpha} & \mapsto & w_{\bar a} \\
v_\beta             & \mapsto & w_b \\
v_{\bar \beta}   & \mapsto & w_{\bar a} \\
 & & \\
 & & \\  
\varepsilon'_1   & \mapsto & \epsilon_1 \\
\varepsilon'_2   & \mapsto & \epsilon_3 \\
\varepsilon'_3   & \mapsto & \epsilon_3 \\
\varepsilon'_4   & \mapsto & \epsilon_1 \\
\varepsilon'_5   & \mapsto & v_{\bar a} \\
\varepsilon'_6   & \mapsto & \epsilon_6 
\end{array}
$};

\end{tikzpicture}
\caption{
We consider the map $\hat f=\Blow(f^*)$, which is the translation of the map $f^*$ from Figure~\ref{fig:short edge} into the blow-up dialect.
For this, we replace the vertices $v$ and $w$ by local graphs $\Gamma(v)$ and $\Gamma(w)$ respectively: a zoom-in of these local graphs is reproduced below each of the two corresponding global graphs.
\newline
We note that not all of the blown-up local vertex graph $\Gamma(w)$ is in the image of $\hat f$: in practice, the dotted edges $\epsilon_2$ and $\epsilon_5$ could have been removed.
Notice also that the local edge $\epsilon'_5$ in $\Gamma(v)$ is illegal, while all other local edges in $\Gamma(v)$ are legal. 
}
\label{fig:blow up}

\end{figure}


\begin{defn}
\label{blow-up-dialect}
(1)
A graph $\hat \Gamma$ is given in {\em blow-up dialect} if the following conditions are satisfied:
\begin{enumerate}
\item[(a)]
The vertices of $\hat \Gamma$ are partitioned into classes:
$$V(\hat \Gamma) = V_1 \disjoint \ldots \disjoint V_q$$
Here $\disjoint$ denotes the disjoint union.
\item[(b)]
The edges of $\hat \Gamma$ are partitioned 
into classes:
$$\Edges^\pm(\hat \Gamma) = \hat E^\pm \disjoint \cal E^\pm_1 \disjoint \ldots \disjoint \cal E^\pm_q$$
Occasionally we will specify this notation to $\hat E^\pm(\hat \Gamma) := \hat E^\pm$ and $\cal E^\pm(\hat \Gamma) := \cal E^\pm_1 \cup \ldots \cup \cal E^\pm_q$.
\item[(c)]
For every $k = 1, \ldots, q$
the edges $\epsilon_j$ from $\cal E^\pm_k$ (called {\em local edges}) form a complete graph (called {\em local vertex graph}) over the vertex set $V_k$.
\item[(d)]
Every vertex of $\hat \Gamma$ is the initial vertex of precisely one edge $\hat e_i$ from the set $\hat E^\pm$ of {\em non-local} edges (and thus also the terminal vertex of precisely one such edge).
\end{enumerate}

\smallskip
\noindent
(2)
A graph map $\hat f: \hat \Gamma' \to \hat \Gamma$ is given in {\em blow-up dialect} if both, $\hat \Gamma$ and $\hat \Gamma'$ are in blow-up dialect, and if the map $\hat f$ maps every local vertex graph of $\hat\Gamma'$ to a local vertex graph of $\hat\Gamma$, such that the following holds:

\begin{enumerate}
\item[(a)]
Every local edge $\epsilon'_j$ of $\hat\Gamma'$ is either mapped to a single local edge $\epsilon_k = \hat f(\epsilon'_j)$ of $\hat\Gamma$, or else $\epsilon'_j$ is contracted by $\hat f$ to a vertex. In the first case 
the local edge 
$\epsilon'_j$ will be termed {\em legal}, while in the second case we call it {\em illegal}, with respect to the map $\hat f$.

\item[(b)]
We also require that for every non-local edge 
$\hat e' \in \hat E^\pm(\hat \Gamma')$
the image edge path $\hat f(\hat e')$ is reduced and non-trivial. Furthermore $\hat f(\hat e')$ 
does not have a local edge as initial or as terminal edge, and $\hat f(\hat e')$ never traverses two consecutive local edges. From the hypothesis that $\hat f(\hat e')$ is reduced it follows that $\hat f(\hat e')$ also never traverses two consecutive non-local edges, so that in fact it is an edge path that alternates between local and non-local edges.
\end{enumerate}
\end{defn}
 
\begin{rem}
\label{contraction}
(1)
Let $\hat \Gamma$ be a graph in blow-up dialect as in Definition \ref{blow-up-dialect}. The graph $\Gamma$ obtained from $\hat \Gamma$ by contracting all local edges of $\hat \Gamma$ (and hence identifying, for each $k = 1, \ldots, q$, all vertices in $V_k$ to define a single quotient vertex $\cal V_k$), is said to be obtained {\em by contraction}. We denote this by:
$$\Gamma = \Contr(\hat \Gamma)$$

\smallskip
\noindent
(2)
Let $\hat f: \hat \Gamma' \to \hat \Gamma$ be a graph map in blow-up dialect.  We say that the map $f: \Gamma' \to \Gamma$ is obtained from $\hat f$ {\em by contraction} if we have $\Gamma = \Contr(\hat \Gamma)$ and $\Gamma' = \Contr(\hat \Gamma')$, and $f$ is the map induced by $\hat f$ on the two quotient graphs. In this case we write:
$$f = \Contr(\hat f)$$
It follows directly from the above definitions that this map $f$ maps every edge of $\Gamma'$ to a reduced non-trivial edge path in $\Gamma$.
\end{rem}

We now want to describe the converse ``translation'' 
(see Figure \ref{fig:blow up}). 
For this purpose we first define a {\em blow-up} procedure at a vertex $v$ of the graph $\Gamma$:  
Let $E(v)$ be the set of oriented edges $e$ which have $v$ as initial vertex 
(to be specific: if some edge $e$ has $v$ as initial and also as terminal vertex, then both, $e$ and $\bar e$ belong to $E(v)$). We define  a {\em local vertex graph} $\Gamma(v)$, which has a vertex $v_e$ for each $e \in E(v)$ and is the full graph over this {\em local vertex set} $\{ v_e\mid e \in E(v) \}$. The edges of such a local graph $\Gamma(v)$ are called {\em local edges} and will be denoted by $\epsilon$ or $\epsilon_k$.

\begin{defn-rem}
\label{blow-up}
(1)
For any graph $\Gamma$ the associated {\em blow-up graph} $\hat \Gamma$ is defined as the union of the local vertex graphs $\Gamma(v)$, for any vertex $v$ of $\Gamma$, together with an edge $\hat e$ for every edge $e$ of $\Gamma$: if $e$ has initial vertex $v_1$ and terminal vertex $v_2$, then the initial vertex of $\hat e$ is the {\em local vertex} $v_e$ of $\Gamma(v_1)$, and the terminal vertex of $\hat e$ is the local vertex $v_{\bar e}$ of $\Gamma(v_2)$. One verifies easily that the conditions (a) - (d) of Definition \ref{blow-up-dialect} (1) are satisfied. We write:
$$\hat \Gamma = \Blow(\Gamma)$$

\smallskip
\noindent
(2)
Given a graph map $f: \Gamma' \to \Gamma$ which maps edges to reduced non-trivial edge paths, we define the associated {\em blow-up map} $\hat f: \hat \Gamma' \to \hat \Gamma$ by passing to the blow-up graphs $\hat \Gamma := \Blow(\Gamma)$ and $\hat \Gamma' := \Blow(\Gamma')$. We now consider any vertex $v'$ of $\Gamma'$ and any edge $e'$ with $v'$ as initial vertex, and for the vertex $v_{e'}$ of the local vertex graph $\Gamma(v')$ we define $\hat f(v_{e'}) = v_e$, where $e$ is the initial edge of the edge path $f(e')$ (so that $v_e$ is a vertex of the local vertex graph $\Gamma(f(v'))$).

For any edge $e'$ of $\Gamma'$ with $f(e') = e_1 e_2 \ldots e_t$ we define $\hat f(\hat e') := \hat e_1 \epsilon_1 \hat e_2 \epsilon_2 \ldots \epsilon_{t-1} \hat e_t$, where $\epsilon_k$ is the local edge that connects the terminal vertex $v_{\bar e_k}$ of $\hat e_k$ to the initial vertex $v_{e_{k+1}}$ of $\hat e_{k+1}$. (Note that such a local edge must exist, since the terminal vertex of $e_k$ agrees with the initial vertex of $e_{k+1}$ in $\Gamma$, and since from the assumption that $f(e')$ is reduced it follows that $v_{\bar e_k} \neq v_{e_{k+1}}$.) 

For any local edge $\epsilon'_j$ of $\hat\Gamma'$ which connects a vertex $v_1$ to a vertex $v_2$, the image $\hat f(\epsilon'_j)$ is either defined to be the local edge connecting $\hat f(v_1)$ to $\hat f(v_2)$, in case that these two vertices are distinct, or else $\hat f(\epsilon'_j)$ is contracted to the single vertex $\hat f(v_1) = \hat f(v_2)$. 

Again, one sees directly that the map $\hat f$ is in blow-up dialect, as set up above in Definition \ref{blow-up-dialect} (2).
We write:
$$\hat f = \Blow(f)$$
\end{defn-rem}

At this point we would like to point out that a vertex $v'$ of $\Gamma'$ is an $f$-folding vertex for the graph map $f: \Gamma' \to \Gamma$ if and only if, after translation into blow-up dialect, all edges of the local vertex graph $\Gamma(v')$ are illegal with respect to $\hat f$ (see Definition \ref{blow-up-dialect} (2)(a)).

\begin{rem}
\label{no-name}
The reader verifies directly from the definitions the following equalities, for any graph $\Gamma$, or for any graph map $f: \Gamma' \to \Gamma$
which maps edges to reduced non-trivial edge paths:
$$\Gamma = \Contr(\Blow(\Gamma))$$
$$f = \Contr(\Blow(f))$$
Similarly, for any graph $\hat \Gamma$ and any graph map $f: \hat \Gamma' \to \hat\Gamma$ in blow-up dialect we have:
$$\hat\Gamma = \Blow(\Contr(\hat\Gamma))$$
$$\hat f = \Blow(\Contr(\hat f))$$
\end{rem}

In the next sections the blow-up dialect will almost always be used in combination with the short-edge dialect, i.e. we will consider, for a given graph map $f: \Gamma' \to \Gamma$ the combined translations $\Blow(\Short(f)): \Blow(\Short(\Gamma')) \to \Blow((\Gamma))$.

However, our above set-up allows also to use the blow-up dialect in combination with the long-edge dialect, i.e. one works with the maps $\Blow(\Long(f)): \Blow(\Long(\Gamma')) \to \Blow(\Gamma)$.

\begin{convention}
\label{translation-paths}
\rm
In the subsequent sections we will often pass in an informal way from one dialect to the other. In this case we use the following convention, for any graph map $f: \Gamma' \to \Gamma$, and any reduced non-trivial edge path $\gamma$ in $\Gamma'$:

The path $\gamma$ will not change name if we pass to long-edge or to short-edge dialect: Indeed, $\gamma$ stays topologically the same, as simply valence 2 vertices will be added or removed. In long-edge dialect it could hence be that $\gamma$ is not any more an edge path in the classical sense, but starts and finishes with a ``partial edge'' (or $\gamma$ may also be entirely contained in a single long edge).

In the case where we pass to the blow-up dialect, the name $\gamma$ still stays, but in addition we impose that in the blow-up graph $\hat \Gamma'$ the corresponding path $\gamma$ never starts or ends with a local edge, and never passes over two consecutive local edges. 

It is a direct consequence of the above conventions that changing back and forth dialects will not change $\gamma$ if after several changes one ends up in the same dialect as started out with. Here we need to assume, if we start out in blow-up dialect, that $\gamma$ does not start or end with a local edge, and does never pass over two consecutive local edges.
\end{convention}

For the subsequent sections we need to generalize
the above defined concept of ``long-edge dialect'' slightly:  

\begin{defn-rem}
\label{rel-long}
Let $I$ be a (possibly infinite) index set, $J$ a subset of $I \times I$, and let $(f_{i, j}: \Gamma_j \to \Gamma_i)_{(i, j) \in J}$ be a family of graph maps which is stable with respect to 
composition. Then we define the {\em relative long-edge dialect (with respect to $J$)} by erasing all non-intrinsic vertices from the $\Gamma_i$, where a vertex $v$ of $\Gamma_i$ is counted as {\em intrinsic} if one of the following conditions is satisfied:
\begin{enumerate}
\item
$v$ has valence $1$ or $\geq 3$.
\item
$v$ is an $f_{k, i}$-folding vertex for some $(k,i) \in J$.
\item
$v = f_{i, j}(v')$ for some $(i, j) \in J$, and $v'$ is an intrinsic vertex as in (1) or (2) above.
\end{enumerate}

Equivalently, a vertex $v$ is non-intrinsic if it is of valence 2, and so are all of the preimage vertices of $v$, and 
furthermore neither $v$ nor any of these preimage vertices is a folding vertex, for any map from the given family $J$ of graph maps.

 As a consequence, intrinsic vertices are always mapped to intrinsic vertices, while non-intrinsic vertices are mapped either to intrinsic or to non-intrinsic vertices.
\end{defn-rem}

 
\section{Graph towers and their used language}\label{towers}

In this section we will properly define and study the main tool of this paper, called ``expanding graph towers''. To every expanding graph tower we canonically associate a symbolic lamination, called the ``used lamination''. We also show that to any symbolic lamination $\Lsig$ there exists an expanding graph tower which has $\Lsig$ as used lamination.

\subsection{Graph towers}

${}^{}$

We will first recall from section \ref{tools-results} the definition of graph towers,  without reference to any of the three dialects introduced in section \ref{graph-dialects}. We comment below about the translation into these dialects.

\begin{defn}
\label{graph-towers-1}
A {\em graph tower} $\bvec \Gamma$ is given by an infinite family $(\Gamma_n)_{n \in \N \cup \{0\}}$ of finite connected {\em level graphs} $\Gamma_n$ and an infinite family $\bvec f = (f_{m, n})_{0 \leq m \leq n}$ of graph maps $f_{m, n}: \Gamma_n \to \Gamma_m$ (sometimes called {\em level transition maps}), where the following properties are satisfied:
\begin{enumerate}
\item[(a)]
$f_{m,n}$ maps vertices to vertices.

\item[(b)]
$f_{m,n}$ maps edges to reduced non-trivial edge paths.

\item[(c)]
The family $\bvec f$ is {\em compatible}:  one has $f_{k,m} \circ f_{m,n} = f_{k, n}$ for all integers $n \geq m \geq k \geq 0$.
In particular we require $f_{n,n} = id_{\Gamma_n}$ for all $n \geq 0$.
\end{enumerate}
For simplicity we will use the abbreviations  $f_n := f_{0, n}$ for all $n \geq 0$.
\end{defn}

\begin{rem}
\label{gt-translation}
(1)
Any graph tower $\bvec \Gamma = ((\Gamma_n)_{n \in \N \cup \{0\}} , (f_{m,n})_{0 \leq m\leq n})$ can be translated canonically into
\begin{enumerate}
\item[$\bullet$]
 a {\em relative long-edge graph tower} $((\Gamma^*_n)_{n \in \N \cup \{0\}}, (f^*_{m,n})_{0 \leq m\leq n})$ 
\newline 
[where ``relative'' is meant with respect to the family $(f_{m,n})_{0 \leq m\leq n}$, see Definition-Remark \ref{rel-long}], 
 \item[$\bullet$]
 a {\em short-edge graph tower} $((\check \Gamma_n)_{n \in \N \cup \{0\}} , (\check f_{m,n})_{0 \leq m\leq n})$, or 
\item[$\bullet$]
 a {\em blow-up graph tower} $((\hat \Gamma_n)_{n \in \N \cup \{0\}} , (\hat f_{m,n})_{0 \leq m\leq n})$.
\end{enumerate}
The translation of $\bvec \Gamma$ into short-edge or blow-up dialect 
is always done ``from the bottom up'', 
where at any step one always follows the detailed instructions given 
in section \ref{graph-dialects}: One first translates $\Gamma_0$, then $\Gamma_1$ together with $f_1$, then $\Gamma_2$ together with $f_{1, 2}$ (and thus also with $f_2 = f_1 \circ f_{1,2}$), and so on. 

\smallskip
\noindent
(2)
As a consequence, we note for any level graph $\Gamma_n$ that, in the process of translating $\Gamma_n$ into short-edge dialect $\check \Gamma_n$ through subdivision of the edges according to any of the maps $f_{m,n}$ (see Remark \ref{long-short}), the set of newly introduced valence 2 vertices is independent of the choice of $m$, since any of the level graphs $\Gamma_m$ with $m \leq n-1$ has (in the procedure ``from the bottom up'') already been translated into short-edge dialect.

\smallskip
\noindent
(3)
Similarly, we note that if in the blow-up dialect any level graph $\hat \Gamma_n$ has a local edge $\epsilon_k$ which is {\em illegal}, by which we mean ``illegal with respect to the map $\hat f_n$'' (see Definition \ref{blow-up-dialect} (2)), then for any level graph $\hat \Gamma_m$ of lower level $m \leq n$ the image $\hat f_{m,n}(\epsilon_k)$ is either degenerated to a single vertex, or else $\hat f_{m,n}(\epsilon_k)$ is a local edge which then must also be illegal (i.e. illegal with respect to $\hat f_m$).

\smallskip
\noindent
(4)
For the translation of $\bvec \Gamma$ into relative long-edge dialect we first determine for any level graph $\Gamma_n$ the set of vertices that have valence $\geq 3$ or $1$, or are $f_n$-folding vertices (where we observe that any $f_{m, n}$-folding vertex must also be an $f_n$-folding vertex, by the compatibility condition (c) in Definition \ref{graph-towers-1}). These vertices together with their $f_{m, n}$-images, for any $1 \leq m \leq n$, are the intrinsic vertices of $\bvec \Gamma$. Erasing all non-intrinsic vertices from any level graph $\Gamma_m$ with $m \geq 1$ (compare Remark \ref{long-short} (2))
now gives directly the long-edge graph tower associated to $\bvec \Gamma$.
\end{rem}

\begin{rem}
\label{pathological}
(1)
The concept of a graph tower $\bvec \Gamma$ given in Definition \ref{graph-towers-1} has been purposefully kept fairly general. As a consequence, there can occur two types of ``pathologies'', which we'd like to point out:
\begin{enumerate}
\item[(a)]
Some of the level transition maps $f_{m,n}$ may have folding vertices.
\item[(b)]
Some of the level transition maps $f_{m,n}$ may not be surjective.
\end{enumerate}
Both of these pathologies don't harm the proofs given below, but they may lead to unnecessary complications, so that the reader may be more happy to work 
with ``non-pathological'' graph towers, i.e. (a) and (b) above do not occur. This is in fact possible without any essential loss, see Remark \ref{not-lam} (2) and Remark \ref{equal-weights1} (2).

\smallskip
\noindent
(2)
At this point we'd like to note that any graph tower $\bvec \Gamma$ possesses a canonical maximal non-pathological subtower, obtained by iteratively removing all (perhaps infinitely many) edges which are not in the image of any level transition map, and also any edge which is adjacent to a folding vertex. One also removes any vertex which in this iterative process eventually gets valence 0.

Of course, this may lead to degeneracy phenomena like graphs homeomorphic to a circle with level maps that are homeomorphisms, or even empty graphs. Such phenomena lie outside the realm of interest of this paper. The latter is captured properly by the subsequent definition, as in the sequel we will restrict our attention exclusively to expanding graph towers. It is easy to verify that for those the above named degeneracy phenomena regarding the maximal non-pathological subtower can not occur.
\end{rem}

\begin{defn}
\label{expanding}
We say that the graph tower $\bvec \Gamma$, given by a family $\bvec f$ of graph maps as in Definition \ref{graph-towers-1}, is {\em expanding} if, 
when considering the long edges $e^*_i$ obtained from deleting the non-intrinsic vertices of the level graphs $\Gamma_n$ (i.e. by passing over to level graphs $\Gamma^*_n$ through translation into relative long-edge dialect) the {\em minimal long-edge length} 
$$
{\mile} _{\tiny\bvec \Gamma}(n) := \min_{e^*_i \in {\rm Edges}^\pm(\Gamma^*_n)} |f^*_{n}(e^*_i)|$$
(see Definition \ref{length-defn}) satisfies
$$\lim_{n \to \infty} ({\mile} _{\tiny\bvec \Gamma}(n))\to \infty \, .$$
\end{defn}

\medskip

\subsection{Languages and laminations associated to an expanding graph tower}

${}^{}$

To any graph tower there are naturally associated two (competing) natural ``languages''; the difference is subtle but in some contexts 
important, so that we first want to set up our terminology with a bit of care:

\begin{defn}
\label{language+}
Let $\Gamma$ be any finite connected graph.
\begin{enumerate}
\item
We denote by $\cal P(\Gamma)$ the set of finite edge paths $\gamma$ in $\Gamma$ that are reduced.
\item
Any subset $\cal L \subset \cal P(\Gamma)$ is called a {\em language over $\Gamma$}.
\item
For any language $\cal L$ over $\Gamma$ we use the following terminology:
\begin{enumerate}
\item
$\cal L$ is called {\em factorial} if for any $\gamma \in \cal L$ one also has $\gamma' \in \cal L$, for any subpath (``factor'') $\gamma'$ of $\gamma$. 
\item
$\cal L$ is called {\em bi-extendable} if for any $\gamma \in \cal L$ there exists a  $\gamma' \in \cal L$ such that $\gamma$ is a non-boundary subpath of $\gamma'$. 

\qquad
Here a subpath $\gamma$ of $\gamma'$ is a ``boundary subpath'' if $\gamma$ and $\gamma'$ have their first or last edge in common, or, in case that $\gamma$ is trivial, if it agrees with the initial or terminal vertex of $\gamma'$.
\item
$\cal L$ is called {\em invertible} if for any $\gamma \in \cal L$ the reversed path $\bar \gamma$ also satisfies $\bar \gamma \in \cal L$.
\item
$\cal L$ is called {\em laminary} if it is 
non-empty, factorial, bi-extendable and invertible.
\end{enumerate}
\item
Every infinite language $\cal L$ over $\Gamma$ defines an {\em associated} symbolic lamination $\Lsig(\cal L)$ on $\Gamma$: it consist precisely of those biinfinite paths $\gamma$ which have the property that every finite subpath $\gamma'$ of $\gamma$ is also a subpath of some $\gamma'' \in \cal L$ or of its inverse $\bar \gamma''$.
\item
Conversely, every symbolic lamination $\Lsig$ on $\Gamma$ defines an {\em associated} language $\cal L(\Lsig)$ over $\Gamma$, which consists of all finite subpaths of any $\gamma \in \Lsig$.
\end{enumerate}
\end{defn}

On easily verifies that for any symbolic lamination $\Lsig$ on a graph $\Gamma$ the associated language $\cal L(\Lsig)$ is laminary, and that one has $\Lsig(\cal L(\Lsig)) = \Lsig$. Conversely, however, $\cal L(\Lsig(\cal L))= \cal L$ only follows if $\cal L$ is laminary; in general neither of the two languages is included in the other.

\begin{defn}
\label{used-lamination+}
Let
$\bvec \Gamma = ((\Gamma_n)_{n \in \N \cup \{0\}} , (f_{m,n})_{0 \leq m\leq n})$
be an expanding 
graph tower. 
\begin{enumerate}
\item
A finite edge path $\gamma$ in some level graph $\Gamma_m$ is called {\em used} if for any index $n \geq m$ there is an 
index $n' \geq n$ and an 
edge $e^*$ in the relative long-edge level graph $\Gamma^*_{n'}$ such that the path $f^*_{m,n'}(e^*)$ contains $\gamma$ or $\bar \gamma$ as subpath. The set of used edge paths in $\Gamma_0$ is called the {\em used language of $\bvec \Gamma$} and denoted by $\cal L_{used}(\bvec \Gamma)$.

\item
An edge path $\gamma$ in some level graph $\Gamma_n$ is {\em legal} if $f_n(\gamma)$ is reduced, or, equivalently, if $\gamma$ translated into in blow-up dialect only crosses over local edges that are legal. The set of legal paths in $\Gamma_n$ is denoted by $\cal L_{legal}(\Gamma_n)$. An edge path $\gamma$ in $\Gamma_0$ is called {\em infinitely legal} if for any $n \geq0$ it is a subpath of $f_n(\gamma')$ or $f_n(\bar \gamma')$, for some legal path $\gamma'$ in $\Gamma_n$. The set of infinitely legal edge paths in $\Gamma_0$ is called the {\em infinitely legal language of $\bvec \Gamma$} and denoted by $\cal L_{legal}(\bvec \Gamma)$.

\item
We define the {\em used symbolic lamination of $\bvec \Gamma$} by 
$$\Lsig_{used}(\bvec \Gamma)= \Lsig(\cal L_{used}(\bvec \Gamma))$$
and the {\em infinitely legal symbolic lamination of $\bvec \Gamma$} by 
$$\Lsig_{legal}(\bvec \Gamma) = \Lsig(\cal L_{legal}(\bvec \Gamma)) \, .$$
\end{enumerate}
\end{defn}

Before proceeding further, we would like to illustrate these definition with an example based on Figure \ref{fig:blow up}: We consider a graph tower $\bvec \Gamma$ where every level transition map is identical to the map considered there (a ``stationary''  graph tower as considered later in section \ref{stationary-pseudo}). To be more concrete, we pass over to blow-up dialect and extend the above definitions in the obvious way.
As pointed out in the caption of Figure \ref{fig:blow up}, the local edge $\epsilon_5$ is illegal, for every level $n \geq 1$, while all other local edges at all levels are legal. In fact, we verify easily that on level $n= 0$ they are actually infinitely legal.
The local edges $\epsilon_3$ and $\epsilon_4$ are used by the images of the non-local edges from the next level, and $\epsilon_1$ is used by the image of $\epsilon_4$, so that it also counts as ``used''. The local edges $\epsilon_5$ and $\epsilon_6$, although legal, are not ``used'', since (even after passing to preimages in arbitrary high level graphs) they are only crossed over by images of other local edges, but not of any non-local edge.

\begin{rem}
\label{camparison}
(1) From the hypothesis in Definition \ref{used-lamination+} that $\bvec \Gamma$ be expanding it follows directly that the used language $\cal L_{used}(\bvec \Gamma)$ and also the infinitely legal language $\cal L_{legal}(\bvec \Gamma)$ are infinite sets, so that $L^\Sigma_{used}(\bvec \Gamma)$ and $L^\Sigma_{legal}(\bvec \Gamma)$ are indeed symbolic laminations (see Definition \ref{language+} (4)).

\smallskip
\noindent
(2)
From the above definitions we derive directly
$$\cal L_{used}(\bvec \Gamma) \subseteq \cal L_{legal}(\bvec \Gamma)\qquad \text{and hence} \qquad
L^\Sigma_{used}(\bvec \Gamma) \subseteq L^\Sigma_{legal}(\bvec \Gamma)
\, .$$
Furthermore, the difference $L^\Sigma_{legal}(\bvec \Gamma) \smallsetminus L^\Sigma_{used}(\bvec \Gamma)$ is often not very large: 
Indeed, it can be shown that for expanding graph towers $\bvec \Gamma$ with uniformly bounded number of intrinsic vertices in any level graph this difference consists only of finitely many shift-orbits.
\end{rem}

\begin{rem}
\label{not-lam}
(1)
The reader should be cautioned about assuming too quickly that any of the above defines languages $\cal L_{used}(\bvec \Gamma), \cal L_{legal}(\Gamma_n)$ or $\cal L_{legal}(\bvec \Gamma)$ is laminary. In general this will not be the case, but for $\cal L_{legal}(\Gamma_n)$ or $\cal L_{legal}(\bvec \Gamma)$ this follows if one assumes that the graph tower $\bvec \Gamma$ is non-pathological in the sense of Remark \ref{pathological}.

\smallskip
\noindent
(2) It turns out that the used lamination of a pathological graph tower $\bvec \Gamma$ (see Remark \ref{pathological}) contains but is not quite equal to 
that of its maximal non-pathological subtower. However, this is of no serious concern for the purposes pursued in this paper, as the difference 
concerns only the ``inactive'' part of the lamination, by which we mean the union of all those leaves of $L^\Sigma_{used}(\bvec \Gamma)$ which are not in the support of any invariant measure on $L^\Sigma_{used}(\bvec \Gamma)$ (see Remark \ref{disconnected} (2)).
\end{rem}

\begin{defn-rem}
\label{name-used-lam}
In some contexts it may seem more natural to work with infinitely legal rather than with used edge paths. However, since used edge paths are closer to standard objects considered in classical symbolic dynamics (for example in the context of substitution subshifts), we have chosen to concentrate in the subsequent sections on $\cal L_{used}(\bvec \Gamma)$ and $L^\Sigma_{used}(\bvec \Gamma)$. For simplicity of the notation we introduce
$$\Lused := L^\Sigma_{used}(\bvec \Gamma)$$
as name for the used lamination of any expanding graph tower $\bvec \Gamma$.
\end{defn-rem}

\medskip

We show now that every symbolic lamination is the used lamination of some expanding graph tower, 
as has been claimed in Proposition \ref{any-lamination} (restated below for the convenience of the reader). 
The construction used 
in the proof 
has some reminiscences with Rauzy graphs, but the latter give a graph tower that in general may not be expanding. 

\medskip
\noindent
{\bf Proposition \ref{any-lamination}.}
Let $\Gamma$ be a finite graph, and let $L^\Sigma \subset \Sigma(\Gamma)$ be an arbitrary symbolic lamination on $\Gamma$. Then there exists an expanding graph tower $\bvec \Gamma = ((\Gamma_n)_{n \in \N \cup \{0\}}, (f_{m,n})_{0\leq m\leq n})$ and a graph map $f: \Gamma \to \Gamma_0$ which induces a bijection
$$f_{L^\Sigma}:L^\Sigma \to 
\Lused
\, .$$

\begin{proof}
In order to define $\bvec \Gamma$ we consider, for any integer $n \geq 1$, the set $\cal L_n \subset \cal L(L^\Sigma)$ of all paths of length $n$ of the language $\cal L(L^\Sigma)$.
For any $n \in \N \cup \{0\}$ we now define $\Gamma_n$ as graph with a single vertex $v_n$ and edges $e_\gamma$ for any $\gamma \in \cal L_{2^{n}}$, with the convention 
$$\bar e_\gamma := e_{\bar \gamma} \, .$$
We define the maps $f_{n, n+1}$ by setting $f_{n, n+1}(v_{n+1}) := v_{n}$ and, for any $\gamma \in \cal L_{2^{n+1}}$, by
$f_{n, n+1}(e_\gamma) := e_{\gamma'} e_{\gamma''}$, where $\gamma'$ and $\gamma''$ are the initial and terminal subpaths of $\gamma$ of length $2^{n}$. The maps $f_{m,n}$ for arbitrary $0\leq m \leq n$ are obtained from composing appropriately the $f_{n, n+1}$. It follows directly that the resulting graph tower $\bvec \Gamma$ is expanding.

The map $f: \Gamma \to \Gamma_0$ is defined by sending every vertex of $\Gamma$ to $v_0$ and any edge $e'$ of $\Gamma$, interpreted as edge path of length $1 = 2^0$, to the edge $e_{e'}$. As this map is nothing else than the quotient map which identifies all vertices of $\Gamma$ to a single vertex, it follows directly that the map induced by $f$ on $\Sigma(\Gamma)$ is injective. From the construction we also see directly that $f(\Sigma(\Gamma))$ is shift-invariant, and that $f$ commutes with the shift operator. In particular we observe the equality
\begin{equation}
\label{silly}
\cal L(f(L^\Sigma)) = f(\cal L(L^\Sigma)) \, .
\end{equation}

From the injectivity of $f$ we obtain a  
bijection
$f_{L^\Sigma}:L^\Sigma \to f(L^\Sigma)$,
so that it only remains to show that $f(L^\Sigma) = \Lused$. This in turn follows if we can show that the corresponding two laminary languages $\cal L(f(L^\Sigma))$ and $\cal L(\Lused)$ agree.

Since every path $\gamma'$ of $\cal L(L^\Sigma)$ is a subpath of some $\gamma \in \cal L_{2^n}$, for all sufficiently large $n$, it follows that $f(\gamma')$ is a subpath of $f_n(e_\gamma)$. 
From the definition of $\Lused$
we 
thus deduce that $f(\cal L(L^\Sigma)) \subset \cal L(\Lused)$, so that we can use the above equality (\ref{silly}) to obtain the inclusion 
$$\cal L(f(L^\Sigma)) \subset \cal L(\Lused)\, .$$

On the other hand, every path $\gamma''$ of $\cal L(\Lused)$ is a subpath of the path $f_n(e_\gamma)$ for some $\gamma \in \cal L_{2^n}$ with $n$ sufficiently large. Hence $\gamma'' = f(\gamma')$, where $\gamma'$ is a subpath of $\gamma$. This shows $\cal L(\Lused) \subset f(\cal L(L^\Sigma))$ and hence, again using (\ref{silly}), the inclusion 
$$\cal L(\Lused) \subset \cal L(f(L^\Sigma)) \, .$$
\end{proof}

We'd like to note that the graph tower $\bvec \Gamma$ constructed in the last proof is indeed non-pathological, in the meaning of Remark \ref{pathological}.

 
\section{Weights 
on graph towers}
\label{weights-currents}

\begin{convention}
\label{ft-short}
\rm
In this section we assume, unless explicitly otherwise specified, that any graph tower $\bvec \Gamma = ((\Gamma_n)_{n \in \N \cup \{0\}} , (f_{m,n})_{0 \leq m\leq n})$ is given in short-edge dialect (see Definition \ref{long-short} and Remark \ref{gt-translation}). In other words: we have $\Gamma_n = \check \Gamma_n$ and $f_{m, n} = \check f_{m, n}$ for all $n \geq m \geq 0$. 
\end{convention}

\begin{defn}
\label{weights+1}
Let $\Gamma$ be a graph and let $\hat \Gamma$ be the associated blow-up graph as in Definition-Remark \ref{blow-up} (1). We first define a {\em weight function} $\hat \omega$ on $\hat \Gamma$: 
This is a non-negative function 
$$
\hat\omega: \Edges^\pm(\hat \Gamma) \to \R_{\geq 0} \quad {\rm with} \quad \hat\omega(\bar{e}) = \hat\omega(e) \quad {\rm for \,\, all} \quad e \in \Edges^\pm(\hat \Gamma)$$
which in addition satisfies the following {\em switch conditions}: For every vertex $v$ of $\hat \Gamma$ one has 
\begin{equation}
\label{switch-conditions}
\hat\omega(\hat e) = \underset{\epsilon_k \in E(v)}{\sum} \hat\omega(\epsilon_k)\, ,
\end{equation}
where $\hat e$
is the non-local edge which has $v$ as initial vertex,
and $E(v)$ is the set of all local edges $\epsilon_k$
which also have $v$ as initial vertex. 

A {\em weight function} $\omega$ on $\Gamma$ is a function 
$$
\omega: \Edges^\pm(\Gamma) \to \R_{\geq 0}
$$
that is induced by some weight function $\hat \omega$ on the associated blow-up graph $\hat \Gamma$, i.e. for every edge $e$ of $\Gamma$ and the associated non-local edge $\hat e$ of $\hat \Gamma$ one has $\omega(e) = \hat\omega(\hat e)$.
\end{defn}


\begin{figure}[!t]

\centering{
\begin{tikzpicture}[scale=0.8]

\draw (-4,5.5) node[below]{$\Gamma_2$};

\draw (-4,2) circle (2) ;
\draw (-4,-2) circle (2) ;

\draw[red, ultra thick] (-4,0) arc (-90:-135:2) ;
\draw[red, ultra thick] (-4,4) arc (90:135:2) ;
\draw[red, ultra thick] (-2,2) arc (0:45:2) ;
\draw[red, ultra thick] (-4,0) arc (90:141.42:2) ;
\draw[red, ultra thick] (-4.86,-3.8) arc (-115.74:-167.16:2) ;
\draw[red, ultra thick] (-3.13,-3.8) arc (-64.26:-12.84:2) ;

\draw[fill=blue!50] (-4,0) circle (0.15cm);
\node (V7)at(-2.58,0.58)  {$\bullet$} ;
\node (V6)at(-2,2)  {$\bullet$} ;
\node (V5)at(-2.58,3.41)  {$\bullet$} ;
\node (V4)at(-4,4)  {$\bullet$} ;
\node (V3)at(-5.41,3.41)  {$\bullet$} ;
\node (V2)at(-6,2)  {$\bullet$} ;
\node (V1)at(-5.41,0.58)  {$\bullet$} ;

\draw[red] (-4.76,0.15) node[rotate=150] {$>$};
\draw (-4.76,0.5) node {$2$};
\draw (-5.84,1.23) node[rotate=120] {$>$};
\draw (-6.2,1.23) node {$2$};
\draw (-5.84,2.77) node[rotate=60] {$>$};
\draw (-6.2,2.77) node {$2$};
\draw[red]  (-4.76,3.85) node[rotate=30] {$>$};
\draw (-4.76,4.2) node {$2$};
\draw (-3.24,3.85) node[rotate=-30] {$>$};
\draw (-3.24,4.2) node {$2$};
\draw[red]  (-2.16,2.77) node[rotate=-60] {$>$};
\draw (-1.8,2.77) node {$2$};
\draw (-2.16,1.23) node[rotate=-120] {$>$};
\draw (-1.8,1.23) node {$2$};
\draw (-3.24,0.15) node[rotate=-150] {$>$};
\draw (-3.24,0.5) node {$2$};

\node (V'7)at(-2.43,-0.75)  {$\bullet$} ;
\node (V'6)at(-5.94,-2.44)  {$\bullet$} ;
\node (V'5)at(-2.58,3.41)  {$\bullet$} ;
\node (V'4)at(-3.13,-3.8)  {$\bullet$} ;
\node (V'3)at(-4.86,-3.8)  {$\bullet$} ;
\node (V'2)at(-2.05,-2.44)  {$\bullet$} ;
\node (V'1)at(-5.57,-0.75)  {$\bullet$} ;

\draw[red] (-4.86,-0.2) node[rotate=210] {$>$};
\draw (-4.7,-0.6) node {$1$};
\draw (-5.94,-1.55) node[rotate=-103] {$>$};
\draw (-6.3,-1.55) node {$1$};
\draw[red] (-5.56,-3.25) node[rotate=-50] {$>$};
\draw (-5.9,-3.4) node {$1$};
\draw (-4,-4) node[rotate=0] {$>$};
\draw (-4,-4.4) node {$1$};
\draw[red] (-2.43,-3.24) node[rotate=50] {$>$};
\draw (-2.1,-3.4) node {$1$};
\draw (-2.05,-1.55) node[rotate=103] {$>$};
\draw (-1.7,-1.55) node {$1$};
\draw (-3.13,-0.2) node[rotate=150] {$>$};
\draw (-3.3,-0.6) node {$1$};

\draw [->,thick] (-1.5,0) -- (-0.5,0) node[midway,above] {$f_{1,2}$} ;

\draw (1.5,3.5) node[below]{$\Gamma_1$};

\draw (1.5,1) circle (1) ;
\draw (1.5,-1) circle (1) ;

\draw[red, ultra thick] (1.5,0) arc (-90:-210:1) ;
\draw[red, ultra thick] (1.5,0) arc (90:270:1) ;

\draw[fill=blue!50] (1.5,0) circle (0.15cm);
\node (V1)at(0.64,1.5)  {$\bullet$} ;
\node (V2)at(2.36,1.5)  {$\bullet$} ;
\node (V3)at(1.5,-2)  {$\bullet$} ;

\draw (1.5,2) node {$>$};
\draw (1.5,2.4) node {$5$};
\draw[red] (0.64,0.5) node[rotate=120] {$>$};
\draw (0.25,0.5) node {$5$};
\draw (2.36,0.5) node[rotate=60] {$<$};
\draw (2.75,0.5) node {$5$};

\draw[red] (0.5,-1) node[rotate=90] {$<$};
\draw (0.1,-1) node {$4$};
\draw (2.5,-1) node[rotate=90] {$>$};
\draw (2.9,-1) node {$4$};

\draw [->,thick] (4,0) -- (5,0) node[midway,above] {$f_{0,1}$} ;

\draw (7,2.5) node[below]{$\Gamma_0$};

\draw (7,0.5) circle (0.5) ;
\draw[red, ultra thick] (7,-0.5) circle (0.5) ;

\draw[fill=blue!50] (7,0) circle (0.15cm);
\draw (7,1) node {$>$};
\draw (7,1.4) node {$14$};
\draw[red] (7,-1) node {$>$};
\draw (7,-1.4) node {$9$};


\node[draw,circle,fill=blue!5,ultra thin,dashed,minimum width=2.5cm] (AA)at(-4,-7) {};
\node[draw,circle,white,minimum width=4cm] (BB)at(-4,-7) {};
\node (AAA)at(-4,-7) {};
\draw[red, ultra thick] (AA.north west) to node[pos=0.5,rotate=-45] {$<$} node[black, pos=1.3] {$2$} (BB.north west); 
\draw (AA.north west) node{$\bullet$} ;
\draw (AA.north east) to node[pos=0.5,rotate=45] {$<$} node[ pos=1.3] {$2$} (BB.north east); 
\draw (AA.north east) node{$\bullet$} ; 
\draw[red, ultra thick]  (AA.south west) to node[pos=0.5,rotate=45] {$<$} node[black, pos=1.3] {$1$} (BB.south west); 
\draw (AA.south west) node{$\bullet$} ; 
\draw (AA.south east) to node[pos=0.5,rotate=-45] {$<$} node[ pos=1.3] {$1$} (BB.south east); 
\draw (AA.south east) node{$\bullet$} ; 

\draw[blue] (AA.north west) to[bend right] node[midway,above] {$1$} node[midway] {$<$} (AA.north east) ;
\draw[blue] (AA.north west) -- node[pos=1.4,right] {$\;1$} node[pos=1.7,rotate=-45] {$<$} (AAA) -- (AA.south east);
\draw[blue] (AA.north east) to node[pos=0.6,left] {$1\;\;$} node[pos=0.75,rotate=45] {$<$} (AA.south west);

\node[draw,circle,fill=blue!5,ultra thin,dashed,minimum width=2.5cm] (CC)at(1.5,-7) {};
\node[draw,circle,white,minimum width=4cm] (DD)at(1.5,-7) {};
\node (CCC)at(1.5,-7) {};
\draw[red, ultra thick] (CC.north west) to node[pos=0.5,rotate=-45] {$<$} node[black, pos=1.3] {$5$} (DD.north west); 
\draw (CC.north west) node{$\bullet$} ; 
\draw (CC.north east) to node[pos=0.5,rotate=45] {$<$} node[ pos=1.3] {$5$} (DD.north east); 
\draw (CC.north east) node{$\bullet$} ; ;
\draw[red, ultra thick]  (CC.south west) to node[pos=0.5,rotate=45] {$<$} node[black, pos=1.3] {$4$} (DD.south west); 
\draw (CC.south west) node{$\bullet$} ; 
\draw (CC.south east) to node[pos=0.5,rotate=-45] {$<$} node[ pos=1.3] {$4$} (DD.south east); 
\draw (CC.south east) node{$\bullet$} ; 

\draw[blue] (CC.north west) to[bend right] node[midway,above] {$2$} node[midway] {$<$} (CC.north east) ;
\draw[blue] (CC.south east) to[bend right] node[midway,below] {$1$} node[midway] {$<$} (CC.south west);
\draw[blue] (CC.north west) -- node[pos=1.4,right] {$\;3$} node[pos=1.7,rotate=-45] {$<$} (CCC) -- (DD.south east);
\draw[blue] (CC.north east) to node[pos=0.6,left] {$3\;\;$} node[pos=0.75,rotate=45] {$<$} (CC.south west);

\node[draw,circle,fill=blue!5,ultra thin,dashed,minimum width=2.5cm] (EE)at(7,-7) {};
\node[draw,circle,white,minimum width=4cm] (FF)at(7,-7) {};
\node (EEE)at(7,-7) {};
\draw (EE.north west) to node[pos=0.5,rotate=-45] {$<$} node[black, pos=1.3] {$14$} (FF.north west); 
\draw (EE.north west) node{$\bullet$} ;
\draw (EE.north east) to node[pos=0.5,rotate=45] {$<$} node[ pos=1.3] {$14$} (FF.north east); 
\draw (EE.north east) node{$\bullet$} ; 
\draw[red, ultra thick]  (EE.south west) to node[pos=0.5,rotate=45] {$<$} node[black, pos=1.3] {$9$} (FF.south west); 
\draw (EE.south west) node{$\bullet$} ; 
\draw[red, ultra thick]  (EE.south east) to node[pos=0.5,rotate=-45] {$<$} node[black, pos=1.3] {$9$} (FF.south east); 
\draw (EE.south east) node{$\bullet$} ; 

\draw[blue] (EE.north west) to[bend right] node[midway,above] {$5$} node[midway] {$<$} (EE.north east) ;
\draw[blue] (EE.north west) -- node[pos=1.4,right] {$\;9$} node[pos=1.7,rotate=-45] {$<$} (EEE) -- (EE.south east);
\draw[blue] (EE.north east) to node[pos=0.6,left] {$9\;\;$} node[pos=0.75,rotate=45] {$<$} (EE.south west);

\end{tikzpicture}
}

\begin{enumerate}
\item[]
Above the lowest three levels $\Gamma_2 \to \Gamma_1 \to \Gamma_0$ of some graph tower $\bvec \Gamma$ are represented.
For each of the level graphs $\Gamma_n$ (with $n \in \{0, 1 ,2\}$) we chose a weight function $\omega_n$ on the edges of $\Gamma_n$: the weight of each edge (for simplicity chosen to be an integer) is stated next to the edge. These choices are arbitrary, except that they must verify both, the switch conditions and the compatibility with respect to passing to a lower level.

If for each of the above pictured graphs $\Gamma_n$ we denote the top edge (in long-edge dialect) by $a_n$ and the bottom edge by $b_n$, the level transition maps are given by $f_{1,2}: a_2 \mapsto a_1 b_1 a_1, \, b_2 \mapsto b_1 b_1 a_1$ and $f_{0,1}: a_1 \mapsto b_0 a_0 a_0, \, b_1 \mapsto b_0 a_0$. In the figure, in addition, the subdivision of the long edges into short-edge dialect is indicated, and the color code illustrates the composed map $f_{0,1} \circ f_{1,2}$.

As in Figure \ref{fig:blow up}, below each graph $\Gamma_n$ the zoom-in of the blown-up local vertex graph is given for the sole vertex in each of the three graphs. The reader verifies easily that both, the switch conditions and the compatibility equalities are satisfied by the given weights on the local and non-local edges.

The given data suffice to compute for sufficiently short paths $\gamma$ the Kolmogorov function which is associated to any weight tower $\bvec \omega$ that agrees on the lowest three levels with the given weight functions. For example, the path $\gamma = a_0 b_0 a_0$ admits 4 distinct lifts to $\Gamma_1$, each crossing over a local edge, thus giving rise to 
$\mu_{\Gamma}^{\tiny \bvec \omega}(a_0 b_0 a_0) = 2 + 3 + 3 + 1 = 9$. The alternative count on $\Gamma_2$ gives $2 + 2 +1 + 1$ for the lifts of $\gamma$ contained in the long edges, and $1+1+1+0$ for the lifts crossing over a local edge. Due to the compatibility of the chosen weights these sums add up to the same value 
as before.

The graph $\Gamma_2$ suffices to compute the Kolmogorov function value for any path of length $\leq 8$: for example $\gamma'= b_0 a_0 a_0 b_0 a_0 b_0$ admits a lift inside $a_2$ and two lifts crossing over local edges, adding up to $\mu_{\Gamma}^{\tiny \bvec \omega}(b_0 a_0 a_0 b_0 a_0 b_0) = 2 + 1 + 0 = 3$.

\end{enumerate}

\caption{
}
\label{fig:tower}

\end{figure}

\begin{defn}
\label{weights+2}
Let
$\bvec \Gamma = ((\Gamma_n)_{n \in \N \cup \{0\}} , (f_{m,n})_{0 \leq m\leq n})$ be a graph tower as in Convention \ref{ft-short}, and let $((\hat\Gamma_n)_{n \in \N \cup \{0\}} , (\hat f_{m,n})_{0 \leq m\leq n})$ be the associated blow-up graph tower (see Remark \ref{gt-translation}).

A {\em tower of weight functions} 
(or simply a {\em weight tower}) $\bvec \omega$ on $\bvec \Gamma = ((\Gamma_n)_{n \in \N \cup \{0\}} , (f_{m,n})_{0 \leq m\leq n})$ is a family of weight functions $\omega_n: \Edges^\pm(\Gamma_n) \to \R_{\geq 0}$ which is induced by a family of weight functions $\hat\omega_n: \Edges^\pm(\hat\Gamma_n) \to  \R_{\geq 0}$ as in part (1) above. The functions $\omega_n$ must furthermore satisfy for all integers $n \geq m \geq 0$ and any edge $e \in \Edges(\Gamma_m)$ the following {\em compatibility condition}: 
\begin{equation}
\label{compatibility-equation}
\omega_m(e) = \underset{\{ e_i \in {\rm Edges}^\pm(\Gamma_n) \, \mid \, f_{m,n}(e_i) = e\}} {\sum}\omega_n(e_i)
\end{equation}
Similarly, for any local edge $\epsilon$ of $\hat \Gamma_m$ one has:
\begin{equation}
\label{local-compatibility}
\hat\omega_m(\epsilon) = 
\underset{\{\epsilon_k \in \cal E^\pm(\hat \Gamma_n) \, \mid \, \hat f_{m,n}(\epsilon_k) = \epsilon\}}{\sum}\hat \omega_n(\epsilon_k)
\end{equation}
\end{defn}

\begin{rem}
\label{equal-weights1}
(1)
From the switch conditions (\ref{switch-conditions}) and the compatibility conditions (\ref{compatibility-equation}) and (\ref{local-compatibility}) together it follows directly 
(see Remark \ref{gt-translation} (3))
that every illegal local edge $\epsilon_i$ 
at any vertex of any level graph $\Gamma_n$ must have weight $\hat \omega_n(\epsilon_i) = 0$. Indeed, any such $\epsilon_i$ is mapped by some $f_{m,n}$ to a single local vertex, and as a result, if the compatibility conditions for $f_{m,n}$ are valid, then $\hat \omega_n(\epsilon_i) \neq 0$ would imply that the switch conditions for $\hat \omega_m$ at this local vertex fail, assuming that for $\hat \omega_n$ they are valid.

\smallskip
\noindent
(2)
Similarly, it follows directly from the above definitions, that any weight tower $\bvec \omega$ on an expanding 
graph tower $\bvec \Gamma$ is indeed carried by the maximal non-pathological subtower (in the sense of Remark \ref{pathological}) of $\bvec \Gamma$, meaning that all edges outside this subtower have weight 0.
\end{rem}

\begin{rem}
\label{equal-weights2}
We observe that any weight function $\omega_n$ on a level graph $\Gamma_n$ induces a weight function $\omega_n^*$ on the relative long-edge dialect level graph $\Gamma_n^*$ associated to $\Gamma_n$, with the property $\omega_n^*(e^*) = \omega_n(e_i)$ for any long edge $e^*$ of $\Gamma_n^*$, and any edge $e_i$ of $\Gamma_n$ which arises from subdividing $e^*$. This is a consequence of the fact that at any subdivision vertex $v_i$ on $e^*$, say equal to the terminal vertex of $e_{i-1}$ and the initial vertex of $e_i$, the local vertex graph $\Gamma(v_i)$ consists only of a single local edge $\epsilon_i$, so that the switch conditions give:
$$\omega_n(e_{i-1}) = \hat\omega_n(\hat e_{i-1}) = \hat\omega_n(\epsilon_i) = \hat\omega_n(\bar \epsilon_i) = \hat\omega_n(\bar{\hat e_{i}}) = \hat\omega_n(\hat e_{i}) = \omega_n(e_{i})$$
As a consequence, we see that in terms of weight functions the local edges at valence 2 vertices of $\Gamma_n$ do not really play any important role.

However, one should keep in mind that, in the compatibility condition (\ref{local-compatibility}) for the local edges, 
for the sum on the right hand side, the summation has to be taken over {\em all} local edges $\epsilon_k$ that are mapped by $\hat f_{m,n}$ to $\epsilon$, which includes also the local edge of the local vertex graph of any blown-up valence 2 vertex.
\end{rem}

\begin{rem}
\label{equal-weights3}
The converse of what has been pointed out in Remark \ref{equal-weights2} is also true:  We first recall that Definition \ref{weights+1} is valid for any graph, so that we can in particular define weight functions $\omega^*_n$ on the level graphs $\Gamma_n^*$ of a graph tower given in long-edge dialect. We then pass to the associated short-edge graph tower, and obtain canonically a weight function $\omega_n$ on its level graphs $\Gamma_n$,
where $\omega^*_n$ and $\omega_n$ are related as spelled out in Remark \ref{equal-weights2}. 

The weights $\omega^*_n$ are now said to define a weight tower $\bvec \omega^*$ if the $\omega_n$ satisfy Definition \ref{weights+2}. We thus conclude that for any graph tower $\bvec \Gamma$ the canonical transition between short-edge dialect and relative long-edge dialect is mirrored by a canonical transition between ``short-edge'' weight towers $\bvec \omega = (\omega_n)_{n \in \N \cup\{0\}}$ on $\bvec \Gamma$ and ``long edge'' weight towers $\bvec \omega^* = (\omega^*_n)_{n \in \N \cup\{0\}}$ on $\bvec \Gamma^*$.
\end{rem}

\begin{convention}
\label{no-hat}
For simplicity, since no confusion is to be feared, we will from now on drop the hat of $\hat \omega_m$ and denote the weight of any local edge $\epsilon$ of any level graph $\Gamma_m$ simply by $\omega_m(\epsilon)$.
\end{convention}

The rest of this section is devoted to performing the first step of what is needed in the next section for the construction of an invariant measure from a given weight tower. For this purpose we first need to extend the definition of weights of edges to weights of edge paths $\gamma$ in a level graph $\Gamma_n$ which are ``sufficiently short''. In some sense this definition is the crucial point of our whole paper.

\smallskip

As a consequence of Remark \ref{equal-weights2} we observe for any reduced edge path $\gamma = e_1 e_2 \ldots e_q$ in $\Gamma_n$ that, if $\gamma$
is entirely contained in some edge $e^*$ from the associated relative 
long-edge dialect graph $\Gamma^*_n$, all edges $e_{i}$ traversed by $\gamma$ have the same weight. Thus the definition 
\begin{equation}
\label{5.4+}
\omega_n(\gamma) := \omega_n(e_i)
\end{equation}
is independent of the edge $e_i$ traversed by the path $\gamma$.

On the other hand, if $\gamma$ traverses any {intrinsic} vertex $v$ of $\Gamma_n$, i.e. a vertex which is inherited from a vertex of $\Gamma_n^*$, then the local edge $\epsilon$ traversed by $\gamma$ at $v$ and the two edges $e$ and $e'$ of $\Gamma_n$ which are adjacent to $\epsilon$ on $\gamma$ satisfy
\begin{equation}
\label{weight-inequality}
\omega_n(\epsilon) \leq \omega_n(e), \quad {\rm and} \quad \omega_n(\epsilon) \leq \omega_n(e')\, ,
\end{equation}
and these inequalities may well be strict. For such $\gamma$, 
if $\epsilon$ is the only intrinsic local edge traversed by $\gamma$, we set:
\begin{equation}
\label{5.5+}
\omega_n(\gamma) := \omega_n(\epsilon)
\end{equation}

For the following definition we recall from Definition \ref{expanding} that for any graph tower $\bvec \Gamma = ((\Gamma_n)_{n \in \N \cup \{0\}} , (f_{m,n})_{0 \leq m\leq n})$ in short-edge dialect (recall Convention \ref{ft-short}) and any edge path $\gamma$ in some level graph $\Gamma_n$ the condition
$$|\gamma| \leq {\mile} _{\tiny\bvec \Gamma}(n)+1$$
implies that $\gamma$ crosses over at most one intrinsic vertex of $\Gamma_n$.
This definition and the subsequent proposition are illustrated by Figure \ref{fig:tower}.

\begin{defn}
\label{weights-of-paths}
Let $\bvec \omega = (\omega_n: \Gamma_n \to \R_{\geq 0})_{n\in \N \cup \{0\}}$ be a weight tower on an expanding graph tower $\bvec \Gamma = ((\Gamma_n)_{n \in \N \cup \{0\}} , (f_{m,n})_{0 \leq m\leq n})$.
\begin{enumerate}
\item
For any non-trivial finite reduced path $\gamma$ in any level graph $\Gamma_n$, with length $|\gamma| \leq {\mile} _{\tiny\bvec \Gamma}(n)+1$, 
the {\em weight} $\omega_n(\gamma)$ is defined by the equalities (\ref{5.4+}) and (\ref{5.5+}).
\item
For any finite reduced path $\gamma$ in $\Gamma_0$ and any level graph $\Gamma_n$ 
we denote by $\cal E_n(\gamma)$ the set of all edge paths $\gamma_i$ in $\Gamma_n$ with $f_n(\gamma_i) = \gamma$.
\item
For any finite reduced path $\gamma$ in $\Gamma_0$ and any integer $n \geq 0$ which satisfies 
$|\gamma| \leq {\mile} _{\tiny\bvec \Gamma}(n)+1$ 
we define the {\em $n$-weight} $\omega_n(\gamma)$ of $\gamma$ through:
$$\omega_n(\gamma) := \sum_{\gamma_i \in \cal E_n(\gamma)} \omega_n(\gamma_i)$$
\end{enumerate}
\end{defn}

\begin{prop}
\label{weights-via-lifts}
Let $\bvec \Gamma = ((\Gamma_n)_{n \in \N \cup \{0\}} , (f_{m,n})_{0 \leq m\leq n})$  be an expanding graph tower, and let $\bvec \omega = (\omega_n: \Gamma_n \to \R_{\geq 0})_{n\in \N \cup \{0\}}$ be a tower of weight functions on $\bvec \Gamma$.

For any finite reduced path $\gamma$ in $\Gamma_0$ the $n$-weights $\omega_n(\gamma)$ are independent of $n$, for any $n$ with $|\gamma| \leq {\mile} _{\tiny\bvec \Gamma}(n)+1$.

In particular, from $\bvec \omega$ one obtains a well defined function $\mu^{\tiny \bvec \omega}_\Gamma: \cal P(\Gamma_0) \to \R_{\geq 0}$, given, for any sufficiently large integer $n \geq 0$, 
by:
\begin{equation}
\label{5.4}
\mu^{\tiny \bvec \omega}_\Gamma(\gamma) := 
\sum_{\gamma_i \in \cal E_n(\gamma)} \omega_n(\gamma_i)
\end{equation}
\end{prop}

\begin{proof}
We first 
recall that the assumption $|\gamma| \leq \mile_{\tiny\bvec \Gamma}(n)+1$ implies that any $\gamma_i \in \cal E_n(\gamma)$ traverses at most one intrinsic vertex of $\Gamma_n$.
We now choose any integer $k \geq n$ and consider any legal path $\gamma'_j$ in $\cal E_k(\gamma)$. We note that  
the path $\gamma_i := f_{n,k}(\gamma'_j)$ belongs to $\cal E_n(\gamma)$. Since the paths in $\cal E_k(\gamma)$ are partitioned according to their $f_{n, k}$-image paths, in order to show that the sum in equation (\ref{5.4}) over $\cal E_k(\gamma)$ gives the same value as the sum over $\cal E_n(\gamma)$, it suffices to show for every $\gamma_i \in \cal E_n(\gamma)$ the following {claim}:
\begin{equation}
\label{claim}
\omega_n(\gamma_i) = \sum_{\{\gamma'_j \in \cal E_k(\gamma) \,\mid\, f_{n,k}(\gamma'_j) = \gamma_i\}} \omega_k(\gamma'_j)
\end{equation}

We distinguish two cases:

\smallskip
\noindent
(1)
If $\gamma_i$ does not cross over an intrinsic vertex of $\Gamma_n$, then any of the $\gamma'_j$ with $f_{n,k}(\gamma'_j) = \gamma_i$ can not cross either over any intrinsic vertex of $\Gamma_k$, as the level maps in any graph tower 
map intrinsic vertices to intrinsic vertices (see Definition-Remark \ref{rel-long}).  In particular, it follows that every short edge $e'$ of 
$\Gamma_k$ 
which is mapped to any short edge $e$ on the path $\gamma_i$, must be part of a unique path $\gamma'_j$ with $f_{n,k}(\gamma'_j) = \gamma_i$.

Furthermore, from (\ref{5.4+}) we obtain $\omega_k(\gamma'_j) = \omega_k(e')$ for any short edge $e'$ contained in $\gamma'_j$, and since we similarly have $\omega_n(\gamma_i) = \omega_n(e)$ for any short edge $e$ contained in $\gamma_i$, the above claim (\ref{claim}) follows now directly from the compatibility condition (\ref{compatibility-equation}).

\smallskip
\noindent
(2)
In the case that $\gamma_i$ crosses over a single intrinsic vertex $v$ of $\Gamma_n$, then we consider the local edge $\epsilon$ at $v$ traversed 
by $\gamma_i$, and observe that $\omega_n(\gamma_i) = \omega_n(\epsilon)$ holds, by the above definition (\ref{5.5+}) of $\omega_n(\gamma_i)$. Since the $f_{n, k}$-image of any intrinsic vertex in $\Gamma_k$ is again an intrinsic vertex of $\Gamma_n$ (see Definition-Remark \ref{rel-long}), we deduce from the hypothesis $|\gamma| \leq \mile_{\tiny\bvec \Gamma}(n)+1$ that for any preimage vertex $v' \in \Gamma_k$ of $v$ and any local edge $\epsilon'$ at $v'$ with $f_{n,k}(\epsilon') = \epsilon$ there is precisely one edge path $\gamma'_j$ crossing over $\epsilon'$ such that $f_{n,k}(\gamma'_j) = \gamma_i$.  Conversely, 
for any path $\gamma'_j$ in $\Gamma_k$ with $f_{n,k}(\gamma'_j) = \gamma_i$ there must be a preimage vertex $v'$ of $v$ and a local edge $\epsilon'$ at $v'$ with $f_{n,k}(\epsilon') = \epsilon$ such that $\gamma'_j$ crosses over $\epsilon'$.

From the previous paragraph we know that the only vertex crossed over by $\gamma'_j$ which is possibly intrinsic must be the preimage vertex $v'$ of $v$. We deduce from the definition of the weights of paths in (\ref{5.4+}) and in (\ref{5.5+}) that in either case, whether $v'$ is intrinsic or not, we have $\omega_k(\gamma'_j)= \omega_k(\epsilon')$. As in case (1), the ``local'' compatibility conditions (\ref{local-compatibility}) give $\omega_n(\epsilon) = \underset{\{\epsilon' \,\mid\, f_{n,k}(\epsilon') = \epsilon\}}{\sum} \omega_k(\epsilon')$, which shows the above claim (\ref{claim}) also in case (2).
\end{proof}

\begin{rem}
\label{almost-weights--}
For any expanding graph tower $\bvec \Gamma = ((\Gamma_n)_{n \in \N \cup \{0\}} , (f_{m,n})_{0 \leq m\leq n})$ and any weight tower $\bvec{\omega} = (\vec \omega_n) _{n \in \N \cup \{0\}}$ on $\bvec \Gamma$ one has:
$$\lim_{n \to \infty} \max\{\omega_n(e) \mid e \in {\rm Edges}^\pm(\Gamma_n)\} = 0$$
This can be derived directly from the above definitions; alternatively it follows from Remark \ref{several}, where also a stronger statement has been exhibited.
\end{rem}

 
\section{Invariant measures via weight towers}
\label{measures-through-weights}
In this section we will explain how to derive from a weight tower $\bvec \omega$ on an expanding graph tower $\bvec \Gamma$ an invariant measure on the used lamination $\Lused$. The central tool for this purpose are Kolmogorov functions. We start the section by briefly recalling its basics.

\smallskip

For any graph $\Gamma$ we denote as before by $\cal P(\Gamma)$ the set of finite reduced edge paths $\gamma$ in $\Gamma$.

\begin{defn}
\label{K-function}
A a non-negative function $\mu_\Gamma: \cal P(\Gamma) \to\R_{\geq 0}$ is called a {\em Kolmogorov function} on $\Gamma$ if it satisfies, for every $\gamma = e_1 e_2 \ldots e_q \in \cal P(\Gamma)$, the equality
\begin{equation}
\label{Kolmo-1}
\mu_\Gamma(\gamma) = \mu_\Gamma(\bar \gamma)
\end{equation}
as well as the {\em Kirchhoff rules}:
\begin{equation}
\label{Kirchhoff}
\mu_\Gamma(\gamma)
= \sum_{\substack{
e_0 \in {\rm Edges}^\pm(\Gamma) \smallsetminus \{\bar e_1\} \\ 
\tau(e_0) = \tau(\bar \gamma)}} 
\mu_\Gamma(e_0 \gamma) 
= \sum_{\substack{
e_{q+1} \in {\rm Edges}^\pm(\Gamma) \smallsetminus \{\bar e_q\} \\
\tau(\gamma) = \tau(\bar e_{q+1})}} 
\mu_\Gamma(\gamma e_{q+1})
\end{equation}
\end{defn}

Every Kolmogorov function $\mu_\Gamma$ on $\Gamma$ {\em generates} a language $\cal L(\mu_\Gamma) \subseteq\cal P(\Gamma)$, sometimes called the {\em support} of $\mu_\Gamma$, which is given by all reduced paths $\gamma$ with $\mu_\Gamma(\gamma) > 0$. A language $\cal L \subset \cal P(\Gamma)$ {\em carries} a Kolmogorov function $\mu_\Gamma$ if $\mu_\Gamma(\gamma) = 0$ for all $\gamma \in \cal P(\Gamma) \smallsetminus \cal L$. It follows directly from the above definition of a Kolmogorov function that $\cal L(\mu_\Gamma)$ is laminary (see Definition \ref{language+} (3)(d)), unless $\mu_\Gamma$ is the zero-function.

\smallskip

Recall from section \ref{sec:invariant-measures} that
an {invariant measure} $\mu_\Sigma$ 
for $\Gamma$ is a finite Borel measure on $\Sigma(\Gamma)$ which is invariant under shift and inversion.
It canonically defines a function $\mu_\Gamma$ on the set $\cal P(\Gamma)$ of all finite reduced edge paths $\gamma$ in $\Gamma$, given by setting:
$$\mu_\Gamma(\gamma) := \mu_\Sigma(C_\gamma)$$
Here $C_\gamma \subseteq \Sigma(\Gamma)$ denotes as before the {cylinder} associated to the edge path $\gamma = e_1 \ldots e_r$ in $\Gamma$, defined as the set of all biinfinite reduced paths $\ldots e'_{n-1} e'_n e'_{n+1} \ldots$ which satisfy $e' _1 = e_1, \ldots, e'_r = e_r$.

The following consequence of Carath\'eodory's extension theorem is 
classical (see 
Theorem 3.2 of \cite{Ka1}):

\begin{prop}
\label{K-function-measure}
Let $\Gamma$ be any finite graph, let $\cal P(\Gamma)$ be the set of finite reduced edge paths in $\Gamma$, and let $\Sigma(\Gamma)$ denote the set of biinfinite reduced edge paths in $\Gamma$.
\begin{enumerate}
\item
For any invariant measure $\mu_\Sigma$ on $\Sigma(\Gamma)$ the associated function $\mu_\Gamma: \cal P(\Gamma) \to \R_{\geq 0}$ is a Kolmogorov function.
\item
Conversely, every Kolmogorov function $\mu_\Gamma$ on $\cal P(\Gamma)$ defines an invariant measure $\mu_\Sigma$ for $\Gamma$
which satisfies $\mu_\Sigma(C_\gamma) = \mu_\Gamma(\gamma)$ for any $\gamma \in \cal P(\Gamma)$.
\qed
\end{enumerate}
\end{prop}

It follows from this proposition that the passage back and forth between invariant measures and Kolmogorov functions is canonical. This canonical passage translates into a canonical passage between their supports $\cal L(\mu_\Gamma)$ and $\Lsig(\mu_\Sigma)$ (see (\ref{support-measure})), in the sense that one has, using the notation from Definition \ref{language+} (4) and (5),
\begin{equation}
\label{support-Kolmo}
\cal L(\mu_\Gamma) = \cal L(L^\Sigma(\mu_\Sigma))
\quad \text{and} \quad
L^\Sigma(\mu_\Sigma) = L^\Sigma(\cal L(\mu_\Gamma)) \, .
\end{equation}
Again we must exclude here the case where $\mu_\Gamma$ and $\mu_\Sigma$ are the zero-functions, as the empty set doesn't count formally as symbolic lamination.

\begin{rem}
\label{total-weight}
Let $\Gamma$ be a finite graph equipped with a Kolmogorov function $\mu_\Gamma$, and let $\mu_\Sigma$ be the associated invariant measure for $\Gamma$.
For any integer $n \geq 0$ denote by $\cal L_n \subset \cal P(\Gamma)$ the set of all reduced edge paths in $\Gamma$ of length $n$. 
We then derive immediately from Definition \ref{K-function} 
that the sum
$$\sum_{\gamma \in \cal L_n} \mu_\Gamma(\gamma)$$
is independent of $n$. It is called the {\em total weight} of the Kolmogorov function $\mu_\Gamma$, and it is equal to the total measure $\mu_\Sigma(\Sigma(\Gamma))$ of $\mu_\Sigma$.
\end{rem}

\smallskip

We now return to weighted graph towers as defined in the previous section:
\begin{prop}
\label{currents-weights}
Let $\bvec \Gamma = ((\Gamma_n)_{n \in \N \cup \{0\}} , (f_{m,n})_{0 \leq m\leq n})$  be an expanding graph tower,
and let $\bvec \omega = (\omega_n: \Gamma_n \to \R_{\geq 0})_{n\in \N \cup \{0\}}$ be a tower of weight functions on $\bvec \Gamma$.
Then the function $\mu_\Gamma^{\tiny \bvec \omega}: \cal P(\Gamma_0) \to \R_{\geq 0}$ defined in Proposition \ref{weights-via-lifts} is a Kolmogorov function on $\Gamma_0$. 

In particular, the weight tower $\bvec \omega$ also defines an associated invariant measure $\mu_\Sigma^{\tiny \bvec \omega}$ for $\Gamma_0$.
\end{prop}

\begin{proof}
In order to 
to verify the Kirchhoff conditions (\ref{Kirchhoff}) for the function $\mu^{\tiny \bvec \omega}_\Gamma$ we consider, according to Proposition \ref{weights-via-lifts}, for any finite reduced path $\gamma$ in $\Gamma_0$ of length $|\gamma| = s$, any level graph $\Gamma_n$ where the minimal length of long edges satisfies ${\mile} _{\tiny\bvec \Gamma}(n) \geq s$. 
Then for any legal path $\gamma'$ in $\Gamma_n$ with $f_n(\gamma') = \gamma$ the switch conditions (\ref{switch-conditions}) show directly, 
for both of the cases considered in (\ref{5.4+}) and (\ref{5.5+}), that one has:
$$\omega_n(\gamma') = \sum_{\gamma_i \in I(\gamma')} \omega_n(\gamma_i) = \sum_{\gamma_j \in T(\gamma')} \omega_n(\gamma_j)$$
where $I(\gamma')$ is the set of all reduced paths $\gamma_i$ of length $s+1$ in $\Gamma_n$ which have $\gamma'$ as initial subpath, 
and $T(\gamma')$ is the set of all reduced paths $\gamma_j$ of length $s+1$ in $\Gamma_n$ which have $\gamma'$ as terminal subpath. The equality (\ref{Kirchhoff}) for the function $\mu^{\tiny \bvec \omega}_\Gamma$ is then a direct consequence of the definition of $\mu^{\tiny \bvec \omega}_\Gamma$ in Proposition \ref{weights-via-lifts}.  

Furthermore, the condition (\ref{Kolmo-1}) is a direct consequence of the equality $\omega(\bar e) = \omega(e)$ from Definition \ref{weights+1}.

Since $\mu_\Gamma^{\tiny \bvec \omega}$ satisfies the condition (\ref{Kolmo-1}) and (\ref{Kirchhoff}), it is a Kolmogorov function on $\Gamma_0$, and hence, by Proposition \ref{K-function-measure} (2), 
it defines an associated invariant measure $\mu_\Sigma^{\tiny \bvec \omega}$ for $\Gamma_0$.
\end{proof}

From Remark \ref{equal-weights1} we know 
$\omega_n(\epsilon_i) = 0$ for any illegal local edge at any level graph $\Gamma_n$, which implies 
directly that the support of $\mu_\Sigma^{\tiny \bvec \omega}$ 
is contained in $L^\Sigma_{legal}(\bvec \Gamma)$.
However, we can actually do a little better:

\begin{prop}
\label{almost-weights-}
Let $\bvec \Gamma = ((\Gamma_n)_{n \in \N \cup \{0\}} , (f_{m,n})_{0 \leq m\leq n})$  be an expanding graph tower, let $\bvec{\omega} = (\vec \omega_n) _{n \in \N \cup \{0\}}$ be a weight tower on $\bvec \Gamma$, and let $\mu_\Sigma^{\tiny \bvec \omega}$ be the associated invariant measure on $\Gamma_0$.

Then the support $\Lsig(\mu_\Sigma^{\tiny \bvec \omega})$ of $\mu_\Sigma^{\tiny \bvec \omega}$ is contained in the used lamination $\Lused$ defined by $\bvec \Gamma$.
\end{prop}

\begin{proof}
Going back to the notation of Definition \ref{weights-of-paths} (2) and Proposition \ref{weights-via-lifts}, 
for every reduced edge path $\gamma$ in 
$\Gamma_0$ and any level graph $\Gamma_n$ with ${\mile} _{\tiny\bvec \Gamma}(n) \geq  4 \, |\gamma|$ 
we introduce the following notation:
\begin{enumerate}
\item
Let $\cal E^{(1)}_n(\gamma)$ be the set of those edge paths $\gamma_i$ from $\cal E_n(\gamma)$ which cross over an intrinsic vertex. 
\item
Let $\cal E^{(2)}_n(\gamma)$ be the set of those $\gamma_i \in \cal E_n(\gamma)$ which don't cross over an intrinsic vertex, but are adjacent to a path from $\cal E^{(1)}_n(\gamma)$ 
(in the sense that the two paths can be concatenated to give a reduced path).
\item
Let $\cal E^{(3)}_n(\gamma)$ be the set of those $\gamma_i \in \cal E_n(\gamma)$ which don't cross over any intrinsic vertex and are not adjacent to any $\gamma_i \in \cal E_n(\gamma)$ which does. 
\end{enumerate}
The hypothesis ${\mile} _{\tiny\bvec \Gamma}(n) \geq  4 \, |\gamma|$ is used here to ensure that any $\gamma_i \in \cal E^{(2)}_n(\gamma)$ is precisely on one of its two sides adjacent to a path from $\cal E^{(1)}_n(\gamma)$, while on the other side it is adjacent to a path from $\cal E^{(3)}_n(\gamma)$.

\smallskip

According to this set-up, for $j \in \{1, 2, 3\}$ we define the sum
$$
\omega^{(j)}_n(\gamma) := \sum_
{\gamma_i \in \cal E^{(j)}_n(\gamma)}
\omega_n(\gamma_i) \, .
$$
From the disjoint union decomposition $\cal E_n(\gamma)$ = $\cal E^{(1)}_n(\gamma) \disjoint \cal E^{(2)}_n(\gamma) \disjoint \cal E^{(3)}_n(\gamma)$ 
and from the definition of $\mu^{\tiny \bvec \omega}_\Gamma$ through Proposition \ref{weights-via-lifts} 
we thus obtain directly:
$$
\mu^{\tiny \bvec \omega}_\Gamma(\gamma)
= \omega^{(1)}_n(\gamma) + \omega^{(2)}_n(\gamma) + \omega^{(3)}_n(\gamma)
$$

For any integer $m \geq 1$  we now consider the set $\cal L_m(\Gamma_0)$ of all reduced edge paths 
in $\Gamma_0$ 
of length $m$. 
Let $n \geq 0$ be any integer which satisfies ${\mile} _{\tiny\bvec \Gamma}(n) \geq  4m$. 
We then deduce 
from the switch conditions  (\ref{switch-conditions}) 
(as specified in the inequalities (\ref{weight-inequality})) that
\begin{equation}
\label{inequ1}
\sum_{\gamma \in \cal L_m(\Gamma_0)}
\omega^{(1)}_n(\gamma) \leq \sum_{\gamma \in \cal L_m(\Gamma_0)} \omega^{(2)}_n(\gamma)\, .
\end{equation}
(Actually one has $2 \underset{\gamma \in \cal L_m(\Gamma_0)}{\sum}\omega^{(1)}_n(\gamma) \leq \underset{\gamma \in \cal L_m(\Gamma_0)}{\sum} \omega^{(2)}_n(\gamma)\,$, since to any path $\gamma_i \in \cal E^{(1)}_n(\gamma)$ there are two paths from $\cal E^{(2)}_n(\gamma)$ adjacent to $\gamma_i$, namely one on each side of $\gamma_i$.)

Furthermore we observe, for any long edge $e^*$ in $\Gamma^*_n$, that only the first $m$ and the last $m$ subpaths of length $m$ of $e^*$ 
contribute to $\sum_{\gamma \in \cal L_m(\Gamma_0)}\omega^{(2)}_n(\gamma)$, while all the other such subpaths contribute to $\sum_{\gamma \in \cal L_m(\Gamma_0)}\omega^{(3)}_n(\gamma)$, and all such contribution are equal (as they are equal to $\omega_n(e^*)$). Since from the definition of ${\mile} _{\tiny\bvec \Gamma}(n)$ we have $|e^*| \geq {\mile} _{\tiny\bvec \Gamma}(n)$ for any such long edge $e^*$, we
obtain for any sufficiently large $n$ the following inequality:
\begin{equation}
\label{inequ2}
\sum_{\gamma \in \cal L_m(\Gamma_0)} \omega^{(2)}_n(\gamma) \leq  \frac{2 \, m}{{\mile} _{\tiny\bvec \Gamma}(n) - 2m}\sum_{\gamma \in \cal L_m(\Gamma_0)}\omega^{(3)}_n(\gamma)
\end{equation}
Together with the inequality (\ref{inequ1})
and the summation
$$
\sum_{\gamma \in \cal L_m(\Gamma_0)} \omega^{(1)}_n(\gamma) + \sum_{\gamma \in \cal L_m(\Gamma_0)}\omega^{(2)}_n(\gamma) + \sum_{\gamma \in \cal L_m(\Gamma_0)}\omega^{(3)}_n(\gamma)
=
\sum_{\gamma \in \cal L_m(\Gamma_0)} 
\mu^{\tiny \bvec \omega}_\Gamma(\gamma)
$$
the inequality (\ref{inequ2}) 
gives:
\begin{equation}
\label{last-inequality}
\sum_{\gamma \in \cal L_m(\Gamma_0)}
\omega^{(1)}_n(\gamma)
\leq  \frac{2 \, m}{{\mile} _{\tiny\bvec \Gamma}(n) -2m} \sum_{\gamma \in \cal L_m(\Gamma_0)} 
\mu^{\tiny \bvec \omega}_\Gamma(\gamma)
\end{equation}
We now recall from 
Remark \ref{total-weight}
that for any $m \geq 1$ the sum 
$$\sum_{\gamma \in \cal L_m(\Gamma_0)} \mu_\Gamma^{\tiny \bvec \omega}(\gamma)$$
is equal to the the total measure $\mu_\Sigma^{\tiny \bvec \omega}(\Sigma(\Gamma_0))$ and hence independent of $m$. Thus from the hypothesis that $\bvec \Gamma$ is expanding, i.e. ${\mile} _{\tiny\bvec \Gamma}(n) \to \infty$ for $n \to \infty$, and from inequality (\ref{last-inequality}) we deduce 
for any $m \geq 1$ that 
$$
\lim_{n \to \infty} \sum_{\gamma \in \cal L_m(\Gamma_0)}
\omega^{(1)}_n(\gamma) = 0
\qquad  \text{and thus} \qquad
\lim_{n \to \infty} 
\omega^{(1)}_n(\gamma) = 0
$$
for any $\gamma \in \cal L_m(\Gamma_0)$.
As a consequence we obtain, for any reduced edge path $\gamma$ in $\Gamma_0$ that 
\begin{equation}
\label{last-sum}
\lim_{n \to \infty} (\omega^{(2)}_n(\gamma) + \omega^{(3)}_n(\gamma)) = \mu^{\tiny \bvec \omega}_\Gamma(\gamma) \, .
\end{equation}

From Definition \ref{used-lamination+} (1) we know that for any non-used path $\gamma$ in $\Gamma_0$ 
there exists a level $n_0 \geq 0$ such that for any level graph $\Gamma_n$ 
with $n \geq n_0$ 
all paths $\gamma_i$ in $\Gamma_n$ with $f_n(\gamma_i) = \gamma$ must cross over at least one intrinsic vertex, and hence, for sufficiently large $n$, over precisely one intrinsic vertex, so that any such $\gamma_i$ belongs to $\cal E^{(1)}_n(\gamma)$. 
Hence we have $\cal E^{(2)}_n(\gamma) = \cal E^{(3)}_n(\gamma) = \emptyset$, 
which shows $\omega^{(2)}_n(\gamma) = \omega^{(3)}_n(\gamma) = 0$, so that equality (\ref{last-sum}) proves
$$\mu^{\tiny \bvec \omega}_\Gamma(\gamma) = 0$$
for any non-used path $\gamma$ in $\Gamma_0$.
Hence the claim follows from the second equality of (\ref{support-Kolmo}).
\end{proof}

\begin{rem}
\label{non-used}
From Proposition \ref{almost-weights-} it follows in particular, for any level graph $\Gamma_n$ of an expanding graph tower $\bvec \Gamma$ with weight tower $\bvec \omega  = (\vec \omega_n) _{n \in \N \cup \{0\}}$, that for any {\em non-used} local edge $\epsilon_j$ in $\Gamma_n$, i.e. a local edge that is not crossed over by any used path (or alternatively: ``that is not contained in the $f_{n, k}$-image of any non-local edge of $\hat \Gamma_k$, for any $k \geq n$'') one has: 
$$\omega_n(\epsilon_j) = 0$$
\end{rem}

\begin{rem}
\label{erase-finite-part}
(1)
For any graph tower $\bvec \Gamma$, and for any subset $M \subset \N \cup \{0\}$ which doesn't contain $0$ and has infinite complement $(\N \smallsetminus M) \cup \{0\}$, we consider the {\em telescoped} graph tower $\bvec \Gamma'$ obtained from $\bvec \Gamma$ through erasing all level graphs $\Gamma_n$ with $n \in M$. 
The one derives directly from the definition of the used lamination (see Definition \ref{used-lamination+}) that
$$L^{\tiny \bvec \Gamma'} = \Lused$$

\smallskip
\noindent
(2)
For any weight tower $\bvec \omega$ on $\bvec \Gamma$ and for the associated telescoped weight tower $\bvec \omega'$, obtained correspondingly from $\bvec \omega$, 
one derives directly from Proposition \ref{weights-via-lifts} and Proposition \ref{currents-weights}:
$$\mu_\Sigma^{\tiny \bvec \omega'} = \mu_\Sigma^{\tiny \bvec \omega}
\qquad {\rm and} \qquad
\mu_\Gamma^{\tiny \bvec \omega'} = \mu_\Gamma^{\tiny \bvec \omega}
$$
\end{rem}

 
\section{Weight towers through invariant measures}
\label{measure-to-weight}

The purpose of this section is to show the converse of the previous section, i.e. that every invariant measure on the used lamination of an expanding graph tower $\bvec \Gamma$ is given by some weight tower on $\bvec \Gamma$.
\begin{lem}
\label{pseudo-Birkhoff}
Let $\Gamma$ be a graph, $L^\Sigma \subset \Sigma(\Gamma)$ a symbolic lamination, and let $\mu_\Sigma$ be an invariant measure for $\Gamma$ 
with support in $L^\Sigma$. Denote by $\cal L = \cal L(L^\Sigma)$ the laminary language associated to $L^\Sigma$, and by $\mu_\Gamma$ the Kolmogorov function associated to $\mu_\Sigma$. For any $n \geq 0$ denote by $\cal L_n$ the set of edge paths $\gamma \in \cal L$ of length $|\gamma| = n$. Then one has, for any $\gamma' \in \cal L$ and any integer $n \geq |\gamma'|$ :
\begin{equation}
\label{measure-summation}
\mu_\Gamma(\gamma') = \frac{1}{n - |\gamma'| +1} \sum_{\gamma \in \cal L_n} |\gamma|_{\gamma'} \,\, \mu_\Gamma(\gamma)
\end{equation}
Here $|\gamma|_{\gamma'}$ denotes the number of occurrences of the path $\gamma'$ as subpath in $\gamma$ 
(see Definition \ref{length-defn} (2)).
\end{lem}

\begin{proof}
Fix any 
integer $k$ with $1\leq k \leq n - |\gamma'| +1$, and take first the sum $S_k$ on the right hand side of equality (\ref{measure-summation}), but only over those occurrences 
of $\gamma'$ as a subpath of any of the $\gamma$ where $\gamma'$ starts
at the $k$-th edge of $\gamma$:
$$S_k = \sum_{\cal L'} \mu_\Gamma(\gamma_1 \cdot \gamma' \cdot \gamma_2) \qquad \text{with}
$$
$$\cal L' = \{(\gamma_1, \gamma_2) \in 
\cal L_{k-1}\times \cal L_{n - k + 1 - |\gamma'|} \mid 
\gamma_1 \cdot \gamma' \cdot \gamma_2 \in \cal L_n\}$$
In a second step these partial sums $S_k$ are summed up over all $k$, to get the right hand side of (\ref{measure-summation}) up to the factor $\frac{1}{n - |\gamma'| +1}$:
$$
\sum_{k = 1}^{n - |\gamma'| +1} S_k
=
\sum_{\gamma \in \cal L_n} |\gamma|_{\gamma'} \,\, \mu_\Gamma(\gamma)
$$
Now, for any fixed $k$ we use the Kirchhoff rules (\ref{Kirchhoff}) to see directly that the sum $S_k$ is equal to the measure $\mu_\Gamma(\gamma')$, and since there are precisely $n - |\gamma'| +1$ possibilities for $k$, we obtain the desired equality.
\end{proof}

\begin{prop}
\label{inverse-curr-w}
Let $\bvec \Gamma = ((\Gamma_n)_{n \in \N \cup \{0\}} , (f_{m,n})_{0 \leq m\leq n})$  be an expanding graph tower, 
and let $\mu_\Sigma$ be an invariant measure 
carried by the used symbolic lamination $\Lused$Let $\mu_\Gamma$ be the associated Kolmogorov function.

Then there exists a tower of weight functions $\bvec \omega = (\omega_n: \Gamma_n \to \R_{\geq 0})_{n\in \N \cup \{0\}}$ on $\bvec \Gamma$ which satisfies: 
$$\mu_\Sigma^{\tiny \bvec \omega} = \mu_\Sigma
\qquad {\rm and} \qquad
\mu_\Gamma^{\tiny \bvec \omega} = \mu_\Gamma
$$
\end{prop}

\begin{proof}
Observe first that by definition of the used lamination, for every finite path $\gamma$ in the language defined by used lamination 
$\Lused$, and for any $m \geq 0$, 
there is a reduced path $\gamma_m$ in $\Gamma_m$ with $f_m(\gamma_m) = \gamma$. 

We now fix $m \geq 0$ and consider for any (large) integer $n$, say $n \geq {\mile} _{\tiny\bvec \Gamma}(m)$, 
the subset $\cal L_n \subset \cal L(\Lused)$ of all used paths $\gamma$ of length $n$, and and for each of them we choose as above some lift $\gamma_m$ in $\Gamma_m$. In this choice of the lifts we respect the restriction that the lift $(\bar \gamma)_m$ of the inverse $\bar \gamma$ of any $\gamma \in \cal L_n$ is equal to the inverse $\bar \gamma_m$ of the lift $\gamma_m$ of $\gamma$.
We now define a ``pseudo-weight function'' $\omega_{m,n}$ on $\Gamma_m$, through setting for any edge $e$ of $\Gamma_m$ (and analogously for any local edge $\epsilon$ of $\Gamma_m$):
\begin{equation}
\label{8-2-1}
\omega_{m,n}(e) := \frac{1}{n} \sum_{\gamma \in \cal L_n} |\gamma_m|_e \,\, \mu_\Gamma(\gamma)
\end{equation}
Similarly, for any reduced edge path $\gamma_0$ in $\Gamma_m$ with $|\gamma_0| \leq {\mile} _{\tiny\bvec \Gamma}(m)$ we define:
\begin{equation}
\label{8-2-2}
\omega_{m,n}(\gamma_0) := \frac{1}{n} \sum_{\gamma \in \cal L_n} |\gamma_m|_{\gamma_0} \,\, \mu_\Gamma(\gamma)
\end{equation}

With respect to these definitions we point out the following three observations, where we recall, using the terminology of Definition \ref{weights-of-paths} (2), 
that for any finite reduced path $\gamma'$ in $\Gamma_0$ the condition $|\gamma'| \leq {\mile} _{\tiny\bvec \Gamma}(m)$ implies that any path $\gamma'_i$ in the set $\cal E_m(\gamma')$ 
of lifts of $\gamma'$ to $\Gamma_m$ crosses over at most one intrinsic vertex of $\Gamma_m$. 

\begin{enumerate}
\item[({\bf A})]
For any reduced edge path $\gamma'$ in $\Gamma_0$ of length $|\gamma'| \leq {\mile} _{\tiny\bvec \Gamma}(m)$ 
we observe that there is a canonical bijection between the occurrences of $\gamma'$ as subpath of $\gamma \in \cal L_n$ on one hand, and the occurrences of any of the $\gamma'_i \in \cal E_m(\gamma')$ as subpath of the lift $\gamma_m$ of $\gamma$ on the other. Hence
we obtain 
$$
\sum_{\gamma'_i \in \cal E_m(\gamma')}\omega_{m, n}(\gamma'_i) = 
\sum_{\gamma'_i \in \cal E_m(\gamma')}\big{(}\frac{1}{n} \sum_{\gamma \in \cal L_n} |\gamma_m|_{\gamma'_i} \,\, \mu_\Gamma(\gamma)\big{)}
= 
\frac{1}{n} \sum_{\gamma \in \cal L_n} |\gamma|_{\gamma'} \,\, \mu_\Gamma(\gamma)
$$
and thus, by Lemma \ref{pseudo-Birkhoff}:
$$
\sum_{\gamma'_i \in \cal E_m(\gamma')}\omega_{m, n}(\gamma'_i) = 
\frac{n-|\gamma'|+1}{n} \mu_\Gamma(\gamma')$$

\item[({\bf B})]
The switch conditions from Definition \ref{weights+1} for the function $\omega_{m,n}$ from (\ref{8-2-1}) may not be satisfied at all vertices of $\Gamma_m$, but at any vertex the error between the two sides of equality (\ref{switch-conditions}) is caused by those of the above chosen lifts $\gamma_m$ which start or end in that vertex. Hence the total error, summed up over all vertices of $\Gamma_m$, 
is bounded by 
$$\frac{2}{n} \sum_{\gamma \in \cal L_n} 
\mu_\Gamma(\gamma)
= \frac{2}{n}\mu_\Sigma(\Lused) =: K(n)
\, ,$$
where the first equality is explained in Remark \ref{total-weight}.

\item[({\bf C})]
We consider in $\Gamma_m$ a single lift $\gamma'_i \in \cal E_m(\gamma')$ of $\gamma'$,
and consider first the case where $\gamma'_i$ doesn't cross over any of the intrinsic vertices of $\Gamma_m$. 
If none of the chosen lifts $\gamma_m$ of the $\gamma \in \cal L_n$ has one of its endpoints in the interior of $\gamma'_i$, then $\omega_{m,n}(\gamma'_i)$ agrees with $\omega_{m,n}(e)$ for any edge $e$ which is traversed by $\gamma'_i$. If some of the $\gamma_m$ have one or both of its endpoints in the interior of $\gamma'$, then $\omega_{m,n}(\gamma'_i)$ and $\omega_{m,n}(e)$ can disagree, but the difference is bounded by $K(n)$. In the other case, where $\gamma'_i$ crosses over a single intrinsic vertex, then for any occurrence of $\gamma'_i$ in any of the $\gamma_m$ the path $\gamma_m$ crosses over the same intrinsic vertex, and indeed at this vertex over the same local edge $\epsilon$ as does $\gamma'_i$. Hence our above definitions give directly:
$$\omega_{m, n}(\gamma'_i) = \omega_{m,n}(\epsilon)
$$
\end{enumerate}

We now pass successively to larger and large $n$, while keeping $m$ fixed, and observe that for $n \to \infty$ we get $K(n) \to 0$. Thus, if we extract a subsequence $(n_k)_{k \in \N}$ of integers such that the functions $\omega_{m, n_k}$ 
converge to a function $\omega_m$, the fact ({\bf B}) noted above shows that $\omega_m$ satisfies the switch conditions (\ref{switch-conditions}). 
To ensure the existence of such a converging subsequence we note that for any of the finitely many edges $e$ of $\Gamma_m$ the value of $\omega_{m, n}(e) \geq 0$ is bounded above by the total measure $\mu_\Sigma(\Lused)$, by the last equality in observation ({\bf A}) above, applied to the special case $|\gamma'| = 1$.

Furthermore, since the Kolmogorov function $\mu_\Gamma$ satisfies (\ref{Kolmo-1}), and since we have chosen above our lifts $\gamma_m$ to respect the restriction $(\bar \gamma)_m = \bar \gamma_m$, we deduce directly that $\omega_m(\bar e) = \omega_m(e)$ for any edge $e$ of $\Gamma_m$. It follows that $\omega_m$ is a weight function on $\Gamma_m$ as given through Definition \ref{weights+1}. 

In addition, for any of the paths $\gamma'_i \in \cal E_m(\gamma')$ we have
$$
\lim_{k \to \infty} \omega_{m, n_k}(\gamma'_i) = \omega_m(\gamma'_i) \, ,
$$
where the $\omega_{m, n_k}(\gamma'_i)$ on the left hand side have been defined in (\ref{8-2-2}) above, while for the right hand side one needs to apply the general definition for weight functions as given through (\ref{5.4+}) and  (\ref{5.5+}). The claimed equality follows hence from the above observation ({\bf C}).

Finally, from ({\bf A}) we deduce that the ``measure'' of the path $\gamma'$, 
defined as on the right hand side of equality (\ref{5.4}) by the weight function $\omega_m$, agrees precisely with the value of $\mu_\Gamma(\gamma')$ given by the measure $\mu_\Sigma$:
\begin{equation}
\label{measure-agrees}
\mu_\Gamma(\gamma') = 
\sum_{\gamma'_i \in \cal E_m(\gamma')}\omega_{m}(\gamma'_i)
\end{equation}
In particular, the sum on the right hand side turns out to be independent from our above choice of the subsequence $(n_k)_{k \in \N}$.

To conclude, it remains now simply to repeat the above explained procedure for any level graph $\Gamma_m$ to define weight functions $\omega_m$ for all integers $m \geq 0$. Through the compatibility conditions (\ref{compatibility-equation}) and (\ref{local-compatibility}) any such $\omega_m$ defines weight functions $\omega_{\ell, m}$ on $\Gamma_\ell$, for any non-negative integer $\ell \leq m$. Thus a standard diagonal argument, 
using as above the total measure $\mu_\Sigma(\Lused)$ as upper bound to the functions $\omega_{\ell, m}$,
allows us to extract from the family of $\omega_{\ell, m}$ a family of weight functions which are compatible, i.e. a weight tower $\bvec \omega$ on the graph tower $\bvec \Gamma$. From (\ref{measure-agrees}) and Proposition \ref{currents-weights} we now obtain that one has indeed
$$\mu_\Gamma^{\tiny \bvec \omega} = \mu_\Gamma$$
and thus, as desired:
$$\mu_\Sigma^{\tiny \bvec \omega} = \mu_\Sigma$$
\end{proof}

\begin{rem}
\label{disconnected}
(1)
As in some contexts this may occur naturally and be relevant, we'd like to point out that in the definition of graph towers one could as well have omitted the condition that the level graphs are connected: all the definitions and results from this section and the previous ones remain valid without this hypothesis.

\smallskip
\noindent
(2)
If $\bvec \Gamma$ is a graph tower with ``pathologies'' as in Remark \ref{pathological}, then its maximal non-pathological subtower $\bvec \Gamma'$ has the property that for any weight tower $\bvec \omega$ on $\bvec \Gamma$ there is a weight tower $\bvec \omega'$ on $\bvec \Gamma'$ which defines the same invariant measure:
\begin{equation}
\label{8.5}
\mu_\Sigma^{\tiny \bvec \omega'} = \mu_\Sigma^{\tiny \bvec \omega}
\end{equation}
Here we interpret both, $\mu_\Sigma^{\tiny \bvec \omega'}$ and $\mu_\Sigma^{\tiny \bvec \omega}$ as invariant measures on all of $\Sigma(\Gamma_0) = \Sigma(\Gamma'_0)$, due to the potential small differences 
in the used laminations pointed out in Remark \ref{not-lam} (2).

The reason for equality (\ref{8.5}) is that both ``pathologies'' pointed out in Remark \ref{pathological} concern edges which cannot carry positive weight for any weight function as defined in section \ref{weights-currents}: In the case of non-surjective level transition maps this follows directly from the compatibility conditions (\ref{switch-conditions}) and (\ref{local-compatibility}). For level transition maps with folding vertices this is an immediate consequence of the switch condition (\ref{switch-conditions}).

\end{rem}

 
\section{Weight vectors}
\label{sec:weight-vectors}

The purpose of this section is to 
introduce a new tool which will allows us for most applications to
greatly simplify the technology introduced in the sections \ref{graph-dialects} -- \ref{measure-to-weight}. We will show below that the
weight towers considered so far can be replaced by much simpler ``vector towers'' (see Definition \ref{vector-towers}). 
The advantage of vector towers with respect to weigh towers is that local edges can be ignored, while one is still able to describe precisely the same set of invariant measures. On the other hand, for certain computational questions (see for instance Corollary 3.5 of \cite{BHL2}), and also for more sophisticated investigations (for example concerning the associated complexity function), the local edges and their weights turn out to be still very useful.

\smallskip

In order 
to present the transition from weight towers to this easier technology with preciseness, we first work with the convention (as in Convention \ref{ft-short}) that the graph tower 
$\bvec \Gamma = ((\Gamma_n)_{n \in \N \cup \{0\}} , (f_{m,n})_{0 \leq m\leq n})$ is given in short-edge dialect, and that
$((\Gamma^*_n)_{n \in \N \cup \{0\}} , (f^*_{m,n})_{0 \leq m\leq n})$
is the associated long-edge graph tower (see Remark \ref{gt-translation}). Below we also need to specify an (arbitrary) choice of ``positively'' oriented edges ${\rm Edges}^+(\Gamma)$ among the set of all edges ${\rm Edges}^\pm(\Gamma)$, for any given graph $\Gamma$, as explained in subsection \ref{S2-graphs}.

\smallskip

For any of the graph maps $f^*_{m,n}$ given by a graph tower $\bvec \Gamma$ there is defined (see subsection \ref{S2-graph-maps})
a non-negative {\em transition matrix} 
$$M(f^*_{m,n}) = (m_{e^*, e'^*})_{e^*\in {\rm Edges}^+(\Gamma^*_m), \, e'^* \in {\rm Edges}^+(\Gamma^*_n)} \, ,$$
with coefficients 
$m_{e^*, e'^*} = |f^*_{m,n}(e'^*)|_{e^*} + |f^*_{m,n}(e'^*)|_{\bar e^*}$ 
equal to the number of times that $f^*_{m,n}(e'^*)$ crosses over $e^*$ or over $\bar e^*$ (in both cases counted positively).
From the compatibility condition for graph towers (Definition \ref{graph-towers-1} (c)) and equality (\ref{matrix-product}) one obtains directly that
$$M(f^*_{k,n}) = M(f^*_{k,m}) M(f^*_{m,n})$$
holds for all integers $n \geq m \geq k \geq 0$.

\smallskip

For any weight function $\omega^*_n$ on a long-edge 
level graph $\Gamma^*_n$, induced as described in Remark \ref{equal-weights2} by a weight function $\omega_n$ on the short-edge level graph $\Gamma_n$,we consider the associated {\em weight vector} 
\begin{equation}
\label{associated-weight-vector}
\vec v(\omega_n) := (\omega^*_n(e_i^*))_{e_i^* \in {\rm Edges}^+(\Gamma^*_n)}\, ,
\end{equation} 
thought of as column vector.  We deduce from the compatibility conditions (\ref{compatibility-equation})
that for any weight tower $\bvec \omega = (\omega_n)_{n \in \N\cup\{0\}}$ on $\bvec \Gamma$, and for any integers $n \geq m \geq 0$, the associated weight vectors satisfy the following equations:
$$\vec v(\omega_m) = M(f^*_{m,n}) \vec v(\omega_n)$$

This gives rise to the following:

\begin{defn}
\label{vector-towers}
A {\em vector tower} $\bvec v$ on a given graph tower $\bvec \Gamma$ is a family $\bvec v = (\vec v_n)_{n \in \N \cup\{0\}}$ of functions $\vec v_n : {\rm Edges}^+(\Gamma^*_n) \to \R_{\geq 0}$ on the set of positively oriented long edges of the level graphs $\Gamma_n^*$ of $\bvec \Gamma$.  The functions $\vec v_n$ are thought of as column vectors $\vec v_n = (\vec v_n(e_i^*))_{e_i^* \in {\rm Edges}^+(\Gamma^*_n)}$, and they must satisfy the 
{\em compatibility equalities}
\begin{equation}
\label{vector-compatibility}
\vec v_{_m} = M(f^*_{m,n}) \vec v_{n}
\end{equation}
for all $n \geq m \geq 0$.
\end{defn}

\begin{rem}
\label{several}
If $\bvec \Gamma$ is an expanding graph tower, and if $\bvec v= (\vec v_n)_{n \in \N \cup\{0\}}$ is a vector tower on $\bvec \Gamma$,
then we have:
\begin{enumerate}
\item
For any constant $C > 0$, any fixed $m \geq 0$ and any sufficiently large 
$n \geq 0$, the matrix $M(f^*_{m,n})$ has in every column (i.e. for every edge $e'^*$ of $\Gamma^*_n$) a coefficient $m_{e^*, e'^*} > C$.
This is a direct consequence of the hypothesis that $\bvec \Gamma$ be expanding (see Definition \ref{expanding}).
\item
As an immediate consequence of (1) we derive from equality (\ref{vector-compatibility}):
$$\lim_{n \to \infty} \max\{\vec v_n(e^*) \mid e^* \in \Edges^+(\Gamma^*_n)\} = 0$$
\item
For fixed $m \geq 0$ and $n \to \infty$, the size of the vectors $\vec v_n = (\vec v(e_i^*))_{e_i^* \in {\rm Edges}^+(\Gamma^*_n)}$ from $\bvec v$ may in general grow arbitrarily large (and hence the matrices $M(f^*_{m,n})$ would become less and less ``square''). Despite of this phenomenon, one obtains the following improvement of the limit result in (2) above:
$$\lim_{n \to \infty} \sum_{e^* \in \text{Edges}^+(\Gamma^*_n)} \vec v_n(e^*)  = 0$$
To derive this from (1) above, it suffices to consider the sum of the coefficients of the vector $\vec v_m$, and to observe from $\vec v_m = M(f^*_{m,n}) \vec v_n$ in (\ref{vector-compatibility}) that this sum is alternatively obtained as ``scalar product'' $\vec \sigma_{m,n}^* \cdot \vec v_n$, where the row vector $\vec \sigma_{m,n}^*$ is obtained from $M(f^*_{m,n})$ by summing up each column.
\end{enumerate}
\end{rem}

\begin{prop}
\label{almost-weights}
(1)
Let $\bvec \Gamma = ((\Gamma_n)_{n \in \N \cup \{0\}} , (f_{m,n})_{0 \leq m\leq n})$  be an expanding graph tower.
Then any vector tower $\bvec{v} = (\vec v_n) _{n \in \N \cup \{0\}}$ on $\bvec \Gamma$
determines a weight tower 
$\bvec \omega(\bvec{v}) = (\omega_n)_{n \in \N \cup\{0\}}$ 
on $\bvec \Gamma$, 
with the property 
that for any $n \geq 0$ one has (using 
Remark \ref{equal-weights3}) 
\begin{equation}
\label{weight-on-edge}
\omega^*_n(e^*) = \omega^*_n(\bar e^*) := \vec v_n(e^*)
\end{equation}
for all long edges $e^* \in {\rm Edges}^+(\Gamma^*_n)$.
In particular, 
the vector tower $\bvec v$ determines an invariant measure 
$$\mu_\Sigma^{\tiny \bvec{v}} := \mu_\Sigma^{\tiny \bvec \omega(\bvec{v})}$$
on $\Lused$.

\smallskip
\noindent
(2)
In the special case where the vector tower $\bvec v$ is given via (\ref{associated-weight-vector}) through a weight tower 
$\bvec \omega = (\omega_n)_{n \in \N \cup\{0\}}$ on $\bvec \Gamma$, one obtains:
$$\mu_\Sigma^{\tiny \bvec{v}} = \mu_\Sigma^{\tiny \bvec \omega} \qquad \text{and} \qquad \mu_\Gamma^{\tiny \bvec{v}} = \mu_\Gamma^{\tiny \bvec \omega}$$
\end{prop}

\begin{proof}
(1)
Every vector $\vec v_n$ determines via $\omega^*_n(e^*) = \omega^*_n(\bar e^*) := \vec v_n(e^*)$ a non-negative function $\omega^*_n$ on the edges $e^*$ of $\Gamma^*_n$.
This in turn defines a non-negative function $\omega_n$ on the corresponding short edges $e$, obtained by the canonical subdivision of $\Gamma^*_n$, by declaring $\omega_n(e) = \omega^*_n(e^*)$ if $e$ is contained in $e^*$.
We extend this function to the local edges $\epsilon_i$ of $\Gamma^*_n$ by defining
$$\omega_{n}(\epsilon_i) := \sup_{t \geq n} \omega_{n, t} (\epsilon_i)$$
with
$$
\omega_{n, t}(\epsilon_i) := 
\sum_{e'^*_k \in {\rm Edges}^+(\Gamma^*_t)} m^{(t)}_{i, k}\, \, \omega_t^*(e'^*_k) \, ,$$
where $m^{(t)}_{i, k} = |f^*_{n,t}(e'^*_k) |_{\epsilon_i} +|f^*_{n,t}(e'^*_k) |_{\bar \epsilon_i}$ 
equals the number of times that $f^*_{n,t}(e'^*_k)$ crosses over $\epsilon_i$ or $\bar \epsilon_i$. We observe that the equality (\ref{vector-compatibility}) implies $\omega_{n, t}(\epsilon_i) \leq \omega_{n, t'}(\epsilon_i)$ for $t \leq t'$, so that one has: 
\begin{equation}
\label{sup-lim}
\omega_{n}(\epsilon_i) := \lim_{t \to \infty} \omega_{n, t} (\epsilon_i)
\end{equation}

We now observe that from the compatibility equalities on the $\vec v_n$ one obtains the compatibility conditions for the functions $\omega_n$. Indeed, the 
equality (\ref{compatibility-equation}) from Definition \ref{weights+2} follows directly from our assumption $\vec v_m = M(f^*_{m,n}) \vec v_n$ in (\ref{vector-compatibility}), while equality (\ref{local-compatibility}) is a direct consequence of the above definition of the $\omega_n(\epsilon_i)$ and equality (\ref{sup-lim}).

The function $\omega_n$ also satisfies the switch conditions (\ref{switch-conditions}),
by the following argument:

From the definition of $\omega_n$ we observe that
any possible error in the switch condition equalities (\ref{switch-conditions}) 
must come from the fact that, for any level graph $\Gamma^*_t$ with $t \geq n$ and any of its edges $e'^*_k$, 
the initial (or terminal) edge of the image path $f^*_{n, t}(e'^*_k)$ does not contribute enough to $\omega_n(\epsilon_i)$, for any of the local edges $\epsilon_i$ adjacent to the initial (or terminal) vertex of $f^*_{n, t}(e'^*_k)$.
It follows that the total error in the switch conditions for the functions $\omega_{n,t}$, summed over all vertices of $\Gamma_n$, is bounded by 
\begin{equation}
\label{error-estimate}
2 \sum_{e^* \in \text{Edges}^+(\Gamma^*_t)} \vec v_t(e^*)\, ,
\end{equation}
which tends by Remark \ref{several} (3) to $0$ if $t \to \infty$. Hence we deduce from the definition
$\omega_{n}(\epsilon_i) := \underset{t \to \infty}{\lim} \,\omega_{n, t} (\epsilon_i)$ 
in (\ref{sup-lim}) above that the function $\omega_n$ satisfies the switch conditions. 

We have thus shown that the family of $\omega^*_n$ defines indeed a weight tower $\bvec \omega =: \bvec \omega(\bvec v)$ on $\bvec \Gamma$, and hence, by Remark \ref{equal-weights3} and Proposition \ref{almost-weights-}, an invariant measure $\mu_\Sigma^{\tiny \bvec{v}} = \mu_\Sigma^{\tiny \bvec \omega}$ on $\Lused$.

\smallskip
\noindent
(2)
This is a direct consequence of the above definition of the weight tower $\bvec \omega(\bvec v)$, which gives $\bvec \omega(\bvec v) = \bvec \omega$.
${}^{}$
\end{proof}

We recall from Definition \ref{length-defn} that 
for any finite reduced path $\gamma$ in $\Gamma_0$ and any long edge $e^*$ 
of some level graph $\Gamma^*_n$ one denotes by $|f^*_n(e^*)|_\gamma$ the number of occurrences 
of $\gamma$ as subpath in the path $f^*_n(e^*)$. We also recall (see Proposition \ref{K-function-measure} (2)) that any such $\gamma \in \cal P(\Gamma_0)$ determines a cylinder $C_\gamma$, and that for any invariant measure $\mu_\Sigma$ for $\Gamma_0$ the measure of $C_\gamma$ is given by the Kolmogorov function $\mu_\Gamma: \cal P(\Gamma_0) \to \R_{\geq 0}$ associated to $\mu_\Sigma$ via
$$\mu_\Sigma(C_\gamma) = \mu_\Gamma(\gamma)\, .$$

\begin{prop}
\label{cylinder-vectors}
For any vector tower $\bvec v  = (\vec v_n) _{n \in \N \cup \{0\}}$ on an expanding graph tower $\bvec \Gamma = ((\Gamma_n)_{n \in \N \cup \{0\}} , (f_{m,n})_{0 \leq m\leq n})$, and for any finite reduced path $\gamma$ in $\Gamma_0$, the measure $\mu_\Sigma^{\bvec v}$ defined by $\bvec v$ takes on the cylinder $C_\gamma$ the value
$$
\mu_\Sigma^{\tiny\bvec v}(C_\gamma) = \lim_{n \to \infty} \sum_{e^* \,\in \text{\rm Edges}^\pm(\Gamma^*_n)} \vec v_n(e^*) \cdot |f^*_n(e^*)|_\gamma \, .
$$
\end{prop}

\begin{proof}
We pick any level $n \geq 0$ which satisfies $|\gamma| \leq {\mile} _{\tiny\bvec \Gamma}(n)+1$, and consider as in Definition \ref{weights-of-paths} (2) the set $\cal E_n(\gamma)$ of lifts $\gamma_i$ of $\gamma$ to $\Gamma_n$. We thus obtain a disjoint-union decomposition
$$\cal E_n(\gamma) = \cal E_n^{\rm non-intr}(\gamma) \disjoint \cal E_n^{\rm intr}(\gamma)$$
into the set $\cal E_n^{\rm non-intr}(\gamma)$ of such lifts $\gamma_i$ which do not cross over any intrinsic vertex, and the set $\cal E_n^{\rm intr}(\gamma)$ of such lifts $\gamma_j$ which cross over precisely one intrinsic vertex (not counting the initial or terminal vertex of $\gamma_j$). 

From Definition \ref{weights-of-paths} (3), Proposition \ref{weights-via-lifts} and Proposition \ref{currents-weights} we know that for any weight tower $\bvec \omega = (\omega_n)_{n \in \N\cup\{0\}}$ on $\bvec \Gamma$ one has 
$$\mu_\Sigma^{\tiny \bvec{\omega}}(C_\gamma) = \mu_\Gamma^{\tiny \bvec{\omega}}(\gamma) = \sum_{\gamma_i \in \cal E_n(\gamma)} \omega_n(\gamma_i)$$
and thus
\begin{equation}
\label{9-5-1}
\mu_\Sigma^{\tiny \bvec{\omega}}(C_\gamma) = \sum_{\gamma_i \in \cal E_n^{\rm intr}(\gamma)} \omega_n(\gamma_i) + \sum_{\gamma_i \in \cal E_n^{\rm non-intr}(\gamma)} \omega_n(\gamma_i)
\end{equation}

For the weight tower $\bvec \omega = \bvec \omega(\bvec v)$  given by the vector tower $\bvec v$ one has $\omega^*_n(e^*) = \vec v_n(e^*)$ for any long edge $e^*$. Furthermore, there is a 1-1 correspondence between the occurrences of any of the $\gamma_i \in \cal E_n^{\rm non-intr}(\gamma)$ as subpath of some long edge $e^*$ (understood as edge path in the subdivided graph $\check \Gamma_n$ in short-edge dialect) on one hand, and the occurrences of $\gamma$ as subpath of $f_n(e^*)$ on the other. 
Since for any such $\gamma_i$ one has $\omega_n(\gamma_i) = \omega^*_n(e^*)$  (and furthermore $\omega^*_n(e^*)= \vec v_n(e^*)$),
this gives:
\begin{equation}
\label{9-5-2}
\sum_{\gamma_i \in \cal E_n^{\rm non-intr}(\gamma)} \omega_n(\gamma_i) = \sum_{e^* \in \text{Edges}^\pm(\Gamma^*_n)} \vec v_n(e^*) \cdot |f^*_n(e^*)|_\gamma
\end{equation}
This sum differs from $\mu_\Sigma^{\tiny \bvec{\omega}}(C_\gamma) = \mu_\Gamma^{\tiny \bvec{\omega}}(\gamma)$ exactly by the term
$$
\sum_{\gamma_j \in \cal E_n^{\rm intr}(\gamma)} \omega_n(\gamma_j) = 
\sum_{\gamma_j \in \cal E_n^{\rm intr}(\gamma)} \omega_n(\epsilon_j)
\, ,$$
where $\epsilon_j$ is the local edge traversed by $\gamma_j$ at the unique intrinsic vertex crossed over by $\gamma_j$.

Since any long edge $e^*$ of $\Gamma^*_n$, again understood as edge path in the subdivided level graph $\check \Gamma_n$, 
doesn't cross over any intrinsic vertex, for any of the paths $\gamma_j \in \cal E_n^{\rm intr}(\gamma)$ there is precisely one edge $e^* \in {\rm Edges}^\pm(\Gamma_n)$ such that the initial vertex $v_0$ of $\gamma_j$ must lie in the terminal boundary subpaths of length $|\gamma|$ of $e^*$. 
Since $v_0$ must be different from the terminal boundary vertex $v_1$ of $e^*$, there are precisely $|\gamma|$ possibilities for the position of $v_0$ on $e^*$. We fix one of them and consider the set $\cal E_n^{v_0}(\gamma)$ of all possible $\gamma_j \in \cal E_n^{\rm intr}(\gamma)$ with the same initial vertex $v_0$ on $e^*$. All these $\gamma_j \in \cal E_n^{v_0}(\gamma)$ cross over the same intrinsic vertex $v_1 = \tau(e^*)$,
so that any $\gamma_j \in \cal E_n^{v_0}(\gamma)$ crosses over a local edge $\epsilon_j$ at $v_1$, and $\gamma_j$ is uniquely determined by $\epsilon_j$, among all paths in $\cal E_n^{v_0}(\gamma)$. From the switch condition (\ref{switch-conditions}) we thus obtain:
\begin{equation}
\label{9-5-3}
\sum_{\gamma_j \in \cal E_n^{v_0}(\gamma)} \omega_n(\gamma_j) = 
\sum_{\gamma_j \in \cal E_n^{v_0}(\gamma)} \omega_n(\epsilon_j)
\leq
\omega^*_n(e^*)
\end{equation}

Since for every $e^* \in {\rm Edges}^\pm(\Gamma_n)$ there are precisely $|\gamma|$ possible choices for the vertex $v_0$ on $e^*$, we obtain from (\ref{9-5-3}) and from $\omega^*(\bar e^*) = \omega^*(e^*)$ the upper bound
$$
\sum_{\gamma_j \in \cal E_n^{\rm intr}(\gamma)} \omega_n(\gamma_j) 
\leq
\sum_{e^* \in \text{Edges}^\pm(\Gamma^*_n)} \omega^*_n(e^*) \cdot |\gamma|
= 2 \, |\gamma| \sum_{e^* \in \text{Edges}^+(\Gamma^*_n)} \omega^*_n(e^*) 
$$
We thus conclude from $\omega^*_n(e^*) = \vec v_n(e^*)$ and from Remark \ref{several} (3) that 
$$\lim_{n \to \infty} \sum_{\gamma_j \in \cal E_n^{\rm intr}(\gamma)} \omega_n(\gamma_j) = 0 \, .
$$
Hence (\ref{9-5-1}) and (\ref{9-5-2}) imply the claim
$$\mu_\Sigma^{\tiny\bvec v}(C_\gamma) = \lim_{n \to \infty} \sum_{e^* \,\in \text{Edges}^\pm(\Gamma^*_n)} \vec v_n(e^*) \cdot |f^*_n(e^*)|_\gamma \, .
$$
\end{proof}

\begin{rem}
\label{increasing}
For any integers $n \geq m \geq 0$ and any finite reduced path $\gamma$ in $\Gamma_0$ the level map $f_{m, n}$ maps the set of lifts $\cal E_n(\gamma)$ to the set of lifts $\cal E_m(\gamma)$. Since intrinsic vertices are mapped to intrinsic vertices, it follows (using the same terminology as in the last proof) that $f_{m,n}$ maps $\cal E_n^{\rm intr}(\gamma)$ to $\cal E_m^{\rm intr}(\gamma)$, while the paths from $\cal E_n^{\rm non-intr}(\gamma)$ can be mapped to either $\cal E_m^{\rm non-intr}(\gamma)$ or to $\cal E_m^{\rm intr}(\gamma)$. From the compatibility conditions \ref{local-compatibility} we thus obtain: 
$$
\sum_{\gamma_i \in \cal E_m^{\rm intr}(\gamma)} \omega_m(\gamma_i)
\geq
\sum_{\gamma_j \in \cal E_n^{\rm intr}(\gamma)} \omega_n(\gamma_j)
$$
We can hence deduce from the equalities (\ref{9-5-1}) and (\ref{9-5-2}) that the sequence of sums $\underset{e^* \,\in \text{\rm Edges}^\pm(\Gamma^*_n)}{\sum} \vec v_n(e^*) \cdot |f^*_n(e^*)|_\gamma$ in Proposition \ref{cylinder-vectors} is monotonically increasing.
\end{rem}

\begin{rem}
\label{not-long-graph}
The fact that in Definition \ref{vector-towers} and in Proposition \ref{almost-weights} one works with the long-edge dialect is immaterial:  Both, the definition of vector towers as well as the construction of the associated weight tower works just as well on the edge set $\Edges^+(\Gamma_n)$ of an arbitrary graph tower $\bvec \Gamma = ((\Gamma_n)_{n \in \N \cup \{0\}} , (f_{m,n})_{0 \leq m\leq n})$, as long as $\bvec \Gamma$ is strongly expanding in the meaning of Definition \ref{strongly-expanding} (a).

The only reason for passing over to the long-edge dialect lies in its advantage that the transition matrices $M(f_{m,n})$ become smaller and hence more tractable.  However, if in certain given situations some more general expanding graph tower 
turns out to be handier, there is no need to abandon it for the associated long-edge graph tower $((\Gamma^*_n)_{n \in \N \cup \{0\}} , (f^*_{m,n})_{0 \leq m\leq n})$.
\end{rem}

We have now all tools ready needed to give a formal proof of the main result as assembled in section \ref{tools-results}. The small differences in the notation used there and in the present chapter are immaterial, as explained in Remark \ref{not-long-graph}.

\begin{proof} [Proof of Theorem \ref{thm1}]
Statement (1) of Theorem \ref{thm1} is the content of Proposition \ref{almost-weights} (1). Statement (2) follows directly from Proposition \ref{inverse-curr-w}, and the fact that via (\ref{associated-weight-vector}) any weight tower $\bvec \omega$ defines a vector tower $\bvec v$ with $\mu_\Sigma^{\tiny \bvec{v}} = \mu_\Sigma^{\tiny \bvec \omega}$, by Proposition \ref{almost-weights} (2). Statement (3) is an immediate consequence that the definition of the functions $\bvec v \mapsto \bvec \omega(\bvec v)$ and $\bvec \omega \mapsto \mu_\Sigma^{\tiny \bvec \omega}$ are clearly linear. Finally, statement (4) has been proved in Proposition \ref{cylinder-vectors} above.
\end{proof}

\begin{rem} 
\label{sharp-error-estimate}
We now assume that for some vector tower $\bvec v = (\vec v)_{n \in \N\cup\{0\}}$ on some expanding graph tower $\bvec \Gamma = ((\Gamma_n)_{n \in \N \cup \{0\}} , (f_{m,n})_{0 \leq m\leq n})$ the $n$-th level graph $\Gamma_n$ and the corresponding weight vector $\vec v_n$ are known. We will discuss here how this information can be used to concretely determine, for any finite reduced path $\gamma$ in $\Gamma_0$ with $|\gamma| \leq {\mile} _{\tiny\bvec \Gamma}(n)+1$, the measure $\mu_\Sigma^{\tiny \bvec v}(C_\gamma)$ of the cylinder $C_\gamma$. Of course, such a ``determination'' will in general only be possible up to a certain error, which 
will become smaller if the difference $n - |\gamma|$ increases.
We'll give in (\ref{9-7-2})  below a concrete formula for an estimation of $\mu_\Sigma^{\tiny \bvec v}(C_\gamma)$, together with a precise (small) upper bound in (\ref{9-7-3}) for the possible error term.

In a first step we determine, for $\Gamma_n$ understood in short-edge dialect, the finite set $\cal E_n(\gamma)$ of lifts $\gamma_j$ of $\gamma$ to $\Gamma_n$. As in the proof of Proposition \ref{cylinder-vectors} the set $\cal E_n(\gamma)$ partitions into the set $\cal E_n^{\rm intr}(\gamma)$ of paths $\gamma_j$ that cross over an intrinsic vertex, and the set $\cal E_n^{\rm non-intr}(\gamma)$ of paths $\gamma_k$ that are entirely contained in some long edge $e^* =: e^*(\gamma_k)$ of $\Gamma_n$ (where $e^*(\gamma_k)$ is understood as edge path with respect to the short-edge subdivision of $\Gamma_n$).

The determination of the set $\cal E_n^{\rm non-intr}(\gamma)$ allows us to calculate the ``estimation''
\begin{equation}
\label{9-7-2}
S_n(\gamma) := 
\sum_{\gamma_k \in \cal E_n^{\rm non-intr}(\gamma)}
\vec v_n(e^*(\gamma_k))  \, .
\end{equation}
From the fact, that the occurrences of any of the paths $\gamma_j \in \cal E_n(\gamma)$ as subpath of a long edge $e^* \in \Gamma_n^*$ are in 1-1 correspondence with the occurrences of $\gamma$ as subpath of $f_n(e^*)$, we obtain
$$
\sum_{e^* \in \text{Edges}^\pm(\Gamma^*_n)} \vec v_n(e^*) \cdot |f^*_n(e^*)|_\gamma = 
\sum_{\gamma_k \in \cal E_n^{\rm non-intr}(\gamma)}
\vec v_n(e^*(\gamma_k))  \, .
$$
Hence equality (\ref{9-7-2}) together with the equalities (\ref{9-5-1}) and (\ref{9-5-2}) give the ``error term''
\begin{equation}
\label{9-7-1}
\mu^{\tiny \bvec v}(C_\gamma) - S_n(\gamma) = 
\sum_{\gamma_j \in \cal E^{\rm intr}(\gamma)} \omega_n(\gamma_j) \, ,
\end{equation}
where $\omega_n$ is the weight function on $\Gamma_n$ of the weight tower $\bvec \omega(\bvec v)$ associated to the vector tower $\bvec v$ through Proposition \ref{almost-weights}.

A first upper bound of the error term from (\ref{9-7-1}) is derived at the end of 
the proof of Proposition \ref{cylinder-vectors}, where
$$
\sum_{\gamma_j \in \cal E_n^{\rm intr}(\gamma)} \omega_n(\gamma_j) 
\,\,\, \leq
\,\,\, |\gamma| \cdot 
\sum_{e^* \in \text{Edges}^\pm(\Gamma^*_n)} \omega_n(e^*) 
$$
has been shown. A better bound can be obtained as follows:
For any of the paths $\gamma_j \in \cal E^{\rm intr}(\gamma)$ the unique intrinsic vertex $v_j$ of $\Gamma_n$ can be determined, over which $\gamma_j$ crosses, as well as the local edge $\epsilon_j$ used by $\gamma_j$ at $v_j$. By Definition \ref{weights-of-paths} (1) we know $\omega_n(\gamma_j) = \omega_n(\epsilon_j)$, so that it suffices to estimate the value of $\omega_n(\epsilon_j)$\,: From (\ref{weight-inequality})
we know that $\omega_n(\epsilon_j) \leq \omega_n(\hat e^*_j(1))$ and $\omega_n(\epsilon_j) \leq \omega_n(\hat e^*_j(2))$, where $\hat e^*_j(1)$ and $\hat e^*_j(2)$ are the two (uniquely determined) non-local edges adjacent to $\epsilon_j$ in the blow-up dialect derived from the long edges in $\Gamma_n^*\,$: there are corresponding long edges $e^*_j(1)$ and $e^*_j(2)$ in $\Gamma_n^*$ which both start at $v_j$ and which satisfy $\omega_n(e^*_j(1)) = \omega_n(\hat e^*_j(1))$ and $\omega_n(e^*_j(2)) = \omega_n(\hat e^*_j(2))$.
The long edges $e^*_j(1)$ and $e^*_j(2)$ are described alternatively by stating that $\gamma_j$ is a subpath of $\bar e^*_j(1) e^*_j(2)$ or of $\bar e^*_j(2) e^*_j(1)$.
An upper bound for the error term in question is hence given by: 
\begin{equation}
\label{9-7-3}
\big{|}\
\mu_\Sigma^{\tiny \bvec{v}}(C_\gamma) 
\, - \, 
S_n(\gamma)
\big{|}
\leq
\sum_{\gamma_j \in \cal E_n^{\rm intr}(\gamma)}
\min(\vec v_n(e^*_j(1)), \vec v_n(e^*_j(2)))
\end{equation}
\end{rem}
 

\section{Applications and discussions}
\label{applications}

In this section we will outline some first applications of the technology of graph towers provided with vector towers, as presented in the previous sections. Further applications and developments of this technology are given in subsequent work of the authors, see \cite{BHL2, BHL3}. 

We would also like to point the readers attention to two subsections in the appendix to this paper, where terminology and known results are recalled that are used in this section: in subsection \ref{cones} convex linear cones and their extremities are discussed, and in subsection \ref{non-negative-matrices} reducible non-negative matrices and their eigenvectors are reviewed.

\medskip

\subsection{Thin and thick graph towers}
\label{thin-thick-towers}

${}^{}$

For any graph tower $\bvec \Gamma = ((\Gamma_n)_{n \in \N \cup \{0\}} , (f_{m,n})_{0 \leq m\leq n})$ we denote by $\cal V(\bvec \Gamma)$ the set of vector towers $\bvec{v} = (\vec v_n) _{n \in \N \cup \{0\}}$ on $\bvec \Gamma$ as in Definition \ref{vector-towers}. It is a subset of the vector space $\underset{n \in \N \cup \{0\}}{\Pi}\R^{\text{Edges}^+(\Gamma^*_n)}$ and as such it is a convex linear {cone}
(see the Appendix \ref{cones}).
From the definition of $\cal V(\bvec \Gamma) =:\cal V$ it follows directly that $\cal V$ is a closed subset of the non-negative cone 
$\underset{n \in \N \cup \{0\}}{\Pi}\R_{\geq 0}^{\text{Edges}^+(\Gamma^*_n)}$
and that 
hence 
$\cal V$ is salient 
(see equality (\ref{salient})). 

From Proposition \ref{almost-weights} we obtain a surjective map
$$\frak m = \frak m_{\tiny \bvec \Gamma} : \cal V(\bvec \Gamma) \to \cal M(\Lused)$$
into the set $\cal M := \cal M(\Lused)$ of invariant measures on the 
used lamination $\Lused \subset \Sigma(\Gamma_0)$. 
 By statement (3) of Theorem \ref{thm1} this map $\frak m$ is linear.
It follows in particular (see Lemma \ref{extremal-preimage}) that for finite dimensional $\cal V$
the set $\cal V_{ext} \subset \cal V$ of extremal points of $\cal V$ maps to a set $\frak m(\cal V_{ext})$ that contains the extremal points of $\cal M$, which is also a salient closed convex linear 
cone (see Remark \ref{cone-of-measures}). The extremal points of $\cal M$ are precisely the ergodic measures on 
$\Lused$, 
and they are linearly independent (see Remark \ref{cone-of-measures}), 
so that the cone $\cal M$, if finite dimensional, is simplicial (see Definition \ref{simplicial}).

On the other hand, even in the special case where $\cal V$ is finite dimensional, it doesn't seem clear whether the extremal directions of $\cal V$ are projectively independent or not, given that the map $\frak m$ is possibly non-injective.

\medskip

For any integer $n \geq 0$ the map $\bvec{v} \mapsto \vec v_n$ defines a {\em level quotient map}
$$\frak m_n: \cal V(\bvec \Gamma) \to \R^{\text{Edges}^+(\Gamma^*_n)}\, ,$$
and the image $\cal C_n := \frak m_n(\cal V(\bvec \Gamma))$ is again a salient closed convex linear 
cone contained in $\R_{\geq 0}^{\text{Edges}^+(\Gamma^*_n)}$. 
From equality (\ref{vector-compatibility}) we see that the transition matrices $M(f^*_{m, n})$ of the level 
maps $f^*_{m, n}$ satisfy
$$\frak m_m = M(f^*_{m, n}) \circ \frak m_n\, ,$$
so that we have:
$$\cal C_m = M(f^*_{m, n})(\cal C_n)$$
As a consequence, we deduce:
\begin{equation}
\label{dimension-decrease}
\dim \cal C_m \leq \dim \cal C_n \qquad \text{for all} \qquad m \leq n
\end{equation}

\begin{defn-rem}
\label{weight-cones}
The image cones $\cal C_n$ should not be confused with the perhaps more natural {\em weight cones} $\cal C^{\omega}_n$, by which 
we mean  the set of all weight functions on the level graph $\Gamma_n$. Of course we have
$$\cal C_n \subset \cal C^{\omega}_n$$
for any level $n \geq 0$, but in general the inclusion will be strict. The 
transition 
matrix $M(f^*_{m, n})$ of any level  
map $f^*_{m, n}$ defines a map 
$$\cal C^{\omega}_n \to \cal C^{\omega}_m$$
which extends the above map on $C_n$, but in general this map $\cal C^{\omega}_n \to \cal C^{\omega}_m$ will not be surjective.
\end{defn-rem}

The cone $\cal C_0$ has a special importance, due to the following:

\begin{prop}
\label{level-0-split}
Let $\bvec \Gamma, \cal V, \cal M, \cal C_0, \frak m$ and $\frak m_0$ be as defined above.

\smallskip
\noindent
(1)
The map $\frak m_0$ splits over the map $\frak m$. More precisely, the map
$$\zeta: \cal M \to \cal C_0, \, \mu \mapsto (\mu(C_{e^*}))_{e^* \in \text{\rm Edges}^+(\Gamma^*_0)}$$
(where $C_{e^*}$ denotes as in Proposition \ref{cylinder-vectors}
 the cylinder defined by ${e^*}$, understood as edge path in short-edge dialect)
satisfies
$$\frak m_0(\bvec v) = \zeta(\frak m(\bvec v))$$
for any $\bvec v \in \cal V$.

\smallskip
\noindent
(2)
The cone $\cal C_0$ is equal to the nested intersection 
$$\cal C_0 = \bigcap_{n \,\in\, \N}
\left \{M(f^*_n)\cdot\vec v \,\,\middle| \,\, \vec v \in \R_{\geq 0}^{\tiny \text{\rm Edges}^+(\Gamma^*_n)} \right\} \, .$$ 
\end{prop}

\begin{proof}
(1)
We first consider any weight tower 
$\bvec \omega = (\omega_n)_{n \in \N \cup \{0\}}$ on $\bvec \Gamma$, and note that for $n = 0$ one has ${\mile} _{\tiny\bvec \Gamma}(0) = 1$. Since any edge $e^*$ of $\Gamma^*_0$ satisfies $|e^*| = 1 \leq {\mile} _{\tiny\bvec \Gamma}(0)+1$, from the definition of the Kolmogorov function $\mu_\Gamma^{\tiny \bvec \omega}$ defined by $\bvec \omega$ in equality (\ref{5.4}) we obtain 
$\mu^{\tiny \bvec \omega}_\Gamma(e^*) 
= \omega_0(e^*)$. 

In the special case where $\bvec \omega = \bvec \omega(\bvec v)$ is the weight tower defined 
by the given vector tower $\bvec v = (\vec v_n)_{n \in \N \cup \{0\}}$ on $\bvec \Gamma$
(see Proposition \ref{almost-weights}), then by equality (\ref{weight-on-edge})
one has $\omega_0(e^*) = \vec v_0(e^*)$, so that one obtains $\mu^{\tiny \bvec v}_\Sigma(C_{e^*}) = \mu^{\tiny \bvec v}_\Gamma(e^*) = 
\mu^{\tiny \bvec \omega}_\Gamma(e^*) = \omega_0(e^*)
= \vec v_0(e^*)$, which gives directly the claimed statement (1).

\smallskip
\noindent
(2)
The fact that the left hand side of the claimed equality is included in the right hand side is a direct consequence of the given set-up.
For the converse inclusion we observe from equality (\ref{vector-compatibility}) that every element in the nested intersection of the image cones $M(f^*_n)(\R_{\geq 0}^{\tiny \Edges^+(\Gamma^*_n)})$ is indeed the $\frak m_0$-image of some vector tower.
${}^{}$
\end{proof}


The inequality (\ref{dimension-decrease}) enables us to establish the following classification of graph towers. 
We first recall  (see Remark \ref{erase-finite-part}) that 
from any infinite family of level graphs $\Gamma_{n_k}$ of a graph tower $\bvec \Gamma$ 
(with $n_0 = 0$ and with $n_k < n_{k'}$ for $k < k'$) 
one obtains a {\em telescoped graph tower} $\bvec \Gamma'$, with level graphs $\Gamma'_k := \Gamma_{n_k}$ and level transition maps $f'_{k, k'} := f_{n_k, n_{k'}}$. It has the same used lamination as does $\bvec \Gamma$.

\begin{defn-rem}
\label{thick-thin}
(1)
Any graph tower $\bvec \Gamma$, after having been properly telescoped, belongs precisely to one of the following three classes:
\begin{enumerate}
\item[(a)]
$\bvec \Gamma$ is {\em thick} if one has $\dim \cal C_n > \dim \cal C_m$ for any levels $n > m$.
\item[(b)]
$\bvec \Gamma$ is {\em thin} if $\dim \cal C_n = \dim \cal C_0$ for any level $n \geq 0$.
\item[(c)]
$\bvec \Gamma$ is {\em thinning} if $\dim \cal C_n = \dim \cal C_1$ for any level $n \geq 1$ and $\dim \cal C_0 < \dim \cal C_1$.
\end{enumerate}

\smallskip
\noindent
(2) For any graph tower $\bvec \Gamma$ we define the {\em tower dimension} $\dim \bvec \Gamma$ to be the smallest integer that occurs for infinitely many levels $n \geq 0$ as the number of long edges in the level graph $\Gamma^*_n$, if such an inferior limit exists. Otherwise we set $\dim \bvec \Gamma = \infty$. 
We note that by definition of $\cal V(\bvec \Gamma)$ one has: 
$$\dim \cal V(\bvec \Gamma) \leq \dim \bvec \Gamma$$
We also observe directly that $\dim \bvec \Gamma < \infty$ implies that $\bvec \Gamma$ is (after properly telescoping) thin or thinning.
\end{defn-rem}

It is natural to ask whether quantitative or qualitative 
invariants of a given graph tower 
have a structural relevance, in the sense that they only depend on the used lamination $\Lused$ of the graph tower $\bvec \Gamma$ in question, and not on the particular choice of $\bvec \Gamma$.
In particular this can be asked for the above classification, or for the tower dimension. Indeed,for both the answer is positive, 
if one imposes on $\bvec \Gamma$ 
the additional condition to be ``non-repeating'';
this will be investigated in detail in \cite{BHL3}.

\begin{cor}
\label{thin-case}
If the expanding graph tower $\bvec \Gamma$ is thin, then the map $\frak m$ is 
bijective, and every level cone $\cal C_n$ is isomorphic to the cone $\cal M$. In particular, $\cal V$ and any of the $\cal C_n$ are simplicial.
\end{cor}

\begin{proof}
From the hypothesis that $\bvec \Gamma$ is thin it follows that all the maps $M(f_n^*)$ and thus also all the $\frak m_n$ are isomorphisms. 
$$
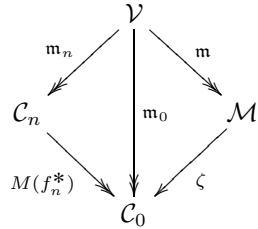
\xymatrix{
       & \mathcal V  \ar@{->>}[ld]_{\mathfrak{m}_n} \ar@{->>}[dd]^{\mathfrak{m}_0} \ar@{->>}[rd]^{\mathfrak m} &  \\
      \mathcal C_n \ar@{->>}[rd]_{M(f_n^*)}  &  & \mathcal M \ar@{->}[ld]^{\zeta} \\
       & \mathcal C_0 & 
    }
$$
Hence part (1) of Proposition \ref{level-0-split} implies that $\frak m$ is injective. From Theorem \ref{thm1} we know that $\frak m$ is surjective.
${}^{}$
\end{proof}

Thin graph towers are implicitly present in several related contexts, such as for ending laminations in the boundary of Teichm\"uller space (see 
\cite{Gabai, Gabai2}), for laminations dual to certain $\R$-trees in Outer space (see \cite{CH, NPR}), and also for $S$-adic subshifts (see \cite{FFT, BD}). This will be further discussed in subsections \ref{stationary-pseudo} and \ref{Outer-space}.

\begin{rem}
\label{ergodic-bound}
(1)
As already pointed out in Corollary \ref{finite-ergodic}, for any expanding 
graph tower $\bvec \Gamma$ with finite tower dimension the number of ergodic probability measures on the used lamination $\Lused$ is bounded above by  $\dim \bvec \Gamma$.

\smallskip
\noindent
(2) Examples of expanding thin $\bvec \Gamma$ where this bound is actually achieved are easy to produce. For example, one can take 
for any level graph $\Gamma_n$ a fixed 1-vertex graph, and define the level maps by mapping each edge $e$ to an edge path $e^t$ with $t \geq 2$. More interesting such examples, with minimal $\Lused$, have recently been constructed by the authors (see \cite{BHL2}).
 \end{rem}

\begin{rem}
\label{invertible-matrix-bound}
${}^{}$
The special situation, where all level transition matrices $M(f^*_{m,n})$ are invertible over $\Z$, occurs frequently in different contexts. In \cite{CH} 
a quick argument is given in a special case, but also valid for general expanding thin $\bvec \Gamma$ with all $M(f^*_{m,n}) \in GL_d(\Z)$ for some $d \geq 1$, which improves the upper bound for the number of ergodic measures (up to rescaling) from $\dim \bvec \Gamma$ to $\dim \bvec \Gamma -1$.
\end{rem}

\medskip

\subsection{Stationary and pseudo-stationary graph towers}
\label{stationary-pseudo}

${}^{}$

We first recall from Definition-Remark \ref{thick-thin} that a graph tower $\bvec \Gamma$ must be thin or thinning, if all of its level graphs $\Gamma_n$ have the same number of edges,  
which then coincides with the tower dimension $\dim \bvec \Gamma$ of $\bvec \Gamma$.

A special case of thin graph towers is given by {\em stationary} graph towers $\bvec \Gamma_{\! \! f}$, which 
are defined by 
a graph $\Gamma$ and a graph self-map $f: \Gamma \to \Gamma$ through the convention 
that $\Gamma_n = \Gamma$ and $f_{m,n} = f^{n-m}$ for all $n \geq m \geq 0$. 

\begin{defn-rem}
\label{train-track-need}
The graph self-map $f: \Gamma \to \Gamma$ of an expanding stationary graph tower $\bvec \Gamma_{\! \! f}$ satisfies:
\begin{enumerate}
\item
$f$ is {\em expanding}\,: For any edge $e$ of $\Gamma$ there exists an integer $t \geq 1$ such that the edge path $f^t(e)$ has length $|f^t(e)| \geq 2$.
\item
$f$ has the {\em train track property}\,: For any integer $t \geq 0$ and any edge $e$ of $\Gamma$ the edge path $f^n(e)$ is reduced.
\end{enumerate}
Conversely, if a graph self-map $f: \Gamma \to \Gamma$ is expanding and has the train track property, then it defines an expanding stationary graph tower $\bvec \Gamma_{\! \! f}$.
\end{defn-rem}

Since any expanding stationary graph tower $\bvec \Gamma_{\! \! f}$ is thin, we obtain from Corollary \ref{thin-case} that the cone $\cal M : = \cal M(\Lused)$ of invariant measures on the used lamination of $\bvec \Gamma_{\! \! f}$ is isomorphic to the cone $\cal V(\bvec \Gamma_{\! \! f})$ of vector towers on $\bvec \Gamma_{\! \! f}$, and also isomorphic to any of the level image cones $\cal C_n := \frak m_n(\cal V(\bvec \Gamma_{\! \! f})) \subset \R_{\geq 0}^{\tiny \Edges^+(\Gamma^+_n)}$. Through the above identification $\Gamma_n = \Gamma$ for any $n \geq 0$ we obtain directly a canonical identification of $\R_{\geq 0}^{\tiny \Edges^+(\Gamma_n)}$ with $\R_{\geq 0}^{\tiny \Edges^+(\Gamma)}$. Thus Proposition \ref{level-0-split} (2) implies natural identifications 
\begin{equation}
\label{core-identification}
\cal C_\infty(M(f)) = \cal C_0 = \cal V(\bvec \Gamma_{\! \! f}) = \cal M \, ,
\end{equation}
where $M(f)$ denotes as before the transition matrix of the map $f$, and $\cal C_\infty$ the nested intersection from Proposition \ref{limit-cone}.

We thus obtain the following theorem, 
which comes close to the main result of \cite{BKMS}. The precise relation between the work of Bezuglyi-Kwiatkowski-Medynets-Solomyak and the work presented here will be described below in Remark \ref{Arnaud's-chouchou}.

Before stating the theorem we recall that any non-negative integer square matrix $M$ possesses up to rescaling only finitely many non-negative eigenvectors which are extremal, and after passing to a suitable power $M^t$ the set of such becomes stable. In Appendix \ref{non-negative-matrices} 
a precise description of these {\em principal eigenvectors} of $M$ is given (see Proposition \ref{principal-ev} and Definition-Remark \ref{barycentric-principal}).

\begin{thm}
\label{mieux-que-Moulinette}
Let $\Lsig_f$ be the used symbolic lamination of a stationary expanding graph tower $\bvec \Gamma_{\! \! f}$, given through some graph map $f: \Gamma \to \Gamma$. Then the set of ergodic measures on $\Lsig_f$ is in 1-1 correspondence with the set of positive scalar multiples of the principal eigenvectors of the transition matrix $M(f)$.  This 1-1 correspondence 
is given by the above identifications (\ref{core-identification}). 

In particular, if $M(f) \in GL_d(\Z)$ (where $d$ is the number of edges of $\Gamma$), then up to scalar multiples the number of distinct ergodic measures on $\Lsig_f$ is bounded above by $\frac{d}{2}$.
\end{thm}

\begin{proof}
The identifications from (\ref{core-identification}) allow a direct translation into the terms defined in subsection \ref{non-negative-matrices}. Hence the claimed statement is a direct consequence of Proposition \ref{limit-cone} and Corollary \ref{half-bound}.
Here non-negative eigenvectors with eigenvalue 0 are excluded by our assumption that $\bvec \Gamma$ be expanding. The hypothesis in Corollary \ref{half-bound} is satisfied, since by telescoping we can replace $f$ by any positive power and then apply Remark \ref{several} (1).
\end{proof}

The above considered situation of stationary graph towers inspires directly the following generalization:

\begin{defn}
\label{pseudo-stationary}
A graph tower $\bvec \Gamma$ is called {\em pseudo-stationary} if for all level transition maps $f_{n,n+1}: \Gamma_{n+1} \to \Gamma_n$ the transition matrices $M(f_{n, n+1})$ are identical (up to permutations), but not necessarily the level maps $f_{n, n+1}$ themselves.
\end{defn}

Since in the proof of Theorem \ref{mieux-que-Moulinette} only the transition matrices but not the transition maps themselves are used, we obtain immediately:

\begin{cor}
\label{Moulinette-eat-your-heart-out}
Theorem \ref{mieux-que-Moulinette} extends verbatim from stationary to pseudo-stationary expanding expanding graph towers.
\qed
\end{cor}

\begin{rem}
\label{Arnaud's-chouchou}
In the most important special case where $\Gamma$ is a 1-vertex graph and the map $f$ respects a preferred ``positive'' orientation on the edges, we have translated in section \ref{dictionary} the set-up from Theorem \ref{mieux-que-Moulinette} 
into that of a given substitution $\sigma: \cal A^* \to \cal A^*$. 

\smallskip
\noindent
(1)
For this case, a bijection between the principal eigenvectors of $M_\sigma$ (called ``distinguished eigenvectors'' in \cite{BKMS}, see Remark \ref{clarification}) and the ergodic measures on the substitution subshift $X_\sigma$ has first been proved in Corollary 5.6 of \cite{BKMS}. It should be noted, though, that due to the particularities of the Bratteli-Vershik machinery (see section \ref{Bratteli-Vershik-presentation}), in the approach presented in \cite{BKMS} one always has to add the additional assumption that the substitution subshift $X_\sigma$ doesn't contain a periodic sequence; this extra assumption is not needed here.

\smallskip
\noindent
(2)
Furthermore, due to the extra effort needed when passing from the substitution to the associated Bratteli-Vershik diagrams (which needs to be made ``proper'' and in this process the incidence matrix in general becomes quite a bit larger), it doesn't seem to be completely obvious whether or not a generalization to pseudo-stationary directive sequences, in analogy to Corollary \ref{Moulinette-eat-your-heart-out} above, is straight forward from the work presented in \cite{BKMS}; a simple quote of the results stated there doesn't quite seem to do it.

\smallskip
\noindent
(3)
In addition, the authors of this paper have to admit a problem in understanding the proof of Corollary 5.6 in \cite{BKMS}: From the given proof we do not quite see that the eigenvectors produced in \cite{BKMS} are actually eigenvectors of the matrix $M_\sigma$ but rather eigenvectors of the much bigger matrix pointed out in (2) above.
\end{rem}

 
\medskip
\subsection{Applications to automorphisms of free groups and current space}
\label{Outer-space}

${}^{}$

The term ``current'' in the context of free groups $\FN$ is the precise analogue of what we have so far called ``invariant measure on some symbolic lamination''. This difference in terminology has in part historic reasons, in part it is due to the absence of a preferred basis in $\FN$, which calls for a basis-free, more algebraic and less combinatorial approach. 
See \cite{Ka1, Ka2} for a general introduction to currents over $\FN$ and to current spaces.

We will try to explain in subsection \ref{cur-and-lam} below the precise relation between the algebraic and the combinatorial approach, without getting lost in technical details. In subsection \ref{out-results} we will exhibit the relevance of Theorem \ref{mieux-que-Moulinette} in the free group context. Finally, in subsection \ref{out-discussion} we will discuss open ends and pose some questions.

\subsubsection{Currents and algebraic laminations}
\label{cur-and-lam}

${}^{}$

Any free group $\FN$ of finite rank $N \geq 2$ has a countable number of bases, and the transition from one to the other is given by an automorphism of $\FN$. The situation is similar to the choice of a marking on a surface $S_g$, and indeed, the group $\Out(\FN)$ of outer automorphisms of $\FN$ is closely related (but even more challenging) than the mapping class group $\text{Mod}_g$ of $S_g$.

On the other hand, 
for any set of letters $\cal A = \{a_1, \ldots, a_N\}$ the free group $F(\cal A)$ over $\cal A$ is naturally isomorphic to $\FN$, 
and the choice of any such isomorphism establishes $\cal A$ as basis of $\FN$.
This defines canonically an embedding of the free monoid $\cal A^*$ into $F(\cal A) \cong \FN$, as well as an identification $\FN \cong \pi_1 (R(\cal A))$, where 
(as in section \ref{dictionary}) the ``rose'' 
$R(\cal A)$ is the 1-vertex graph with $N$ oriented edges labelled by the $a_i$. More generally, the fundamental group of any finite connected graph $\Gamma$ is a free group $\FN$ of finite rank $N \geq 0$. However, in general $\Gamma$ may have more than one vertex, and hence up to $3N - 3$ edges (assuming $N \geq 2$; the cases $N = 1$ and $N = 0$ will not be treated here).

As indicated already above, a {\em current} on $\FN$ can now be viewed as a ``letter-free'' version of a subshift equipped with a shift-invariant measure. Indeed, once a {\em marking isomorphism} $\pi_1 \Gamma \cong \FN$ is specified, 
then a current $\mu$ on $\FN$ gives canonically rise to a shift-invariant measure $\mu_\Sigma$ for $\Gamma$, and conversely.
For the space $\Curr(\FN)$ of currents on $\FN$ this defines a canonical bijection:
\begin{equation}
\label{canonical-currents}
\Curr(\FN) \longleftrightarrow \left\{\mu_\Sigma \,\,\middle| \,\, \mu_\Sigma \text{ shift-invariant measure on }\Sigma(\Gamma) \right\}
\end{equation}
For more details see
\cite{CHL1-III, Ka2}.

A similar ``letter free'' approach to symbolic laminations has been given in 
\cite{CHL1-I}, where the space $\Lambda^2(\FN)$ of {\em algebraic laminations} $L^{F_N}$ 
over a free group $\FN$ has been introduced and studied. 
As before, for any graph $\Gamma$, provided with an identification $\FN \cong \pi_1 \Gamma$, there is a canonical bijection: 
\begin{equation}
\label{canonical-laminations}
\Lambda^2(\FN)  \longleftrightarrow \{\Lsig \subset \Sigma(\Gamma) \mid \Lsig \text{ symbolic lamination on } \Gamma \}
\end{equation}
Any current $\mu \in \Curr(\FN)$ determines a {\em support} $\supp(\mu) \in \Lambda^2(\FN)$, and the naturality of the above described set-up effects that this support map commutes via the two given bijections with the map that associates to any shift-invariant measure $\mu_\Sigma$ on $\Sigma(\Gamma)$ the symbolic lamination $\Lsig(\mu_\Sigma)$ 
from (\ref{support-measure})
on the given graph $\Gamma$.

\begin{rem}
\label{known}
The following facts are well-known, see 
\cite{CHL1-I, CHL1-III}:
\begin{enumerate}
\item
Any automorphism $\phi \in \Out(\FN)$ induces a homeomorphism $\phi_\Lambda$ on the space of algebraic laminations $\Lambda^2(\FN)$.
\item
Any automorphism $\phi \in \Out(\FN)$ induces homeomorphisms 
$\phi_{C}$ on the space $\Curr(\FN)$ and  
$\phi_{\Pr C}$ on it projectivization $\PCurr$. 
The homeomorphism $\phi_{C}$ is furthermore ``linear'', in that for any $\mu = \sum \lambda_i \mu_i$ it satisfies $\phi_{C}(\mu) = \sum \lambda_i \phi_C(\mu_i)$.
\item
The two homeomorphisms 
$\phi_\Lambda$ and $\phi_{\Pr C}$
commute via the support map $\mu \mapsto \supp(\mu)$. 
The latter, however, is not continuous.
\end{enumerate}
\end{rem}

$$
  \xymatrix{
    \Curr(\FN)     \ar[r] \ar[d]^{\phi_{C}} & \PCurr \ar[rr]^{\text{Supp}} \ar[d]^{\phi_{\Pr C}} &  & \Lambda^2(\FN) \ar[d]^{\phi_\Lambda} \\
    \Curr(\FN)     \ar[r]  & \PCurr \ar[rr]^{\text{Supp}}  & & \Lambda^2(\FN)  
    }
$$

\bigskip

The space $\Curr(\FN)$ of such currents has been much studied in recent years (see e.g. \cite{ClayPettet-twisting-currents, CHL1-III, GiHo, Ka1, Ka2, KL09, IlyaTatjana-PattersonSullivan, Caglar}), and even more so its projectivization $\PCurr$, which is compact. 
Although not finite dimensional (as is Outer space $\CVN$), the space $\PCurr$ -- or a suitable open and dense subspace (see \cite{KL5}) -- is generally accepted as a useful analogue of Teichm\"uller space (second to $\CVN$), with $\Out(\FN)$ playing the role of the mapping class group $\text{Mod}_g$.

\subsubsection{Fixed currents under the action of single automorphisms}
\label{out-results}

${}^{}$

Individual automorphisms $\phi \in \Out(\FN)$ can be much more complicated than mapping classes, and they are up today not understood in full generality. 
There is a special interest in the fixed point set of the homeomorphism $\phi_{\Pr C}$ induced by any such $\phi$ on $\PCurr$:
among other, these fixed points, i.e. projectivized $\phi$-invariant currents $[\mu]$, are often a useful tool for the study of the intrinsic structure of $\phi$.

The most powerful technology to analyze individual automorphism $\phi \in \Out(\FN)$ is based on train track maps, which exist in a variety of absolute, relative, improved, etc versions (see \cite{BFH1, BH92, Bridson-Groves, FH, Lu_conj-pr_1, Lu-alpha}). We restrict our attention here to the following:

\begin{defn-rem}
\label{train-track-map}
(1)
A graph self-map $f: \Gamma \to \Gamma$ is called an {\em expanding train track map} if it is expanding and has the train track property (as specified in Definition-Remark \ref{train-track-need}).

\smallskip
\noindent
(2) 
The map $f: \Gamma \to \Gamma$ is said to {\em represent an automorphism} $\phi \in \Out(\FN)$ if for a suitable identification $\pi_1 \Gamma \cong \FN$ the map $f$ induces $\phi$. In this case $f$ is a homotopy equivalence.
\end{defn-rem}

As pointed out in Definition-Remark \ref{train-track-need}, any expanding train track map $f: \Gamma \to \Gamma$ gives rise to an expanding stationary graph tower $\bvec \Gamma_{\! \! f}$, where every level map $f_{n, n+1}$ is identified with the given map $f$. 
More specifically, let us fix a marking isomorphism 
$\theta: \pi_1 \Gamma \overset{\cong}{\longrightarrow} \FN$ and 
assume that the map $f$ represents (with respect to this marking $\theta$) some automorphism $\phi \in \Out(\FN)$. 
Then the set $\cal M(\Lsig_f)$ of shift-invariant measures on the used symbolic lamination 
$\Lsig_f := L^{\tiny \bvec \Gamma_{\! \! f}}$ is via (\ref{canonical-currents}) in
natural bijection with the set $\cal M_f \subset \Curr(\FN)$ of currents with support in the {\em used algebraic lamination 
$L^{F_N}_f$}. This algebraic lamination is canonically defined by the train track map $f$, see
\cite{KL-7steps}, Definition 3.35 and Lemma 3.36, and it corresponds via (\ref{canonical-laminations}) precisely to the used symbolic lamination 
$\Lsig_f$. 
We denote by $\Pr\cal M_f \subset \PCurr$ the set of projectivized currents $[\mu]$ defined by 
any element $\mu$ of $\cal M_f$.

\begin{thm}
\label{thm-BHL-old-new}
Let $f: \Gamma \to \Gamma$ be an expanding train track map which represents an automorphism $\phi \in \Out(\pi_1(\Gamma))$.
After possibly replacing $f$ and $\phi$ by some positive power, we obtain:
\begin{enumerate}
\item
There is a canonical 1-1 correspondence between the set of principal eigenvectors $\vec v_i$ of the transition matrix $M(f)$ on one hand, and the set of projectivized $\phi$-invariant ergodic currents $[\mu_i] \in \Pr \cal M_f$ on the other.
\item
Every convex combination $\sum c_i \vec v_i$ 
of principal eigenvectors $\vec v_i$ of $M(f)$ with same eigenvalue defines a projectivized $\phi$-invariant current $[\sum c_i \mu_i] \in \Pr \cal M_f$.
\end{enumerate}
\end{thm}

(See Proposition \ref{principal-ev} and Definition-Remark \ref{barycentric-principal} for a precise description of the ``principal eigenvectors'' of a non-negative matrix.)

\begin{rem}
\label{BHL-old-comment}
The above Theorem \ref{thm-BHL-old-new} and Proposition \ref{principal-ev} show that there is a canonical injection from the set of non-negative eigenvectors of $M(f)$ into the set of projectively $\phi$-invariant currents in $\cal M_f$.
Through developing our theory of expanding graph towers and vector towers a bit further, it will be shown in 
\cite{BHL3} that this injection is actually a bijection.
\end{rem}

The next proof as well as part of the discussion in the next subsection is written for readers with 
some expert knowledge about Outer space and automorphisms of free groups.

\begin{proof}[Proof of Theorem \ref{thm-BHL-old-new}]
From basic train track theory it is known that if a train track map $f$ represents the automorphism $\phi \in \Out(\FN)$, then one has 
$\phi_\Lambda(L^{F_N}_f) \subset L^{F_N}_f$
 (see for instance \cite{KL-7steps}). 
Furthermore, the algebraic lamination 
$L^{F_N}_f$ contains only finitely many sublaminations, each given by a stratum of the expanding train track map $f$ (see \cite{BFH1}). Recalling that $\phi_\Lambda$ is a homeomorphism and hence bijective, we thus deduce from the last inclusion:
$$\phi_\Lambda(L^{F_N}_f) = L^{F_N}_f$$

Recalling that $\cal M_f \subset \Curr(\FN)$ denotes the set of currents $\mu$ with $\supp(\mu) \subset 
L^{F_N}_f$, we can apply the commutativity between $\phi_\Lambda$ and $\phi_{C}$ from Remark \ref{known} (3) to obtain:
$$\phi_{C} (\cal M_f) = \cal M_f  \qquad \text{and} \qquad \cal M_f  = \phi^{-1}_{C} (\cal M_f)$$
Since any convex combination of currents from $\cal M_f$ gives again a current in $\cal M_f$, the latter is a cone. Furthermore we know from Remark \ref{known} (2) that the map $\cal M_f \overset{\phi_C}{\longrightarrow} \cal M_f$ is linear, and since it is invertible, it is an isomorphism between cones. It follows that $\phi_C$ maps extremal currents of the cone $\cal M_f$ to extremal currents of $\cal M_f$, or in other words: $\phi_C$ permutes the projectivized ergodic currents $[\mu] \in \Pr \cal M_f$.

The space $\cal M_f \subset \Curr(\FN)$ is via (\ref{canonical-currents}) canonically 
identified with the set $\cal M(\Lsig_f)$ of shift-invariant measures $\mu_\Sigma$ on the symbolic lamination $\Lsig_f$. We can hence apply Theorem \ref{mieux-que-Moulinette} to obtain a bijection between the set of ergodic projectivized currents $[\mu] \in \Pr \cal M_f$ on one hand, and the set of principal eigenvectors of the transition matrix $M(f)$ on the other. This gives, for suitable positive powers of $f$ and $\phi$, the bijection stated in part (1) of our our claim.

In order to prove part (2) it suffices to observe that among currents $\mu_i$, which all satisfy $\phi_C(\mu_i) = \lambda \mu_i$ for some fixed $\lambda > 0$, all linear dependencies are preserved under application of $\phi_C$ (see Remark \ref{known} (2)). This shows directly that any eigenvector $\vec v$ of $M(f)$, given as convex combination $\sum c_i \vec v_i$ of non-negative principal eigenvectors $\vec v_i$ of $M(f)$ with same eigenvalue as $\vec v$, defines a projectively $\phi$-invariant current which can be expressed precisely in the same way as convex combination  $\sum c_i \mu_i$ of the $\mu_i$ defined by each $\vec v_i$.
\end{proof}

\begin{rem}
\label{thm-BHL-old-new-3}
For $f, \Gamma$ and $\phi$ as in Theorem \ref{thm-BHL-old-new} it 
follows that the number of projectivized $\phi$-invariant ergodic currents $[\mu_i] \in \Pr \cal M_f$ is bounded above by $\frac{3}{2}(N-1)$, for $N = \text{rank}(\pi_1(\Gamma))$. This follows from statement (1) of Theorem \ref{thm-BHL-old-new}, together with the upper bound 
$\frac{d}{2}$ for the number of principal eigenvectors of the 
expanding non-negative $d \times d$-matrix $M(f)$.

Indeed, the matrix size $d$ is here bounded above by the maximal number of edges in any graph $\Gamma$ with $\pi_1 \Gamma = \FN$,
which gives $d \leq 3N - 3$. 
Furthermore, any principal eigenvector of $M(f)$ must correspond to a distinct primitive diagonal block of $M(f)$, given by an
expanding stratum of the train track map $f$.
But any such stratum must involve at least 2 edges, as $f$ is a homotopy equivalence, thus 
lowering the bound from $3N-3$ to $\frac{3}{2}(N-1)$.
\end{rem}

\subsubsection{Further discussion and questions}
\label{out-discussion}

${}^{}$

Theorem \ref{thm-BHL-old-new} applies to a large class of automorphisms 
$\phi \in \Out(\FN)$, among which, most importantly, are all {\em hyperbolic} (or, equivalently, {\em atoroidal}) automorphisms other than those satisfying a certain technical obstruction (an essential non-closed INP-path in some relative train track representative of $\phi$). 
Hyperbolic $\phi \in \Out(\FN)$ can be characterized by the absence of any non-trivial conjugacy class in $\FN$ that grows 
in length polynomially (or equivalently ``subexponentially''), 
under iteration of $\phi$. The class of hyperbolic automorphisms contains the class of non-geometric iwip (= ``irreducible with irreducible powers'') automorphisms of $\FN$, which are the natural 
strict analogue of pseudo-Anosov mapping classes in $\text{Mod}_g$. 

Such non-geometric iwip $\phi$ are known (see \cite{Caglar}) 
to have North-South dynamics on $\PCurr$, a result that has recently been extended (see \cite{LU, U2}) to a ``generalized North-South'' dynamics on $\PCurr$ for the action of any hyperbolic $\phi \in \Out(\FN)$. In the case where $\phi$ is represented by an expanding train track map $f: \Gamma \to \Gamma$, a ``forward limit simplex'' $\Delta_\phi$ for this dynamics has been exhibited in \cite{LU} in terms of frequencies. We believe that the following natural question has a positive answer:

\begin{question}
\label{limit-simp}
Does the above forward limit simplex $\Delta_\phi$ coincide with $\Pr \cal M_f$ from Theorem \ref{thm-BHL-old-new} ?
\end{question}

\medskip

The situation considered in Theorem \ref{thm-BHL-old-new} fits into a 
somehow more general scheme, where we assume that for a graph tower $\bvec \Gamma$ the level transition maps $f_{m, n}$ are all homotopy equivalences. Once a marking isomorphism $\pi_1 \Gamma_0 \cong \FN$ is specified, such a graph tower $\bvec \Gamma$ defines a sequence of non-metric metric graphs which gives rise to an infinite ``unfolding path'' $\bvec \gamma = \bvec \gamma(\bvec \Gamma)$ in Outer space $\CVN$. Due to the lack of a specified metric on the graphs $\Gamma_n$, this path $\bvec \gamma$ is only determined up to homotopy within the simplicial neighborhood of $\bvec \gamma$, given by the canonical simplicial structure of $\CVN$. Despite the possible local perturbations of $\bvec\gamma$ issuing from this non-determinacy, it follows from Remark \ref{several} (2), together with elementary considerations based on the geometric intersection form between currents and $\R$-trees from \cite{KL09}, 
that for expanding $\bvec \Gamma$ any limit point $[T] \in \partial \CVN$ of the unfolding path $\bvec\gamma(\bvec \Gamma)$ is given by a tree $T$ which contains the active part of the used lamination $\Lused$ in its dual lamination  $L(T)$. Here we mean by the {\em active part} of a symbolic lamination $\Lsig$ the union of the supports of any invariant measure on $\Lsig$.

This set-up has been investigated previously in \cite{CH} and \cite{NPR}.
In \cite{CH} Coulbois and Hilion have used unfolding techniques from their previous work to derive 
$3N-3$ as upper bound\footnote{\,\, The slightly sharper bound stated in \cite{CH} is based on an argument on ``pseudo-invertible matrices'' that unfortunately only holds in special cases.}
for the number of projectivized ergodic currents with support in by $L(T)$,
assuming that the $\FN$-action on $T$ in $\partial \CVN$ is free. Alternatively, the upper bound $3N - 3$ follows from Remark  \ref{ergodic-bound}, since in the above set-up the graph tower $\bvec \Gamma$ is thin of tower dimension $\leq 3N-3$, as this is the maximal number of edges of any graph with fundamental group isomorphic to $\FN$.

However, to our knowledge no example of a tree $T$ in $\partial \CVN$ 
with free $\FN$-action (or with cyclic 
point stabilizers) is known, where the number of projectivized ergodic currents carried by $L(T)$ exceeds the bound $\frac{3}{2}(N-1)$. This bound has been established through Remark \ref{thm-BHL-old-new-3} 
for the special case of certain trees $T$ that are projectively fixed by properly chosen hyperbolic automorphisms $\phi \in \Out(\FN)$: 
There one considers a stationary graph tower $\bvec \Gamma_{\!\! f}$, defined by a train track map $f$ that represents some 
$\phi$.

Note though that in this case the $\R$-tree $T$ in question is not the ``usual'' forward limit tree determined by a left-eigenvector of $M(f)$. In order to describe $T$ by this well developed train track technology (see \cite{GJLL, Lu_conj-pr_1}), one would first have to find a train track map $f_-: \Gamma_- \to \Gamma_-$ which represents $\phi^{-1}$, and then one has to impose further conditions on the size of the stretching factors of the strata of $f_-$ to ensure the existence of some such $T$ with trivial point stabilizers.

\begin{question}
\label{maximal-ergodic}
Let $T$ be an $\R$-tree in $\partial \CVN$ with free $\FN$-action (or with cyclic point stabilizers).
\begin{enumerate}
\item
What is the maximal number of projectivized ergodic currents carried by the dual lamination $L(T)$ ?  
Is it bounded above by $\frac{3}{2}(N-1)$~?
\item
Same question under the additional assumption that $L(T)$ is minimal 
(or minimal up to diagonal leaves).
\end{enumerate}
\end{question}

In the context of the last question (2)
we'd like also to point the reader's attention to the work of D. Gabai on ending laminations for 
orientable surfaces, see \cite{Gabai, Gabai2}. 
The latter has also inspired recent work of Leininger, Lenzhen and Rafi \cite{LLR} and of Brock, Leininger, Modami and Rafi \cite{BLMR} on 
limit sets of particular Teichm\"uller geodesics in the Thurston boundary of Teichm\"uller space, which accumulate on more than one ergodic measure on a given minimal surface lamination.

In Theorem 9.1 of \cite{Gabai} Gabai exhibits 
minimal laminations with $\frac{3}{2}N - 2$ projectively distinct ergodic measures, 
for surfaces $S_{g,1}$ of genus $g$ with one puncture and free fundamental group $\FN$ of rank $N = 2g$.
This amounts to the same bound as given in Remark \ref{thm-BHL-old-new-3} and as proposed above in Question \ref{maximal-ergodic}. Indeed, Gabai's laminations are, in our terminology, used laminations of thin graph towers of tower dimension $d = 3N - 3$ (or smaller), 
with level transition maps that are homotopy equivalences. These graph towers, however, will not be stationary (as are those from 
Remark \ref{thm-BHL-old-new-3}), or else the lamination couldn't be minimal: The corresponding unfolding path in Gabai's case must 
-- just as the Teichm\"uller geodesics in \cite{LLR} and \cite{BLMR} --
vanish into the thin part of Outer space, while any stationary graph tower gives a periodic path which is hence contained in some thick part.

In the special case where the surface lamination is orientable, it can be alternatively described by an interval exchange transformation on $d$ intervals. In this case an upper bound for the number of ergodic measures on a minimal lamination, together with a realization result, is given by the classical work of \cite{Katok, Keane, Veech, Yoccoz}; this bound corresponds to Gabai's results mentioned before. 

This also provides a positive answer to Question \ref{maximal-ergodic} (2) for this special case, since Gabai's laminations are (as are all surface laminations) dual to some $T \in \partial \CVN$. Note here that, due to the puncture in $S_{g,1}$, the $\FN$-action on $T$ will not be free, but has cyclic point stabilizers.
 On the other hand, for most $[T]\in \partial \CVN$ the dual lamination $L(T)$ does not come from a surface lamination.



\bibliographystyle{abbrv}
\bibliography{biblio-BHL1}

\newpage
\section{APPENDIX}

The three subsections of this appendix are independent from each other.

\subsection{Bratteli-Vershik technology}
\label{Bratteli-Vershik-presentation}

${}^{}$

In this appendix we will give a brief account of an alternative technique to investigate subshifts on finite alphabets and their invariant measures, as well as a bit of systematic comparison (see Remark \ref{BV-translation}) between the results obtained via this alternative approach and the analogous results coming from our new tower technology. 
This alternative 
technique is based on Bratteli diagrams provided with a Vershik map. A systematic introduction to these objects and their application can be found for instance in \cite{Dur}. Since none of us authors is an expert on Bratteli-Vershik methods, we limit ourselves here to a rough sketch of the aspects close to the results presented here.

\smallskip

A {\em Bratteli diagram} is defined by an infinite graph $\cal B$ with edge set  $E$ and vertex set $V$, partitioned into finite ``levels'' $E_n$ (for $n \geq 1$) and $V_n$  (for $n \geq 0$). The level $V_0 = \{v_0\}$ is a singleton, and each $E_n$ is a set of oriented edges joining vertices of $V_{n-1}$ to vertices of $V_n$. 

There is a natural {\em level matrix} $M_n$ (denoted in \cite{BKMS2} by $F_{n-1}$), with coefficients indexed by a vertex pair $(v, v') \in V_{n}\times V_{n-1}$ which counts the number of edges in $E_n$ that connect $v'$ to $v$. Denoting for any $v \in V_n$ by $h_v$ the number of paths of length $n$ from $v_0$ to $v$, one obtains immediately for the vectors $\vec h_n := (h_v)_{v \in V_n}$ the matrix equality
$$\vec h_{n+1} = M_{n+1} \vec h_n \, .$$

One denotes by $X_{\cal B}$ the set of infinite paths which start at $v_0$ and move monotonically up through the levels. Every finite initial subpath $\underbar e := e_1 e_2 \ldots e_r$ of some $x \in X_{\cal B}$ defines a {\em cylinder} $[\underbar e]$ which consists of all $x' \in X_{\cal B}$ with $\underbar e$ as initial subpath. There is a natural topology on $X_{\cal B}$ generated by the cylinders as clopen sets.

Two paths in $X_{\cal B}$ are {\em tail-equivalent} if they agree up to some finite initial subpaths.
We denote the tail equivalence class of any $x \in X_\cal B$ by $\langle x \rangle$.

\smallskip

The vertices of $V_n$ partition naturally the set $E_n$ into classes of edges that have the same terminal vertex. Each such {\em incoming edge class} is provided with a total order, and the collection of all these {\em local orders} defines canonically a lexicographical total order on any tail equivalence class of $X_{\cal B}$. The {\em Vershik map} $V_{\cal B}: X_{\cal B} \to X_{\cal B}$ is defined by sending every $x \in X_{\cal B}$ to its successor with respect to this total order.

It follows that every $V_{\cal B}$-orbit is a tail-equivalence class, and conversely. This can be visualized as follows: The total order on the incoming edge classes gives (modulo a ``left-right'' convention) a canonical way to unfold a tail equivalence class $\langle x \rangle$ into an infinite but locally finite one-ended tree $T_{\langle x \rangle}$ embedded in a level-preserving way into the plane (see \cite{FFT}). Any path $x' \in X_\cal B$ belongs to the tail equivalence class of $x$ if and only if $x'$ can be realized by a (unique) path in $T_{\langle x \rangle}$. The Vershik map on the tail equivalence class of $x$ is visualized for any path realized in $T_{\langle x \rangle}$ by simply passing to the right (or ``left'', according to the above taken convention) neighboring path in $T_{\langle x \rangle}$. 

\begin{rem}
\label{Vershik-problems}
There is an obvious problem in the ``careless'' definition of $V_\cal B$ as presented above: there might be maximal (and minimal) elements in $X_\cal B$ with respect to the lexicographic order. There are several ways in the literature how to overcome this problem:

In the approach presented in \cite{Dur} it is axiomatically postulated that there is precisely one such maximal and one such minimal element, and the latter is defined to be the $V_\cal B$-image of the former. It is shown there that under suitable extra hypothesis (going by the name of ``proper'') this approach works well, and in particular it leads to convincing results for simple stationary Bratteli diagrams (which correspond to primitive substitutions).

An alternative approach is pursued in \cite{BKM, BKMS, BKMS2}, where for finite sets of minima and maxima the (countable) set of points in the $V_\cal B$-orbits of the latter are excluded from $X_\cal B$, and measures on all of $X_\cal B$ are considered that are invariant under the restriction of $V_\cal B$ to the ``non-exceptional'' complement in $X_\cal B$ of this countable set. 

Finally, a more general (and perhaps more natural, but less practical) way to deal with this problem is described in \cite{FFT}.
\end{rem}

Bratteli-Vershik theory now proceeds by introducing a number of refinements through additional technical conditions, leading to a vast terminology like {\em simple diagrams}, {\em simple hat}, {\em finite rank}, {\em proper}, etc. We have to leave it to the interested reader to pick this up through the papers which are cited in this 
appendix, or through the references given there.

\smallskip

The main point in our context, however, seems to be the following: From the definition of the Vershik map it follows directly that any two cylinders $[\underbar e]$ and $[\underbar e']$ as defined above are in the same $V_\cal B$-orbit if the paths $\underbar e$ and $\underbar e'$ have the same terminal vertex. Hence any $V_\cal B$-invariant measure $\mu$ on $X_\cal B$ must associate to any cylinder $[\underbar e]$ a value $\mu([\underbar e]) =: p(v))$ which only depends on the terminal vertex $v$ of $\underbar e$. Hence $\mu$ defines for any level $n$ a vector $\vec p_n = (p(v))_{v \in V_n}$, and $\mu$ is determined by the family of the $\vec p_n$, for all $n \in \N$. The disjoint union property of cylinders, that arises 
from distinct prolongations of a path $\underbar e$ into the next level, gives 
directly the matrix equality
$$\vec p_n = \,^t\!M_{n+1} \vec p_{n+1}$$
where $^t\!M_{n+1}$ denotes the transpose of $M_{n+1}$. Conversely, this equality assures that the cylinder measures defined by the coordinates of the $\vec v_n$ add up properly to define a measure on $X_\cal B$.  The precise details and more elaborate arguments for these facts are given in Theorem 2.9 of \cite{BKMS} and its proof.

\medskip

There is clearly a very strong formal ressemblance between the family of vectors $\vec p_n$ described above, and the vector towers $\bvec v$ on a graph tower as described in section \ref{tools-results}. This similarity is further enhanced by the strong overlap of our results in the case of stationary graph towers and the main result of \cite{BKMS} for stationary Bratteli diagrams, see Remark \ref{Arnaud's-chouchou}.

It is hence natural to seek for a direct formal translation of one formalism into the other. This, however, seems to be by no means 
an easy matter; we'll give now an account of the situation as presently known, to the best of our understanding:

\begin{rem}
\label{BV-translation} 
There is a natural translation of any level of a Bratteli diagram, provided with the local order at the incoming edge classes as described above, into a substitution: One defines $\cal A_n := V_n$ and $\cal A_{n+1} := V_{n+1}$, and obtains the substitution $\sigma_{n+1}: \cal A^*_{n+1} \to \cal A^*_n$ by defining, for any $v' \in V_{n+1}$, the image $\sigma_{n+1}(v') = v_1 \ldots v_r$, where each $v_i$ is the initial vertex of some edge terminating at $v'$, and the linear order of the $v_i$ in the word $\sigma_{n+1}(v')$ reflects precisely the total order defined on the set of incoming edges at $v'$.

In this way one can ``read off'' from any ordered Bratteli diagram $\cal B$ a directive sequence $\sigma$ of substitutions $\sigma_n$
(and conversely), with the natural hope that the associated subshift $X_\sigma$ gives back the path space $X_\cal B$, where the shift operator ought to be conjugated to the Vershik map $V_\cal B$. That things are not quite as simple transpires already from considering the 1-letter substitution given by $a \mapsto a a$ and the stationary directive sequence built on it. The corresponding Bratteli diagram has at every level a single vertex and two edges, which gives an uncountable space $X_\cal B$, while $X_\sigma$ consists only of the single element $\ldots a a a \ldots$.

As a consequence, there are a number of issues where the two technologies, for Bratteli-Vershik systems on one hand, and for graph towers on the other, diverge:

\begin{enumerate}
\item
In their work \cite{BKM, BKMS, BKMS2} 
Bezuglyi, Kwiatkowski, Medynets and Solomyak  consistently  assume that the Vershik map has no periodic orbits. As a consequence, their main result for substitution subshifts (Theorem 2.9 of \cite{BKMS}) only holds for substitutions that don't have infinite periodic words $\ldots w w w \ldots$ in their associated subshifts. This excludes for instance the example from Remark \ref{ergodic-bound} (2) below.
 
 In comparison, for our ``parallel'' result, Theorem \ref{mieux-que-Moulinette}, no such restriction is needed.
\smallskip
\item
While for the graph tower approach (modified appropriately as described in in section \ref{dictionary})
{\em any} stationary graph tower defines a substitution, and conversely (see section \ref{stationary-pseudo}), 
in the above mentioned work of Bezuglyi-Kwiatkowski-Medynets-Solomyak the situation is more complicated: 
Stationary Bratteli diagrams typically do define a substitution subshift, but alternatively (this seems to correspond to the case where $V_\cal B$ has periodic orbits) such a stationary system can also describe what is called an ``odometer''.

\smallskip
\item
In the graph tower approach the passage from stationary towers to substitutions and back is an honest 1-1 relation, see section \ref{stationary-pseudo}. 
On the other hand, in \cite{BKMS}, pp. 992-996, it has been made explicit that starting with a stationary Bratteli diagram, reading off a substitution, and then deriving from that again a stationary Bratteli diagram, does in general not yield the original diagram, but typically a much larger (still stationary) diagram, with Vershik map conjugated to the original Vershik map. 
Correspondingly (see Remark \ref{Arnaud's-chouchou}), the matrix for which one has to calculate the eigenvectors in order to determine the invariant measures on the substitution subshift (compare Corollary 5.6 of \cite{BKMS} and Theorem \ref{mieux-que-Moulinette}) may in general be quite a bit larger than the original substitution incidence matrix.

\smallskip
\item
The formula (3.1) in \cite{BKMS} has bit of a formal ressemblance with our formula (\ref{9-7-2}) which is used to calculate the measure of any given cylinder (see also \S5 of \cite{BHL2} for an improvement in the stationary case). However, since most words in the subshift $X_\sigma$ are not represented by cylinders of the above type $[\underbar{e}]$ in the Bratteli diagram, a true analogue to (\ref{9-7-2}) in the Bratteli-Vershik setting could be tricky.

\smallskip
\item
The efforts to extend the Bratteli-Vershik technology from \cite{BKMS} beyond stationary systems (see \cite{BKMS2}) seem to be restricted to finite rank systems (i.e. with uniform bound on the number of vertices at every level). For graph towers, however, the analogous case of towers with finite tower dimension (see subsection \ref{thin-thick-towers}) certainly deserves special attention, but there is no reason whatsoever not to employ the very same techniques also to investigate measures on infinite dimensional graph towers.

\smallskip
\item
In \cite{BKM, BKMS, BKMS2} the authors consider also measures that may take on an infinite value on some cylinders. Since in the context of currents in free groups (see section \ref{Outer-space}) such measures are not permitted, we have decided to keep things simple here and to refrain from such an extension into uncharted territory. In principle, however, there doesn't seem to be any obstruction against admitting in sections \ref{weights-currents} - \ref{sec:weight-vectors} 
also weight functions or vector functions which on some edges take on the value $\infty$.
\end{enumerate}
\end{rem}

\subsection{Convex cones in a vector space}
\label{cones}

${}^{}$

Let $V$ denote a (possibly infinite dimensional) vector space over $\R$.  
For any subset $V_0 \subset V$ we define the {\em dimension} $\dim V_0$ to be the dimension of the 
subspace generated by $V_0$.

Given a family $F$ of vectors $\vec v_i \in V$, a vector $\vec v \in V$ is a {\em convex 
combination} from $F$ if for any $\vec v_i$ of some finite subfamily $F_0$ of $F$ there exist $\lambda_i \geq 0$ such that:
$$\vec v = \sum_{v_i \in F_0} \lambda_i \vec v_i$$
We denote by $\cal C(F) \subset V$ the {\em convex hull} of $F$, i.e. 
the set of all convex combinations from $F$. A subset $\cal C \subset V$ is called a {\em convex linear cone} if it is stable under taking convex combinations:
For any two $\vec{v_1}, \vec{v_2} \in \cal C$ and any two $\lambda_1, \lambda_2 \in \R_{\geq 0}$ one has $\lambda_1 \vec{v_1} + \lambda_2 \vec{v_2} \in \cal C$. 
In particular we always have $\vec 0 \in \cal C$ for any convex linear cone $\cal C$. By the {\em projectivization} $\Pr\cal C$ of any convex linear cone $\cal C \subset V$ we always mean the image of $\cal C \smallsetminus \{\vec 0\}$ in the projective space $\Pr V$.

A convex linear cone $\cal C \subset V$ is {\em closed} if it is a closed subset of $V$.
We observe that for any subset $V_0 \subset V$ the convex hull $\cal C(V_0)$ is a 
convex linear cone in $V$ (sometimes called {\em the cone spanned by $V_0$}). 
If $V_0$ is finite, then $\cal C(V_0)$ is closed.

A convex linear cone $\cal C$ is called {\em salient} if the following equality is satiesfied:
\begin{equation}
\label{salient}
\cal C \cap -\, \cal C = \{\vec 0\}
\end{equation}
A point $\vec v$ in a convex linear cone $\cal C$ is called {\em extremal} if it is not a convex combination from 
$\cal C \smallsetminus\{\lambda \vec v \mid \lambda \in \R_{\geq 0}\}$.
Denoting by $\cal C_{ext} \subset \cal C$ the subset of extremal points of $\cal C$, 
the finite dimensional version of the theorem of Krein-Milman, applied to $\Pr C$, gives:

\begin{lem}
\label{obvious}
Any finite dimensional salient closed convex linear 
cone $\cal C \neq \{\vec 0\}$ is the convex hull of its extremal points:
$\cal C = \cal C(\cal C_{ext})$
\qed
\end{lem}

A finite dimensional salient closed convex linear cone $\cal C \subset V$ with up to scalar multiples only finitely many extremal points 
has as {projectivization} $\Pr\cal C$ a finite convex polyhedron. A special case is given if the extremal points of $\Pr\cal C$ are {\em projectively independent}, i.e. they are represented by a family of linearly independent vectors in $\cal C$. This is captured by the following:

\begin{defn}
\label{simplicial}
A finite dimensional salient closed convex linear cone 
$\cal C$ is {\em simplicial} if the projectivization $\Pr\cal C \subset \Pr V$ of $\cal C \smallsetminus \{\vec 0\}$ is a simplex (of dimension $\dim \cal C -1$). 
\end{defn}

An important example of a salient closed convex linear cone is given by the set of invariant mesures on any dynamical system, see 
\cite{Walters}, \S 6. In our context, this gives:

\begin{rem}
\label{cone-of-measures}
Let $\Lsig \subset \Sigma(\Gamma)$ be a symbolic lamination on a graph $\Gamma$. Then we have:
\begin{enumerate}
\item
The set $\cal M := \cal M(\Lsig)$ of shift-invariant finite measures on $\Lsig $ is a 
salient closed convex linear cone different from $\{\vec 0\}$. Here one can take as ambient vector space the space $V = \R^{\cal P(\Gamma)}$, for $P(\Gamma)$ as in Definition \ref{language+} (1).
\item
The set $\cal M_{ext} \subset \cal M$ of its extremal points consists precisely of the ergodic measures on $\Lsig $. 
\item
The set $\cal M_{ext}$ maps under projectivization to a set $\Pr \cal M_{ext} \subset \Pr\cal M$ of points that are projectively independent. In particular, if $\cal M$ is finite dimensional, then  it is simplicial.
\item
For arbitrary subshifts $\Lsig$ the 
(possibly infinite dimensional) projectivization $\Pr \cal M(\Lsig)$ is compact, and hence (by Krein-Milman) it is the closure of the convex hull of its extremal points. Hence every measure $\mu_\Sigma \in \cal M$ is a (non-unique) infinite sum of ergodic measures. 
It turns out that the ergodic measures are dense in $\cal M(\Sigma(\Gamma))$.
\end{enumerate}
\end{rem}

We say that a map $f: \cal C_1 \to \cal C_2$ between convex linear cones $\cal C_i$ is {\em linear} if it is induced by some linear map between ambient vector spaces for the $\cal C_i$. A bijective linear map between 
convex linear cones is called a {\em cone isomorphism}. Two 
convex linear cones $\cal C_1$ and $\cal C_2$ are {\em isomorphic} if there exists a cone isomorphism $f: \cal \cal C_1 \to C_2$.
For finite dimensional convex linear cones $\cal C_i$ this is equivalent to demanding that
there exist surjective linear maps $f: \cal C_1 \to \cal C_2$ and $f': \cal C_2 \to \cal C_1$. Isomorphic 
convex linear cones have in particular the same dimension, even if this is not true for their ambient vector spaces.

The following observation follows directly from Lemma \ref{obvious}: 

\begin{lem}
\label{extremal-preimage}
Let $\cal C_1$ and $\cal C_2$ be two finite dimensional salient closed convex linear cones, and let $f: \cal C_1 \to \cal C_2$ be a linear map. If $f$ is surjective, then any extremal point of $\cal C_2$ has an $f$-preimage point in $\cal C_1$ which is also extremal.
\qed
\end{lem}

\medskip
\subsection{Eigenvectors of non-negative matrices}
\label{non-negative-matrices}

${}^{}$

Everything in this appendix is known (or even well-known) 
and can be found in the literature. We indicate some sources, but sparcingly, as non-negative matrices occur too frequently in too many diverse parts of mathematics for us to be able to do justice to all parties involved. In fact, parallel developments have sometimes lead to conflicting terminologies; we chose here the one which is most convenient for the purposes of this article.

\smallskip

A non-negative $(d \times d)$-matrix $M$ is called {\em reducible} if there exist an $M$-invariant coordinate subspace $\R^k \subset \R^d$ with $1 \leq k \leq d-1$. Here ``coordinate subspace'' means that $\R^k$ is generated by a subset of the standard basis $\vec e_1, \ldots, \vec e_d$ for $\R^d$.

It follows in particular that any $(d \times d)$-zero-matrix $M = 0$ is reducible, as long as $d \geq 2$. 

If $M$ is not reducible it is called {\em irreducible}. This includes formally the special case where $M$ is the $(1 \times 1)$-zero matrix $[0]$.
(Watch out: some authors like \cite{Seneta} exclude $[0]$ from what they call ``irreducible'', while others like \cite{Tam-Schneider-94} don't.)

If $M \neq [0]$ and furthermore any positive power $M^t$ is irreducible, then the matrix $M$ is called {\em primitive}. 
A  primitive matrix $M$ is characterized by the property that it possesses a power $M^t$ with $t \geq 1$ which is {\em positive}, i.e. any of the coefficients of $M^t$ satisfies $m^{(t)}_{i, j} > 0$, see \cite{Seneta}.

\begin{defn-rem}
\label{block-decomposition}
It is well-known (and easy to prove) that for any non-negative $(d \times d)$-matrix $M_0$ 
there is a conjugation with a permutation matrix such that the resulting matrix $M$ 
is in {\em Frobenius form}: $M$ admits a decomposition into matrix blocks, i.e. submatrices $M_{i, j}$, which have the following properties:
\begin{enumerate}
\item
For every upper diagonal block (i.e. $i < j$) one has $M_{i, j} = 0$.
\item
Every diagonal block $M_{i, i}$ is an irreducible 
(possibly $(1\times 1)$-zero) square matrix.
\end{enumerate}
We denote by $B(M) = \{B_1, \ldots, B_s\}$ the partition of the standard basis $\vec e_1, \ldots, \vec e_d$ 
into {\em coordinate blocks} $B_i$, defined by the above block decomposition of $M$.
It is easy to derive from the above conditions (1) and (2) that the partition $B(M)$ 
as well as the corresponding matrix block decomposition of $M$ is uniquely determined by $M_0$, up to a possible permutation of some coordinate blocks. 

Since it is convenient, we will from now on tacitly assume that any non-negative square matrix has been conjugated with a permutation matrix so that it is in Frobenius form.
\end{defn-rem}

From the above definitions of ``irreducible'' and ``primitive'' it follows (see \cite{Seneta}, Theorem 1.4) that any irreducible $M \neq [0]$ 
possesses a positive power $M^t$ which is a block diagonal matrix with primitive diagonal blocks. This block decomposition is stable with respect to passing to further positive powers of $M^t$. We obtain (compare 
\cite{LU1}):

\begin{defn-rem}
\label{block-primitive}
Let $M$ be a non-negative $(d \times d)$-matrix.
\begin{enumerate}
\item
The matrix $M$ is said to be {\em in primitive Frobenius form} if it is in Frobenius form, and if every diagonal block in the associated matrix block decomposition is primitive or 
$(1 \times 1)$-zero.
\item
Any non-negative $M$ possesses a positive power $M^t$ which is in primitive Frobenius form. The canonical block decomposition $B(M^t)$ is a refinement of $B(M)$. It is stable with respect to passage to further positive powers $(M^t)^{t'}$.
\end{enumerate}
\end{defn-rem}

For any non-negative square matrix $M$ we consider again the canonical decomposition $B(M)$ of $\{\vec e_1, \ldots, \vec e_d\}$ into {coordinate blocks} $B_i$, and we set $B_j > B_i$ if $M_{i, j}$ is non-zero. We then denote by $\succeq$ the preorder among the $B_i$ generated by the relation ``$>$'', and observe, since $M$ is block lower triangular, that $\succeq$ is indeed a partial order.

We write 
$B_j \succ B_i$ if $B_j \succeq B_i$ and $B_j \neq B_i$. This includes the possibility that $B_j \succ B_i$ holds, but not $B_j > B_i$. However, from the above set-up it is not hard to derive  
(see Lemma 4.6 of \cite{LU1}\,\footnote{\,\, The ``only if'' part of Remark \ref{partial-order} has been slightly overstated in Lemma 4.6 of \cite{LU1}, but fortunately only the ``if'' part is ever used there.}) 
the following analogue to the irreducible case:

\begin{rem}
\label{partial-order}
If $M$ is in primitive Frobenius form, then there is a bound $t_0 \geq 1$ such that for any integer $t \geq t_0$ the positive power $M^t$ has the following property: For the canonical block decomposition $B(M) = B(M^t)$ any
of the matrix blocks $M^{(t)}_{i, j}$ of $M^t$ is either positive or zero, and the latter case occurs if and only if for any coordinate block $B_k$ with
$$B_j \succeq B_k \succeq B_i$$
the corresponding diagonal matrix block satisfies $M^{(t)}_{k,k} = [0]$.
\end{rem}

It is well-known that any non-negative irreducible $(d \times d)$-matrix $M$ (including possibly the $(1\times1)$-zero matrix)
possesses up to scalar multiples a unique non-negative eigenvector. This {\em Perron-Frobenius} eigenvector $\vec v_{PF}(M)$, assumed 
here to be normalized so that the sum of its coefficients is equal to $1$, has {\em Perron-Frobenius} eigenvalue $\lambda_{PF}(M) \geq 0$, with $\lambda_{PF}(M) > 0$ if $M \neq [0]$. If $M$ is actually primitive, then for any non-negative vector $\vec v \in \R_{\geq 0}^d$ the sequence of vectors $M^t \vec v$ converges projectively to $\vec v_{PF}(M)$.

If $M$ is reducible,
the situation is of course more complicated: 
We consider again the canonical coordinate block decomposition $B(M)$ and define 
$B_j \in B(M)$ to be a {\em distinguished} coordinate block if $\lambda_{PF}(M_{j, j}) > \lambda_{PF}(M_{i, i})$ holds whenever one has $B_j \succ B_i$. In this case the corresponding diagonal matrix block $M_{j,j}$ of $M$ is also called {\em distinguished}.

\begin{rem}
\label{alert-error}
We'd like to emphasize the danger for confusion which comes from the two ``conflicting'' partial orders that are defined on the diagonal blocks  $M_{j,j}$ of a reducible matrix $M$ through $\lambda_{PF}(M_{j, j}) \geq \lambda_{PF}(M_{i, i})$ on one hand, and through $B_j \succeq B_i$ on the other. In particular, the ``top'' blocks in $M$ are in general not automatically distinguished\footnote{\,\, Unless we have overlooked an extra hypothesis, the named confusion seems to have occurred in \cite{BKMS}, p. 983 and Theorem 3.2 (a).} !
\end{rem}

The following is classic (essentially known already to Frobenius, see Theorem 2.1 of \cite{Tam-Schneider-94}).

\begin{prop}
\label{principal-ev}
For any non-negative square matrix $M$ in primitive Frobenius form and any distinguished block $B_j \in B(M)$ there exists a ``principal'' eigenvector $\vec v_j$ of $M$ which has the following properties:
\begin{enumerate}
\item
$\vec v_j$ is non-negative.
\item
$\vec v_j|_{B_j} = \vec v_{PF}(M_{j,j})$, where by the vector 
$\vec v_j|_{B_j}$ we mean the ``subvector'' obtained from $\vec v_j$ by only considering the coordinates of $B_j$.
\item
For any block $B_i \in B(M)$ the similarly defined vector $\vec v_j|_{B_i}$ for $B_i$ satisfies:
\begin{enumerate}
\item
 $\vec v_j|_{B_i}$ is positive if $B_j \succ B_i$, and 
 \item
 $\vec v_j|_{B_i} = \vec 0$ otherwise.
\end{enumerate}
\item
For any other eigenvector $\vec v$ of $M$ with the properties (1) and (3)(b), if $\vec v|_{B_j} \neq \vec 0$, then $\vec v$ is a scalar multiple of $\vec v_j$.
\end{enumerate}
Furthermore, any non-negative eigenvector of $M$ has eigenvalue $\lambda = \lambda_{PF}(M_{j, j})$ for some distinguished block $B_j$ of $M$, and $\vec v$ is a non-negative linear combination of the set of all principal eigenvectors $\vec v_i$ of $M$ with this eigenvalue $\lambda$.
\end{prop}

Since any primitive integer square matrix with Perron-Frobenius eigenvalue $> 1$, if it is invertible over $\Z$, requires a size of at least 2, one obtains as direct consequence:

\begin{cor}
\label{half-bound}
If a non-negative matrix $M \in {GL}_d(\Z)$ in primitive Frobenius form has no column of coefficient sum $\leq 1$, then the number of its principal eigenvectors is bounded above by~$\frac{d}{2}$.
\qed
\end{cor}

If a non-negative matrix $M \neq [0]$ 
is irreducible but not primitive, then $B(M)$ consists of a single block, but for some $k \geq 2$ and any positive multiple $t$ of $k$, 
the canonical block decomposition $B(M^t)$ for $M^t$ has more than one block $B_1, \ldots, B_k$, and $M^t$ is a block diagonal matrix where any of the $k$ diagonal blocks $M^{(t)}_{i, i}$ is primitive, and they are all conjugate to each other.
In particular, the Perron-Frobenius eigenvector $\vec v_{PF}(M)$ of $M$ is the barycentric linear combination of the principal eigenvectors 
$\vec v_i$ of $M^t$.

The very same reasoning is true for any reducible $M$ which is not in primitive Frobenius form. We consider some power $M^t$ of $M$ in primitive Frobenius form, and for any distinguished diagonal block $M_{j,j}$ of $M$ we consider those diagonal blocks $M^{(t)}_{i,i}$ of $M^t$ that derive from $M_{j,j}$ when the canonical block partition $B(M)$ is refined to $B(M^t)$. If $M_{j,j} \neq [0]$ 
is irreducible but not primitive, then there is a non-negative eigenvector $\vec v_j$ of $M$ associated to $M_{j,j}$, which satisfies $\vec v_j|_{B_j} = \vec v_{PF}(M_{j,j})$, and just as in the previous paragraph, the vector $\vec v_j$ is the barycentric linear combination of the principal eigenvectors $\vec v_i$ of $M^t$ associated to the $M^{(t)}_{i,i}$ (which are themselves distinguished diagonal blocks of $M^t$). 

\begin{defn-rem}
\label{barycentric-principal}
For any reducible matrix $M$ which is not in primitive Frobenius form
we call any eigenvector $\vec v_j$ as above a {\em barycentric principal eigenvector} of $M$, since we want to reserve the terminology ``principal eigenvector'' for the primitive Frobenius case:

Indeed, by a {\em principal eigenvector $\vec v$ of $M$} we will denote 
a principal eigenvector of some power $M^t$ (with $t \geq 1$) of $M$ 
which is in primitive Frobenius form. 
It follows from this convention and the above facts that this definition is independent of $t$, that there are only finitely many such $\vec v$, and that the principal eigenvectors for $M$ agree with those for any positive power of $M$.
(Watch out: A ``principal eigenvector'' of $M$ may thus actually not be an eigenvector of $M$ but only of some positive power of $M$ !)
\end{defn-rem}

\begin{rem}
\label{verbatim-true}
It follows that Proposition \ref{principal-ev} is verbatim true if one
omits the hypothesis ``in primitive Frobenius form'' and 
simultaneously replaces ``principal eigenvectors'' by ``barycentric principal eigenvectors''.
\end{rem}

\medskip

To terminate this appendix, 
we want to consider the {\em non-negative limit cone} (also called the {\em core}) of any non-negative $(d \times d)$-matrix $M$. By this we mean the infinite intersection $\cal C_\infty(M)$ of the nested sequence of {\em iterated $M$-image cones} $\cal C_t(M) := M^t(\R_{\geq 0}^d)$ for all $t \in \N$.

Contrary to the above pointed out similarity between the primitive and the irreducible case regarding the statement of Proposition \ref{principal-ev},  in order to determine the limit cone $\cal C_\infty(M)$ it is unavoidable to first pass to some positive power $M^t$ of $M$ in primitive Frobenius  form. 
Indeed, simultaneously omitting in Proposition \ref{limit-cone} below the hypothesis ``in primitive Frobenius from'' and replacing ``principal eigenvectors'' by ``barycentric principal eigenvectors'', as done above in Remark \ref{verbatim-true}, would lead to a distinctly incorrect statement (as can be seen already in the irreducible non-primitive case).

The following result from (applied) linear algebra goes way back; Pullman \cite{Pullman}
considered already some special cases, and a proof scheme is given in \cite{Schneider-85}. A complete proof in a slightly more general context can be found in \cite{Tam-Schneider-94}, which can also serve as standard reference and introduction/survey to this area. An
elementary (but elaborate) proof of a slightly stronger result has recently been provided in 
\cite{LU1}.

\begin{prop}
\label{limit-cone}
Let $M$ be a non-negative $(d \times d)$-matrix, let $M'$ be any positive power of $M$ which is in primitive Frobenius form, and let $V = \{\vec v_1, \dots, \vec v_s\}$ the set of principal eigenvectors of $M'$ with non-zero eigenvalue.

Then the limit cone 
$$\cal C_\infty(M) = \underset{t \in \N}{\bigcap}\, M^t(\R_{\geq 0}^d)$$
is equal to the cone spanned by the set $V$:
$$\cal C_\infty(M) = \left\{\sum_{i = 1, \ldots, s} c_i \vec v_i 
\,\,\middle| \,\, c_i \geq 0 \,\,\, \text{for all} \,\,\, i = 1, \ldots, s \right\}$$
In particular, $\cal C_\infty(M)$ is a simplicial cone with extremal rays given precisely by $V$.
\end{prop}

It remains to clarify some discrepancies in the terminology used in disparate areas of the literature:

\begin{rem}
\label{clarification}
(1)
In applied linear algebra, in particular in the papers by H. Schneider and his students or collaborators, the terminology {\em distinguished eigenvector} is used for any eigenvector which is a non-negative linear combination of ``our'' principal eigenvectors. As a consequence, every non-negative eigenvector of a non-negative square matrix $M$ for them is distinguished.

\smallskip
\noindent
(2)
In dynamics, however, the term {\em distinguished eigenvector} seems to be reserved for scalar multiples of our principal eigenvectors of $M$, if $M$ is in primitive Frobenius form, or of the barycentric principal eigenvector of $M$ in the general case.

\end{rem}

\end{document}